\documentclass[a4paper, 12pt, leqno]{article}
\usepackage{amsmath,amsfonts}
\usepackage[T1]{fontenc}
\usepackage{ae,aecompl} 
\begin{document}
\bibliographystyle{plain}
\newtheorem{theo}{Theorem}[section]
\newtheorem{lemme}[theo]{Lemma}
\newtheorem{cor}[theo]{Corollary}
\newtheorem{defi}[theo]{Definition}
\newtheorem{prop}[theo]{Proposition}
\newtheorem{problem}[theo]{Problem}
\newcommand{\beq}{\begin{eqnarray}}
\newcommand{\enq}{\end{eqnarray}}
\newcommand{\be}{\begin{eqnarray*}}
\newcommand{\en}{\end{eqnarray*}}
\newcommand{\Td}{\mathbb T^d}
\newcommand{\Rd}{\mathbb R^d}
\newcommand{\R}{\mathbb R}
\newcommand{\N}{\mathbb N}
\newcommand{\Zd}{\mathbb Z^d}
\newcommand{\Linf}{L^{\infty}}
\newcommand{\dt}{\partial_t}
\newcommand{\Dt}{\frac{d}{dt}}
\newcommand{\Dtt}{\frac{d^2}{dt^2}}
\newcommand{\demi}{\frac{1}{2}}
\newcommand{\vf}{\varphi}
\newcommand{\epu}{_{\epsilon}}
\newcommand{\ep}{^{\epsilon}}
\newcommand{\bfi}{{\mathbf \Phi}}
\newcommand{\bpsi}{{\mathbf \Psi}}
\newcommand{\bx}{{\mathbf x}}
\newcommand{\ds}{\displaystyle}

\let\cal=\mathcal
\title{The reconstruction problem \\ for the Euler-Poisson system in cosmology}
\author{Gr\'egoire LOEPER\footnotemark[1]}
\maketitle
\footnotetext[1]{Laboratoire J.A.Dieudonn\'e, Universit\'e de
Nice-Sophia-Antipolis \& 
Ecole polytechnique f\'ed\'erale de Lausanne}

\begin{abstract}
The motion of a continuum of matter subject to gravitational interaction
 is classically described by the Euler-Poisson system.  
Prescribing the density of matter at initial and final times, we are able to  obtain weak solutions for this equation by minimizing the action of the Lagrangian which is a convex functional. 
Through this variational formulation, the reconstruction problem becomes very similar to an optimal transportation problem.
Then we see that such minimizing solutions are consistent with smooth  solutions of the Euler-Poisson system and enjoy some special                 
 regularity properties.
\end{abstract}

\section{Introduction}
The Euler-Poisson system describes the motion of a self-gravitating fluid. It is used in cosmology, to model the evolution of the primitive universe. 
In the classical (non-relativistic) description, the gravitational field generated by a 
continuum of matter with density $\bar\rho$, is the gradient of a potential $\bar p$ satisfying the Poisson equation. The system is thus the following
\be
&&\dt\bar\rho+\nabla\cdot(\bar\rho\bar v)=0,\\
&&\dt (\bar\rho\bar v)+ \nabla\cdot (\bar\rho \bar v\otimes\bar v) = - \bar\rho\nabla\bar p,\\
&&\Delta \bar p = \bar \rho,
\en
where $\bar v$ is the velocity field.
The model is expected to be valid at scales where no collisional effects enter into account, and no multi-streaming occurs, otherwise the Vlasov-Poisson system would offer a preferable description. 

Assuming that the Universe is described by this set of equations, we are now interested in solving the so-called reconstruction problem, that we introduce here. We can infer, from red-shift catalogs, the repartition of matter in the present Universe. 
We also know that at times very close to the Big Bang, the matter was highly concentrated, but with very small relative fluctuations of density.
From this knowledge, can we reconstruct the intermediate states of the Universe, as well as the initial and present velocities ? 
The problem amounts to reconstruct a solution of the above Euler-Poisson system based on the knowledge of the initial and final density fields, and seems ill-posed. Another condition, known as {\it slaving} reduces the number of unknowns, and renders the problem well-posed: we impose that the reconstructed solution has a potential velocity field, {\it i.e.} $\bar v=\nabla\bar \phi$ for some potential $\bar \phi$. 
This paper proposes a way to solve the reconstruction problem, by transforming it into a minimization problem, in many ways similar to an optimal transportation problem. Most of the results will be obtained by use of the Monge-Kantorovitch duality, a tool widely used in optimal transportation problems.

\subsection{Equations of motion in co-moving coordinates}
The system of equations we will look at, is not exactly the one displayed above.
We will first express the equations of motion taking into account the global expanding motion of the Universe. We consider the expansion factor $a(t)$, a scalar function of time. Following the global expansion, particles positions $\bx(t)$ are given by $\ds\bx(t)=a(t)\bx_0$ 
and $\bx_0$ is called the co-moving coordinate.
The local velocity corresponding to this uniform expansion is then 
$\ds v_u(t,x)=(\dot a/a) x.$
The mean value of the density $\rho_m$ is supposed to be close to the critical density, defined as the highest value of $\rho_m$ allowing an infinite expansion; this condition is known to imply that $a(t)=\left(t/t_0\right)^{2/3}$. Performing then the change of variables
\be
\left\{\begin{array}{lll}
\bar \rho = \rho \bar \rho_m, & \bar v = ({\dot a}/{a}){\bf x} + a(t) v,\\
\bar p = \bar p_m + p, & \tau = (t/t_0)^{2/3},
\end{array}\right.
\en
we obtain the new Euler-Poisson system
\beq
&&\partial_{\tau}\rho + \nabla\cdot (\rho v)=0,\nonumber \\
&& \partial_{\tau} v + v\cdot \nabla v = -(3/2\tau) ( v+ \nabla p),\label{epcosmo}  \\
&& \Delta p =(\rho-1)/\tau\nonumber.
\enq
Note that the homogeneity condition at time 0 reads $\rho(0)=1$.
The slaving constraint appears as a necessary condition for the right hand side of the momentum equation not to be singular as time goes to 0. 
In this work, we will consider a simplified version of this system where we
do not include the time dependence in the coupling between $\rho$ and $p$, as well as the drag term $-\frac{3}{2\tau} v$ in the second equation.
We will also restrict ourselves to the flat torus $\Td=\Rd/\Zd$ with $\rho_m=1$.
In this framework the Euler-Poisson system hereafter referred to as $(E-P)$ takes the following form:
\beq
&&\dt\rho+\nabla\cdot(\rho v)=0\label{4continuite},\\
&&\dt (\rho v)+ \nabla\cdot (\rho v\otimes v) = - \rho\nabla p\label{4euler-poisson},\\
&&\Delta p=\rho-1\label{4poisson},
\enq
with the additional constraint
\be
\int_{\Td}\rho(.,x)dx\equiv 1.\label{4contrainte}
\en
The only modification from the first set of equations is the neutralizing background effect, that transforms the Poisson equation $\Delta p =\rho$ into $\Delta p = \rho - 1$.
One of the main consequences is that in this model, $\rho=1, v=0$ is a solution. This corresponds in the
physical coordinates to a uniform expansion.
As explained in paragraph \ref{marcheaussi} below, the results and techniques used for this simplified system
extend naturally to the full  system (\ref{epcosmo}).
Note that this form is the cosmological one, {\it i.e.} that the potential is attractive. In the case
of a repulsive potential (used for the description of a plasma) the
associated Poisson equation would be $\Delta p = -[\rho-\rho_m].$

\subsection{Definition of the reconstruction problem}
Given the $(E-P)$ system, one can try to solve the Cauchy problem, {\it i.e.} given $\rho$ and $v$ 
at time $t=0$ find a solution  to (\ref{4continuite}, \ref{4euler-poisson}, \ref{4poisson}) 
on a time interval $[0,T[$. Another approach is to look for  a solution over the time interval
 $[0,T]$ satisfying the two conditions:
\beq
&&\rho|_{t=0}=\rho_0,\label{4debut}\\
&&\rho|_{t=T}=\rho_T.\label{4fin}
\enq 
This approach has been used by Brenier in \cite{Br4}, \cite{Br5} for the incompressible Euler 
equation 
and allows to introduce variational techniques.
Indeed the system  (\ref{4continuite}, \ref{4euler-poisson},  \ref{4poisson}) is hamiltonian,
with Hamiltonian (or energy) given by:
\be
H(\rho,v)=\demi\int_{\Td} |v(t,x)|^2d\rho(t,x) - |\nabla p(t,x)|^2 dx.
\en
Solutions of hamiltonian systems
are critical points for the action of the Lagrangian, here defined by
\beq\label{4defI}
I(\rho,v,p)=\demi\int_{0}^{T}\int_{\Td}  |v(t,x)|^2d\rho(t,x) + |\nabla p(t,x)|^2 dx dt,
\enq
under the constraints (\ref{4continuite}, \ref{4poisson}, \ref{4debut}, \ref{4fin}).
Expressed in terms of $(\rho, J=\rho v)$ the action becomes
\be
F(\rho,J,p)=\demi\int_{0}^{T}\int_{\Td} \left|\frac{J}{\rho}\right|^2d\rho(t,x) + |\nabla p(t,x)|^2 dx dt.
\en
We will see that this function is convex in $(\rho, J,p)$, therefore the critical point will necessarily be a minimizer. The goal of this paper will thus  be to solve the following problem:
\bigskip
\noindent
\begin{problem}\label{4infI1}
Find $\bar\rho,\bar v,\bar p$ such that
\be\label{4goal}
I(\bar\rho, \bar v,\bar p)=\inf I(\rho,v,p)
\en
over all the triple $(\rho,v,p)$ satisfying
\be
&&\dt \rho +\nabla\cdot (\rho v) =0,\\
&&\Delta p = \rho-1,\\
&&\rho|_{t=0}=\rho_0,\\
&&\rho|_{t=T}=\rho_T.
\en
\end{problem}
The problem is here formulated in a very vague way: we do not mention in what space lie $\rho$ and $v$ when we perform the minimization. This will be made precise in the next subsection. 

\subsection{Motivations}

The  interest of studying this boundary problem is twofold.
The natural motivation is its direct application in cosmology, for the reconstruction of the early Universe.
A first approach had been made in \cite{BF}, where the authors solved a variant of the problem, assuming the Zel'dovich approximation. This approximation turned the reconstruction problem into an optimal transportation problem with quadratic cost.
The dual Monge-Kantorovitch problem was solved numerically using an algorithm due to Bertseakas.  
A more detailed discussion about the physical aspects of the problem was also presented in that paper (see also \cite{F}).
 The reader can also refer to the PhD thesis of J.Bec (\cite{Bec}, p.11,12).

On the other hand, from a mathematical point of view, 
it has been observed first by Brenier in \cite{Br5} in the case of the Euler incompressible
equation and by 
Evans and Gomes in \cite{EG} in the smooth finite dimensional case that solutions of
Hamiltonian flows that minimize the action of the Lagrangian 
present additional interesting regularity features. 
Some aspects of the present work can be seen as a continuation of their contribution in one special case of infinite dimensional Hamiltonian system. 
In a very vague way, one of the common features of between those three approaches is the following: the dual formulation of the problem lets appear a Hamilton-Jacobi equation, coupled with a transport equation. Solutions of Hamilton-Jacobi equations are not expected generally to present much regularity (usually they are not differentiable), due to the formation in finite time of caustics. However, thanks to our variational approach, the solutions found here will behave in a nicer way than one can expect. On a sets that contains the dynamics ({\it i.e.} that contains the support of $\rho$) the solution of the Hamilton-Jacobi equation will be differentiable, and even $C^{1, L\log L}$ (see point 7 of Theorem \ref{4main2}). Moreover, the dynamics will be reversible on this set. Therefore, this study shows some links between the weak KAM theory developed by Fathi 
(\cite{FaSi}), Evans (previous reference and \cite{Ev}), optimal transportation and Monge-Kantorovitch duality.

Another interest of this study is to generalize the ''Lagrangian minimizing'' approach developed by Benamou and Brenier in \cite{BB1}, in which the authors gave a continuum mechanics interpretation of the Monge-Kantorovitch problem involving the concept of interpolation between two measures, induced by the Wasserstein distance. This interpolation was also introduced earlier by McCann in \cite{Mc1} to develop the useful concept of displacement convexity, that will also appear in this case.
Otto in \cite{O} also used it to endow the set of probability measures with a formal Riemannian metric, in which the interpolation plays the role of geodesics, then allowing rich interpretations of some dissipative equations in terms of gradient flows.
Here our variational problem induces an interpolation that has somehow more regularity than the one of \cite{BB1} where the Lagrangian is only
\be
\demi\int_{0}^{T}\int_{\Td} |v(t,x)|^2 \rho(t,x) .
\en
Indeed we will see that the additional Dirichlet term forces the intermediate densities to be in $\Linf(\Td)$ independently of the initial and final densities. 
One can give the following heuristic interpretation of this result: the gravitational force being attractive, 
the system has to go first through an expansion so that a concentration does not appear before the final time $T$; this fact is expressed through the  differential inequality (\ref{4ineq2}). Meanwhile, some interesting displacement convexity properties, similar to the one found in \cite{Mc1},  will also appear.

The techniques we will present here can be naturally adapted to minimize the functionals
\be
\int_0^T\int_{\Td} \demi  |v(t,x)|^2 d\rho(t,x) + \int_0^T\int_{\Td} {\cal F}(\rho)dt dx
\en
provided ${\cal F}$ is a convex functional in $\rho$. Examples are
${\cal F}(\rho) = \rho^{\gamma}$, yielding solutions of the gas dynamics system, with attractive pressure term. This system
has an application in the theory of large deviations for random matrices, see {\cite{Gui}}.

Another example is the simple case ${\cal F}(\rho) = \rho V$, for some potential $V(t,x)$. This yields solutions
of $\ds \dt(\rho v) + \nabla\cdot(\rho v\otimes v) = - \rho \nabla V$.

We finally mention that the first steps of this study had been already done in E. Camalet's PhD \cite{C} under supervision of Y.Brenier. 

\subsection{Organization of the paper}
The paper is organized as follows: we first give a slightly different form to the minimization problem (Problem \ref{4infI2}), which we show is equivalent to the formulation of Problem \ref{4infI1}. We then state our results: existence, uniqueness of the minimizer, optimality equations (Theorem \ref{4main}), and regularity properties of the optimal path (Theorem \ref{4main2}).

The rest of the paper is devoted to the proof of those results: in section 3 we show existence and uniqueness, in section 4 we derive the important formula (\ref{4complete}), this will yield that the  optimal path solves the Euler-Poisson system (section 5), and some partial regularity results (section 6). In section 7, we then investigate the consistency of our minimizing solution with other solutions of the Euler-Poisson system.
This ends the proof of Theorem \ref{4main}.
The  section 8 is devoted to the proof of Theorem \ref{4main2}: for this we introduce a time discretization of the problem (Problem \ref{4infIN}). For this approximate problem we are able to obtain rigorous estimates, that yield the result when the time step goes to 0.

Formal proofs are given that present all the arguments necessary for the main results: in paragraph 4.1 for Theorem \ref{4main}, and in paragraph 8.1 for Theorem \ref{4main2}.

\subsection{Precise definition of Problem \ref{4infI1}}
We introduce the domain $D=[0,T]\times\Td$. We also define  the flux of matter  $J$ by 
$J=\rho \,v$. 
Given $J\in \Rd, \rho \in \R^+$ we use the fact that 
\beq
\sup_{c\in \R, m\in \Rd, c+|m|^2/2\leq 0}\{\rho c + J\cdot m\}=\left\{\begin{array}{lll}+\infty \textrm{ if } \rho=0, J\neq 0, \\ 
0 \textrm{ if } J=0,\\
J^2/2\rho \textrm{ if }  \rho >0.\end{array}\right.\label{4J2sur2rho}
\enq
Notice that as a supremum of affine functions, this is a 
(possibly infinite) convex functional in
$(\rho,J)$.
Given $\rho_0,\rho_T$ as in Theorem \ref{4main}, 
the functional $I$ can thus be formulated as
\beq\label{4defItilde}
&&\tilde I(\rho,J,p)=
\sup_{c+|m|^2/2\leq 0}\left\{\int_{D}c(t,x) d\rho(t,x)  +  m(t,x)\cdot d J(t,x) \right\}+\demi\int_{D}|\nabla p(t,x)|^2\,dt dx,
\enq 
where the supremum is taken over all $(c,m) \in C(D)\times \left(C(D)\right)^d$. 
As we shall see in the next proposition, this formulation is consistent with the formulation (\ref{4defI}) in the case where $v \in L^2(D, d\rho)$ and well defined (although leading to possibly infinite value) for $\rho \in C([0,T];{\cal P}(\Td)-w*)$, $J\in ({\cal M}(D))^d$, $\nabla p \in L^2(D)$ where ${\cal M}(D)$ denotes the set of bounded
measures on $D$ and  ${\cal P}(\Td)$ the set of probability measures on $\Td$.

\begin{prop}\label{4coherence}
Let $\rho$ be a measure $[0,T]\times\Td$, $v$ be a $d\rho$ measurable vector field, and $p$ be a measurable function on $[0,T]\times \Td$.  
Let the functionals $I$ and $\tilde I$ be defined by (\ref{4defI}) and (\ref{4defItilde}) respectively, with the convention that $I=+\infty$ whenever $v\notin L^2(D,d\rho)$ or $\nabla p\notin L^2(D,dx)$.
Then \linebreak $I(\rho, v, p)= \tilde I(\rho,v,p)$ as functionals valued in $]-\infty, +\infty]$.
\end{prop}

{\it Proof.} 
First suppose that $J$ is not absolutely continuous with respect to $\rho$. Then there exists a set $S$ such that $|J|(S)>0$ and $|\rho|(S)=0$.
Then one constructs a sequence $c\epu, m\epu$ such that 
\be
\int c\epu d\rho + m\epu\cdot d J \to  +\infty.
\en
Indeed there exists a sequence ${\cal O}\epu$ of open sets containing $S$ such that
$|\rho|({\cal O}\epu) \leq \epsilon$ and $|J|({\cal O}\epu)  \geq \delta$ with $\delta>0$ fixed.
Then for each $\epsilon$, there exists $f\epu\in \left(C^{\infty}_c({\cal O}\epu)\right)^d, \|f\epu\|_{\Linf}\leq 1$ such that $\int f\epu\cdot dJ \geq \delta/2$. We then have, taking $m\epu=\lambda f\epu, c\epu=- |m\epu|^2/2$, 
\be
\int -\demi|\lambda f\epu |^2 d\rho + \lambda f\epu\cdot d J \geq \lambda\delta/2-\lambda^2  \epsilon/2 
\en 
which is greater than $\delta \epsilon^{-1/2}/2 - 1/2$ for $\lambda = \epsilon^{-1/2}$. Letting $\epsilon$ go to 0, we conclude that definitions (\ref{4defI}) and (\ref{4defItilde}) yield $+\infty$.

If $J<<\rho$, $J=\rho v$ for some $\rho$ measurable function $v$; we consider
$\ds v_P=\frac{v}{|v|}\max\{|v|, P\}$. By standard smoothing arguments, 
there exists a smooth sequence $(c_n=-|m_n|^2/2, m_n)$ such that
\be
\int c_n d\rho+ m_n\cdot v_P d\rho &=& \int d\rho|v_P|^2/2-\int d\rho|m_n-v_P|^2/2\\ 
&\to& \int d\rho|v_P|^2/2,
\en 
hence we can  construct a sequence $(c_P, m_P)_{P\in \N}$ such that $\ds\int d\rho\frac{|v_P-m_P|^2}{2}\to 0$, and this implies 
\be
\int c_Pd\rho+ m_P\cdot v d\rho 
&=&\int d\rho |v_P|^2/2 - \int d\rho |v_p-m_P|^2/2 + \int d\rho(v-v_P)\cdot v_p\\
&\geq& \int d\rho v\cdot v_p - d\rho |v_P|^2/2 - \epsilon(P)\\
&\geq& \int d\rho |v_P|^2/2 -\epsilon(P), 
\en
with $\epsilon(P) \to 0$ as $P\to +\infty$.
Hence $\tilde I \geq I$.

Then, if $\int \rho|v|^2$ is bounded, for any $c\leq -|m|^2/2$,
\be
\int c d\rho+ m\cdot v d\rho \leq  \int d\rho\frac{|v|^2}{2} - d\rho\frac{|v-m|^2}{2}\leq \int d\rho\frac{|v|^2}{2}.
\en
This shows $\tilde I \leq I$ whenever $I< +\infty$, and ends the proof of Proposition \ref{4coherence}.

$\hfill \Box$
 
The new formulation of the Problem \ref{4infI1} is then:
\begin{problem}\label{4infI2}
Minimize 
\be
\tilde I(\rho,J,p)&=&\sup_{c+|m|^2/2\leq 0}\left\{\int_{D} c(t,x) d\rho(t,x) + m(t,x)\cdot d J(t,x) \right\} \\ 
&&+\demi\int_{D}|\nabla p(t,x)|^2\,dt dx
\en 
among all $(\rho,J,p)$ that satisfy
$\rho \in {\cal M}(D)\cap C([0,T];{\cal P}(\Td)-w*)$, $J\in ({\cal M}(D))^d$, $\nabla p \in L^2(D)$, and satisfy in the distribution sense
\beq
&&\dt\rho+\nabla\cdot J = 0\label{4c1},\\
&&\Delta p=\rho-1\label{4c2},\\
&&\rho(t=0)=\rho_0\label{4c3},\\
&&\rho(t=T)=\rho_T.\label{4c4}
\enq
We denote 
$$K=\inf_{\rho,J,p} \tilde I(\rho,J,p)$$
among all such $(\rho,J,p)$.
\end{problem}

\subsection{Variational problem for the full equations in co-moving coordinates}\label{marcheaussi}
As we have seen above, the full system for which we want to solve the reconstruction problem is the following
 \be
\left\{\begin{array}{lll}
\partial_{\tau}\rho + \nabla\cdot (\rho v)=0, \\
 \partial_{\tau} v + v\cdot \nabla v = -3/(2\tau) ( v+ \nabla p),      \\
 \Delta p =(\rho-1)/\tau.
\end{array}\right.
\en 
Solutions for this system can be sought as minimizers of the action
\be
I_a=\demi\int_{0}^{T}\int_{\Td} \tau^{3/2} \left( |v(\tau,x)|^2d\rho(\tau,x) + |\nabla p(\tau,x)|^2 dx d\tau \right),
\en 
under the constraints
\be
&&\partial_{\tau}\rho + \nabla\cdot (\rho v)=0, \\
&&\rho(\tau=0)=1, \rho(\tau=T)=\rho_T.
\en
The same techniques as the one that will be exposed in this paper adapt to this minimization problem.

\bigskip


\section{Results}

\subsection{Notation}
\begin{itemize}
\item[]The space dimension will be denoted $d$. 
\item[]
We shall hereafter use the notation $A:B = {\rm trace }\, (A^t B)$ for $A,B$ two $d\times d$ matrices.
\item[] 
For a function $\varphi: \Rd \to \R^m$ we denote by $D\varphi$ its first derivative matrix equal to $\left( \frac{\partial \varphi^i}{\partial x_j} \right)_{\scriptsize{\begin{array}{ll}1\leq i \leq m \\  1\leq j \leq d\end{array}}}$. For $\varphi: \Rd \to \R$, we denote by $D^2\varphi$ its Hessian matrix equal to 
$\left( \frac{\partial^2 \varphi}{\partial x_i \partial x_j}\right)_{1\leq i,j\leq d}$.
\item[] We recall that $D=[0,T]\times \Td$, and ${\cal M}(D)$ is the space of bounded measures on $D$. 
\end{itemize}

\subsection{Definition of weak solutions for $(E-P)$}
\begin{defi}\label{4defiweakep}
A triple $(\rho, v, p)$ is said to be a weak solution of $(E-P)$ if:

\begin{enumerate}

\item $\rho  \in L^2([0,T]; H^{-1}(\Td))\cap C([0,T]; {\cal P}(\Td)-w*)$,  
$v \in  L^2(D, d\rho)$, 

\item  for any
$\varphi=(\varphi^j)_{j\in[1..d]}\in \left(C^{\infty}_c(]0,T[\times\Td)\right)^d$ one has
\beq\label{4weakep}
\int_{[0,T]\times\Td}&&\dt\varphi \cdot v \  d\rho +   D\varphi : v \otimes v \  d\rho - \varphi \cdot \nabla p\\ 
&&+ D\varphi : \nabla p\otimes\nabla p - \frac{1}{2}(\nabla\cdot \varphi)  |\nabla p|^2  =0,\nonumber
\enq

\item for any
$\varphi\in C^{\infty}([0,T]\times\Td)$:
\be
&&\int_{[0,T]\times\Td}\dt \varphi d\rho +  \nabla\varphi\cdot v d\rho =\int_{\Td}\rho_T\varphi|_{t=T}-\int_{\Td}\rho_0\varphi|_{t=0}\label{4weakconti0},\\
&& \int_{[0,T]\times\Td}(d\rho-1)\varphi +\nabla p \cdot\nabla\varphi =0\label{4weakpoisson0}.
\en
\end{enumerate}
\end{defi} 
Equation (\ref{4weakep}) is equivalent to equation (\ref{4euler-poisson}) for smooth $p$ using the identity 
$$(1+\Delta p)\nabla p=\nabla\cdot (\nabla p \otimes \nabla p)-\demi\nabla|\nabla p|^2  + \nabla p.$$
The right hand side of this identity is well defined in the sense of distribution    if we only know that  $\nabla p \in L^2 (D)$.

\subsection{Statement of the Theorems}

We first have an existence/uniqueness result for the minimizer of the action $\tilde I$ defined in (\ref{4defItilde}). Note that from Proposition \ref{4coherence} this minimizer yields also the minimizer of the action $I$ defined by (\ref{4defI}).
\begin{theo}\label{4main}
Let $\rho_0, \rho_T $ be two  probability measures in $L^{\frac{2d}{d+2}}(\Td)$, then there exists a unique \linebreak 
$(\rho,  J,   p)\in ({\cal M}(D) \times \left({\cal M}(D)\right)^d \times L^2([0,T]; H^1(\Td)))$ 
with $\Delta p=\rho-1$ in ${\cal D}'$ minimizer of the  Problem \ref{4infI2}.
The flux $J$ has a density $v$ with respect to $\rho$, 
$(  \rho,   v,   p)$ is a weak solution of the Euler Poisson system $(E-P)$ in the sense of Definition \ref{4defiweakep}
and coincides with any smooth solution of $(E-P)$  
satisfying (\ref{4c3}, \ref{4c4}) and having a potential velocity; such solution is therefore unique.
Moreover  

\begin{enumerate}
\item there exists $\phi \in L^2_{loc}(]0,T[;H^1(\Td))\cap \Linf_{loc}(]0,T[ \times \Td) $ such that $v=\nabla \phi \ d\rho$  a.e.  and we can thus extend the definition of $v$  to all of $\Td$ as  a function belonging to $L^2_{loc}(]0,T[ \times \Td))$,

\item any such  extension satisfies 
\beq\label{4mainreg}
\int_{\Td}\int_{\tau}^{T-\tau} \left|  v(t,x+y)-  v(t,x)\right|^2 d\rho(t,x) \leq C_\tau|y|^2 
\enq 
for all $\tau$ in $]0,T/2]$, $y$ in $\Rd$,

\item the density $\rho$ belongs to $L^2_{loc}(]0,T[\times \Td)\cap C(]0,T[; L^{p})$ for any $p\in [1,3/2[$.
\end{enumerate}

\end{theo}

Then we have the regularity result:

\begin{theo}\label{4main2}
If $\rho_0$ and $\rho_T$ are in $L^{\frac{2d}{d+2}}$ the unique  solution $(\rho,J,p)$ of Problem \ref{4infI2} has the following regularity 
properties:

\begin{enumerate}
\item  The density $\rho$  is in 
$\Linf_{loc}(]0,T[\times\Td)\cap  C(]0,T[;L^k(\Td))$ for every $1\leq k <\infty$:
for every \linebreak $\tau \in ]0,T/2[$ there exists $C_{\tau}$ such that for every $t$ in $[\tau, T-\tau]$,
\be
\|\rho(t,\cdot)\|_{\Linf(\Td)}\leq C_{\tau},
\en 
 and there exists $C$ such that
\be
-C(1+\frac{1}{t})\leq\Dt\log \left(\|\rho(t,\cdot)\|_{L^k(\Td)}\right)\leq C(1+\frac{1}{T-t}),
\en 
moreover the constants $C_{\tau}, C$ are independent of the choice of $\rho_0$ and $\rho_T$.

\item The velocity $v=\nabla\phi$ can be chosen in $\Linf_{loc}(]0,T[\times\Td)$,
 this bound is also independent of the choice of $\rho_0$ and $\rho_T$.

\item The functions $\ds\int_{\Td}[\rho]^k(t,x)\,dx,\, k\geq 1,\, \int_{\Td}[\rho\log\rho](t,x)\,dx$
are convex with respect to time.

\item The velocity potential $\phi$ can be chosen  in  $W^{1,\infty}_{loc}(]0,T[\times\Td)$ 
and to be viscosity solution of \mbox{$\dt\phi +\demi|\nabla\phi|^2 + p=0$} on every $[s,t] \subset ]0,T[$.

\item If $\rho_T$ is in $L^p(\Td)$ with $p>d$ then  point 4  extends up to $t=T$.

\item  One can also choose $\phi$ such that $(\psi,q)(t,x)=(-\phi,p)(T-t,x)$ is a viscosity solution of
\mbox{$\dt\psi +\demi|\nabla\psi|^2 + q=0$} on every $[s,t] \subset ]0,T[$ and point 5 applies. Moreover both choices will co\"incide $d\rho$ almost everywhere.

\item For each $t\in ]0,T[$, there exists a closed set ${\cal S}_t$ of full measure for $\rho(t)$ such that $\phi(t)$ is differentiable with respect to the space variable at every point of ${\cal S}_t$. Moreover, for all $t\in [\tau, T-\tau]$, \linebreak $\tau>0$, for all $(x,y)$ in ${\cal S}_t$, with $|x-y|\leq 1/2$, 
\be
  \left| \nabla\phi(t,x)-\nabla\phi(t,y) \right| \leq C(\tau) |x-y| 
\log(\frac{1}{|x-y|}).
\en

\end{enumerate}

\end{theo} 
The reader can refer to the books of Evans \cite{E} and Barles \cite{Ba} for the definition of viscosity solution of $\dt\phi +\demi|\nabla\phi|^2 + p=0$.

\bigskip

{\it Remark 1.} The assumption that the final and initial densities are in $L^{\frac{2d}{d+2}}(\Td)$ 
is technical: it allows us to show that there exists at least one admissible path with finite action transporting $\rho_0$ on $\rho_T$, see section \ref{4ilexiste}. Actually all the results are true assuming that there exists a $(\rho, J,p)$ satisfying all the constraints (\ref{4c1},.., \ref{4c4}) such that $\tilde I(\rho,J,p)$ is finite.

{\it Remark 2.} The bound (\ref{4mainreg}) is a finite difference version of the formal (but non rigorous since $\rho$ has no regularity~)  assertion
$\int_{[\tau,T-\tau]\times{\Td}}d \rho|\nabla v|^2  <+\infty$. See  \cite{BoBu} where the authors look at an appropriate definition of the tangent space related to a measure.  

{\it Remark 3.} One may observe that the regularity obtained for $\rho$ is stronger in the second theorem than in the first. The third point of Theorem \ref{4main} is obtained just by using the regularity result of $v$ (point 2 of Theorem \ref{4main}), while in the second theorem we use
in a crucial way the gravitational coupling $\Delta p =\rho- 1$. The techniques are thus different, and the one employed in the first theorem could probably be used for other types of coupling.

{\it Remark 4.} The consistency with smooth solutions is detailed in Theorem \ref{4main3}. We try there to examine in what class our variational solution is the unique solution of the $(E-P)$ system satisfying the boundary conditions.  

{\it Remark 5.} In many assertions we only state that ``$\phi$ can be chosen in such a way that...''; this is because $\phi$ is uniquely determined only in the $d\rho$ a.e. sense. However, $\rho$ is uniquely defined, and $v$ is unique $d\rho$ a.e..

\section{Existence and uniqueness of a minimizer for the action}
This section is devoted to the proof of 
\begin{prop}\label{4prop1}
Under the assumption that $\rho_0$ and $\rho_T$ are in $L^{\frac{2d}{d+2}}$, there exists a unique minimizer $(\rho, J, \nabla p)$ in
$C([0,T];{\cal P}(\Td)-w*)\times \left({\cal M}(D)\right)^d \times L^2(D)$
for the Problem \ref{4infI2} under the constraints (\ref{4c1}, \ref{4c2}, \ref{4c3}, \ref{4c4}).
\end{prop} 

The proof of this result will be obtained in two ways: the first will be based and convexity / l.s.c. properties of the functional we are minimizing. However, despite the simplicity of this proof, we will introduce another way of obtaining the minimizing solution, by studying a dual problem similar to the  the Monge-Kantorovitch problem, classically used in optimal transportation problems. This technique will then be a useful tool to obtain regularity results. Therefore we will show the following 
\begin{prop}\label{4propdual}
Let $K=K(\rho_0, \rho_T)$ be the infimum of Problem \ref{4infI2}, and assume that $K<+\infty$, then
\be
K = \sup_{\scriptsize{\begin{array}{ll}\phi \in C^1(D), q\in C^1(D) \\ \dt\phi+q+|\nabla\phi|^2/2\leq 0\end{array}}} &&\left\{ \int_{\Td} \phi(T)\,d\rho_T  - \phi(0)\,d\rho_0 \right.
+ \left.\int_D q -|\nabla q|^2/2 \ dtdx\right\}. 
\en
\end{prop}

{\it Remark.} In the next paragraph, we are going to show that indeed $K<+\infty$, assuming some integrability conditions on the initial and final densities.

\subsubsection{Existence of an admissible solution}\label{4ilexiste}
We are now going to prove that the infimum of Problem $\ref{4infI2}$ is finite.
This proof will use the results obtained in \cite{BB1} and \cite{Mc1}, concerning the time continuous formulation of the Monge-Kantorovitch transport problem. The optimal solution for the optimal transport problem with quadratic cost turns out to have a finite action for our functional $\tilde I$, provided the initial and final densities are sufficiently integrable.
We will use the displacement convexity property obtained by McCann, in order to show that at intermediate times the density remains in $H^{-1}(\Td)$.

\begin{lemme}\label{4lemme-il-existe}
Under the assumption that $\rho_0, \rho_T$ are in $L^{\frac{2d}{d+2}}(\Td)$, there exists 
$(\rho, J=\rho v, p)$ satisfying (\ref{4c1}, \ref{4c2}, \ref{4c3}, \ref{4c4}) and such that $\tilde I(\rho, J, p)$ is finite. Moreover $\rho\in \Linf([0,T]; L^{\frac{2d}{d+2}}(\Td))$
and $v\in \Linf(D, d\rho)$.
\end{lemme}

{\it Proof.}  
We use the following result that combines \cite{BB1} and \cite{Mc1}:
\begin{prop}\label{mccann}
Let $\rho_0$ and $\rho_T$ belong to ${\cal P}(\Td)\cap L^k(\Td)$ for some $1 \leq k \leq  \infty$. 
There exists a unique pair $(\bar\rho(t,x), \bar J=\bar\rho \bar v(t,x))$ that minimizes the action
\be
A(\rho,\rho v)=\int_{[0,T]\times\Td}|v(t,x)|^2  \ d\rho(t,x)
\en
among all $(\rho,J)$ that satisfy
$\rho \in C([0,T];{\cal P}(\Td)-w*)$, $J\in ({\cal M}([0,T]\times\Td))^d$ and 
\be
&&\dt\rho+\nabla\cdot J = 0,\\
&&\rho(t=0)=\rho_0,\\
&&\rho(t=T)=\rho_T.
\en
$A(\bar\rho,\bar\rho\bar v)$ is finite, and for $k\geq 1-1/d$, the function
$t \rightarrow (k-1)\|\bar\rho(t,.)\|^k_{L^k}$ is convex for $t\in[0,T]$. Hence for $k\geq 1$, $\|\rho(t,\cdot)\|^k_{L^k}$  is bounded by 
$\max\{\|\rho_0\|^k_{L^k},\|\rho_T\|^k_{L^k}\}$. Finally $\bar v\in \Linf(D, 
d\bar\rho)$.
\end{prop}

The minimization problem stated here is one formulation of the classical optimal transportation problem with quadratic cost. It is described in more details in paragraph \ref{4optitrans}.

Using  classical elliptic regularity we have, for $1<k<\infty$,
$\|D^2 \Delta^{-1}[\rho(t)-1]\|_{L^k(\Td)} \leq C_k\|\rho(t)\|_{L^k(\Td)}$.
The Gagliardo-Nirenberg inequality gives
$\|\nabla p\|_{L^{k*}(\Td)} \leq C_k \|D^2 p \|_{L^k(\Td)}$ with $k*=\frac{dk}{d-k}$.
Therefore, if $k=\frac{2d}{d+2}$ we have 
$$\|\nabla \Delta^{-1}[\bar\rho(t,\cdot)-1]\|_{L^2(\Td)}\leq C\|\bar\rho(t,\cdot)\|_{L^k(\Td)}.$$ It follows that $(\bar\rho, \bar J=\bar\rho \bar v,\bar p)$ is an admissible pair for which $\tilde I$ is finite. This completes the proof of Lemma \ref{4lemme-il-existe}.

$\hfill \Box$

\subsection{Proof of Proposition \ref{4prop1}}

We present here a proof of existence and uniqueness of the minimizer of Problem \ref{4infI2}. As we have already seen, the functional
\be
\demi \int_D \left|\frac{J}{\rho}\right|^2 \ d\rho  + \demi \int_D |\nabla p|^2 \ dxdt
\en
is convex. The second term is strictly convex in $\rho$ (if we restrict to probability measures) and lower semi-continuous (l.s.c.) with respect to the weak-$*$ convergence of $\rho$.
To see that the first term is also l.s.c. with respect to the weak-$*$ convergence of the measures $(\rho, J)$, we refer to \cite{hut} where the following result is shown:

\begin{theo}\label{4hutchinson}
Let $(\mu_n, f_n)_{n\in \N}$ be a sequence such that
for all $n\in \N$, $\mu_n$ is a Radon measure on $D$, and $f_n$ is $d\mu_n$ measurable.
Assume that $\mu_n \rightharpoonup \mu$ for the weak-$*$ topology, with $\mu$ a Radon measure.
Assume  that $\int_{D} d\mu_n|f_n|^2$ is uniformly bounded. 
Then there exists a subsequence (still labeled by $n$) and a pair $(\bar \mu, \bar f)$, with $\bar f$ $d\bar \mu$ measurable, such that
\begin{enumerate}
\item for all $\varphi \in \left(C^0(D)\right)^d$,
 $\ds\int_{D} d\mu_n f_n \cdot \varphi \to \int_{D} d\bar\mu \bar f \cdot \varphi$,
\item $\ds\int_{D} d\bar \mu  |\bar f|^2 \leq \liminf \int d\mu_n |f_n|^2 $.
\end{enumerate}
\end{theo}

Considering a minimizing sequence $(\rho_n, J_n)$ for Problem \ref{4infI2}, we get also a sequence $(\rho_n, v_n)$ where $v_n$ is the density of $J_n$ with respect to $\rho_n$, that belongs to $L^2(D,d\rho_n)$. Thus $(\rho_n, v_n)$ will be as $(\mu_n, f_n)$ in the Theorem above.

We will show that the sequence $\rho_n$ is equicontinuous in $C([0,T]; {\cal P}(\Td)-w*)$.
For this we use the mass conservation equation  $\dt\rho_n + \nabla\cdot J_n = 0$. 
For any time $t$, we can estimate the total mass of $J_n$ by Cauchy-Schwartz inequality:
\be
{\mathbf M}(J_n(t)) \leq \left( \int_{\Td} d\rho_n(t,x) |v_n(t,x)|^2 \right)^{1/2},
\en  
hence $J_n\in L^2([0,T]; \left({\cal M}(\Td)\right)^d) \subset L^2([0,T]; H^{-s}(\Td))$ for $s$ large enough.
Therefore from the mass conservation equation, $\dt \rho_n$ is uniformly bounded in $L^2([0,T]; H^{-s'}(\Td))$ form some $s'$. 
Of course, $\rho_n$ is also bounded in $\Linf([0,T];{\cal M}(\Td))$. Using standard arguments of functional analysis (see\cite{Li}), we obtain that the sequence $\rho_n$ is equicontinuous in $C([0,T]; {\cal P}(\Td)-w*)$.

The first point shows that a subsequence $(\rho_n, J_n=\rho_nv_n)$ will converge weakly to $(\bar \rho, \bar J=\bar \rho \bar v)$. This implies that one can pass to the limit in the continuity equation $\dt \rho + \nabla\cdot J =0$.
Moreover the equicontinuity of the sequence $\rho_n$ yields that $\bar\rho(0)=\rho_0$, $\bar\rho(T)=\rho_T$.
We thus obtain that the pair $(\bar\rho, \bar J)$ satisfies (\ref{4c1}, \ref{4c3}, \ref{4c4}).
  
The second point shows that $\tilde I(\bar \rho, \bar J, \bar p)$ is smaller than $
\liminf \tilde I(\rho_n, J_n, p_n)$, and therefore equal to $K=\inf \tilde I$, since we consider a minimizing sequence.

For the uniqueness part, since the term $\int |\nabla p|^2$ is strictly convex with respect to $\rho$, two minimizing solutions must have the same density $\rho$.
Once the density is set, the velocities must co\"incide $d\rho$ a.e. again by convexity with respect to $v$ of $v\to \rho |v|^2$.
This proves Proposition \ref{4prop1}.

$\hfill\Box$

\subsection{Introduction of the dual problem}
We use here standard convex analysis arguments that can be found in \cite{B} and the proof is an adaptation of the one found in \cite{Br5}.
The constraints (\ref{4c1}, \ref{4c2}, \ref{4c3}, \ref{4c4}) can be formulated
in the following weak way:
\beq
&&\forall \phi \in C^{\infty}(D),\; \int_D \dt \phi (d\rho-d\bar\rho) + \nabla\phi\cdot (dJ-d\bar J)
 = 0\label{4newweakconstr1},\\
&& \forall q \in C^{\infty}(D),\; \int_D (d\rho-d\bar\rho) q 
 = - \int_D (\nabla p-\nabla\bar p) \cdot\nabla q \ dtdx,
\label{4newweakconstr2}
\enq 
where $(\bar \rho, \bar J, \bar p)$ satisfies the constraints (\ref{4c1}) to (\ref{4c4}), hence the triple of Lemma \ref{4ilexiste} works.
Minimizing $\tilde I$ under the constraints 
(\ref{4c1}) to (\ref{4c4}) is thus equivalent to find 
\be
K=\inf_{\rho,J,p}\sup_{\phi,q,c,m} &&\left\{\int_{D} d\rho\, c + dJ\cdot m - \dt \phi (d\rho-d\bar\rho) -\nabla\phi\cdot (dJ-d\bar J)\right.  \\
&& \left.+ \int \demi|\nabla p|^2 -\nabla q\cdot (\nabla p-\nabla\bar p)  \ dtdx - \int q(d\rho-d\bar\rho)\right\}, 
\en
with the supremum taken over all the continuous functions $c,m$ with $c:D \rightarrow \R$ and $m:D \rightarrow \Rd$ satisfying  $c+|m|^2/2\leq 0$.

\bigskip
\noindent
Here $C(D)$ is the space of continuous functions on $D$ and $C_{\#}(D)$ is defined by the additional constraint that  the integral over $\Td$ vanishes for all $t \in [0,T]$ . 
On $C(D)$ we have the usual duality bracket $<f,g>$ denoted by $\displaystyle\int_D f\,dg$ 
with $g\in {\cal M}(D)$ the set of bounded measures on $D$.
The dual space of $C_{\#}$ is reduced to the set of bounded measures $g$ on $D$ whose total mass at any time is zero 
({\it i.e.} for all $g\in C_{\#}'$, for all $z\in C^0[0,T], \  \int_D z(t)dg = 0$) 
and denoted by ${\cal M}_{\#}(D)$.
We introduce the functionals $\alpha$ and $\beta$ defined on $(c,m,r) \in C(D)\times (C(D))^d\times C_{\#}(D)$. It will be convenient to denote $r=\Delta q$, and this is possible since the mean value of $r$ is zero.
\be
&&\alpha(c,m,r)=\demi\int_D |\nabla \Delta^{-1} r|^2 \ dtdx = \demi\int_{D} |\nabla q|^2 \ dtdx\\
&&\textrm{ if }c+|m|^2/2\leq 0,\\
&&\alpha(c,m,r)=  +\infty \mbox{ otherwise};\\
&&\beta(c,m,r)=  \int_D d\bar\rho \, c + d\bar J \cdot m  + \int_D \bar p r  \ dtdx\\
&&\textrm{ if } \exists\; \phi\in C^1(D) \textrm{ such that } c+\dt\phi+q=0, \; m +\nabla_x\phi=0,\\
&&\beta(c,m,r)= +\infty \mbox{ otherwise}
\en
with $(\bar\rho,\bar J,\bar p)$ as above.

\bigskip
\noindent
We compute $\alpha^*$ and $\beta^*$ the Legendre-Fenchel transform (see \cite{B} for definition) 
 of respectively $\alpha$ and $\beta$. They are defined on $(\rho, J, p)\in {\cal M}(D)\times
\left({\cal M}(D)\right)^d\times \cal{M}_{\#}(D)$ the dual space of $C(D)\times (C(D))^d\times C_{\#}(D)$. We have for $\alpha$
\be
\alpha^*(\rho,J,p)=\sup_{c+|m|^2/2\leq 0,\;r=\Delta q} \left\{ \int_D d\rho\, c + dJ\cdot m  
+ \int_D  r p -|\nabla q|^2/2 \ dtdx\right\},
\en
which is equal to
$$\alpha^*(\rho,J,p)=\demi \int_D \left|\frac{J}{\rho}\right|^2 d\rho + \demi \int_D|\nabla p|^2 \ dtdx.$$
The term $\left|\frac{J}{\rho}\right|^2d\rho$ is defined through (\ref{4J2sur2rho}); note that this can possibly be $+\infty$. 
Then for $\beta$ we have:
\be
\beta^*(\rho,J,p)=\sup_{c,m,r}\left\{\int_D (d\rho-d\bar\rho) c + (dJ-d\bar J)\cdot m  
+\int_D (p- \bar p)r  \ dtdx\right\},
\en
the supremum being restricted to all the $c,m,r=\Delta q$ such that there exists  $\phi$ satisfying:
\be
&&c+\dt\phi+q=0,\\
&&m+\nabla_x\phi=0.
\en 
Thus in terms of $\phi,q$ we have
\be
&&\beta^*(\rho,J,p)=\\
&&\sup_{\phi,q}\left\{\int_D (d\rho-d\bar\rho) (-\dt\phi -q) - (dJ-d\bar J)\cdot \nabla_x\phi 
-\int_D \nabla q \cdot(\nabla p- \nabla\bar p) \ dtdx\right\}.
\en
Using the fact that $(\bar \rho, \bar p)$ satisfy (\ref{4c1}, \ref{4c2}) we find that
$\beta^*(\rho,J,p)=0$ if $(\rho,J,p)$ satisfies (\ref{4newweakconstr1}, \ref{4newweakconstr2}) and $\beta^*(\rho,J,p)=+\infty$ otherwise.
\bigskip
\noindent
It follows that 
\be
K=\inf_{\rho, J,p}\left\{\alpha^*(\rho,J,p)+\beta^*(\rho,J,p)\right\},
\en
where we now compute the infimum over all $(\rho, J,p)$. We have just relaxed the constraints \linebreak
(\ref{4c1}, \ref{4c2}, \ref{4c3}, \ref{4c4}) by adding the convex functional $\beta^*$ which is $+\infty$ if they are not satisfied and $0$ if they are satisfied.

\subsubsection{The duality theorem}
The functions $\alpha,\beta$ are convex with values in $]-\infty, +\infty]$.
At point $c=-1, m=0, r=0$, \mbox{$\alpha(-1,0,0)=0$},  $\alpha$ is continuous with respect to the norm
of $C(D)\times (C(D))^d\times C_{\#}(D)$,
and \linebreak  $\beta(-1,0,0)=-\int_D \bar \rho=-1$ is finite. 
The conditions to apply Fenchel-Rockafellar duality Theorem (see \cite[ch. 1]{B} ) 
are thus fulfilled and we obtain that

\be
&&\inf\{\alpha^*(\rho,J,p)+\beta^*(\rho,J,p)\}\\
=&&\sup\{-\alpha(-c,-m,-r)-\beta(c,m,r)\}\\
=&&K,
\en   
and the infimum is attained. Hence we have
\be
&&K=\sup_{c,m,r=\Delta q} \{\int_D - d\bar\rho\, c -  d\bar J \cdot m +\int_D \nabla  \bar p\cdot  \nabla q -|\nabla q|^2/2\ dtdx\}\\
&&c=-\dt\phi-q,\\
&&m=-\nabla\phi,\\
&&-c+|m|^2/2\leq 0,
\en
which is equivalent to
\be
&&K=\sup_{\phi,q}\{\int_D d\bar \rho \, (\dt\phi+q)+  d\bar J\cdot \nabla\phi  + \int_D -|\nabla q|^2/2 +  \nabla  \bar p\cdot \nabla q \ dtdx\},
\en
the supremum being performed over all $(\phi,q)$ such that
\be
\dt\phi+q+|\nabla\phi|^2/2\leq 0.
\en
To highlight the analogy with the Monge-Kantorovitch duality, and using (\ref{4c1}, \ref{4c2}, \ref{4c3}, \ref{4c4}), we rewrite the dual problem as follows: 
\be
K = \sup_{\dt\phi+q+|\nabla\phi|^2/2\leq 0} &&\left\{ \int_{\Td} \phi(T)\,d\rho_T  - \phi(0)\,d\rho_0 
+ \int_D q -|\nabla q|^2/2 \ dtdx  \right\}.
\en  
The first integral is as in the Monge-Kantorovitch problem, and the second takes into account the non-linear effect of the gravitational  coupling. This shows Proposition \ref{4propdual}.

We choose $ \bar \rho=\rho=1+\Delta p$ and $  
\bar J= \rho   v$ to be any optimal  solution ({\it i.e.} any minimizing solution).
Note that necessarily $J$ has a density $v$ with respect to $\rho$ and $v\in L^2(D,d\rho)$. This justifies the notation $J=\rho v$. 
 Then for all $\epsilon>0$ there exists $(\phi_{\epsilon},p_\epsilon)\in C^1(D)$ with $\dt\phi_\epsilon{}+p_\epsilon{}+|\nabla\phi_\epsilon|^2/2\leq 0$ such that
\beq
&&K=\frac{1}{2}\int_D d\rho   |v| ^2  + \demi\int_D|\nabla  p|^2 \ dtdx\nonumber\\
&&\leq \int_D d\rho\, (\dt\phi_{\epsilon}+p_{\epsilon})+  dJ\cdot\nabla\phi_{\epsilon}  + \int_D -|\nabla p_{\epsilon}|^2/2 +  \nabla  p\cdot \nabla p_{\epsilon} \ dtdx+\epsilon^2,\label{4supinf}
\enq
thus
\be
 &&\frac{1}{2}\int_D d\rho\, |  v-\nabla\phi_{\epsilon}|^2 +\demi \int_D |\nabla  p-\nabla p_{\epsilon}|^2 
\ dtdx\\
\leq &&\int_D  d\rho \,(\dt\phi_{\epsilon}+p_{\epsilon}+|\nabla\phi\epu|^2/2) +\epsilon^2,
\en
and we obtain
\be
\int_D  \demi  d\rho |  v-\nabla\phi_{\epsilon}|^2 + d\rho\left|\dt\phi_{\epsilon}+p_{\epsilon}+|\nabla\phi\epu|^2/2\right|   + \int_D \demi |\nabla  p-\nabla p_{\epsilon}|^2 \ dtdx \leq \epsilon^2.\label{4supinfeps}
\en
It follows that, as $\epsilon \to 0$, 
\begin{itemize}
\item $\nabla\phi_{\epsilon}$ converges to $  v$ in $L^2(D, d\rho)$,
\item $\nabla p_{\epsilon}$ converges to $\nabla  p$ in $L^2(D, dt dx)$,
\item $\dt\phi_{\epsilon}+p_{\epsilon}+|\nabla\phi\epu|^2/2$ converges to $0$ in $L^1(D, d\rho)$. 
\end{itemize}

\subsubsection{Uniqueness of the minimizer}
This property has already been proved directly, we just explain how to recover it from the dual formulation.
Notice that the sequence $(\phi_\epsilon, p_\epsilon)$ does not depend on the optimal
solution $(\rho,v)$ we have chosen, thus if we have $\nabla p_1$ and $\nabla p_2$ two optimal solutions,
then $\nabla p_\epsilon$ converges to both $\nabla p_1$ and $\nabla p_2$, and they are equal Lebesgue a.e.. It follows 
then  that two optimal solutions have the same density $\rho$. Then since 
$\nabla\phi_\epsilon$ converges to both $v_1$ and $v_2$ in $L^2(d\rho)$, $v_1$ and $v_2$
are equal $d\rho$ a.e. which proves the uniqueness of the optimal solution.
This ends the ``dual'' proof of Proposition \ref{4prop1}.

$\hfill \Box$

\section{Variations near optimality}
\subsection{Outline}
The results of the following two sections require to establish the formula (\ref{4complete}).
We present first the general approach of the proof, before entering into the rigorous details.
We first state a simple property of convex functions near their infimum:
Let $\Psi$ be a convex function, and suppose it reaches its minimum value, say $\underline \Psi$, at $x=x_0$. Assume that $D^2\Psi \geq \lambda I$. Then
\be
\demi |x-x_0|^2 \leq \frac{1}{\lambda}(\Psi(x)-\underline\Psi).
\en
The formula (\ref{4complete}) expresses this fact.
Then bounding by 0 the left hand side, we get that $\nabla\Psi(x_0) = 0$.
This yields the optimality equation as well as the conservation of energy.
Keeping the left hand side, we can also estimate $|x-x_0|$ in terms of $\Psi(x)-\underline\Psi$. This yields the formal $H^1(d\rho)$ estimate of Proposition \ref{4prop3}.
In the case studied here, a formal calculation yields
\beq
\nonumber \tilde I(\rho_2, J_2, p_2) - \tilde I(\rho_1, J_1, p_1) &=& \demi\int \frac{|J_2|^2}{\rho_2}+ |\nabla p_2|^2 
- \demi\int \frac{|J_1|^2}{\rho_1}+ |\nabla p_1|^2\\ 
&=&\nonumber\int \frac{J_1}{\rho_1}\cdot (J_2-J_1) - \demi \frac{|J_1|^2}{\rho_1^2}(\rho_2-\rho_1)+ \nabla p_1\cdot(\nabla p_2-\nabla p_1)\\
\label{4complete0} &+& \int \demi\rho_2\left|\frac{J_2}{\rho_2}-\frac{J_1}{\rho_1}\right|^2+ \demi\left|\nabla p_2-\nabla  p_1\right|^2.
\enq
This equation is the formal analogue of formula (\ref{4complete}).
The second line contains the first order terms, and the third line is a positive function.
Assume that $(\rho_1, J_1, p_1)$ is a critical point of $\tilde I$, the second line vanishes. First for $J_2=J_1 + w$ where $\nabla\cdot w=0$, we have $\rho_2=\rho_1$, this gives that $\int v_1\cdot w = 0$ for all divergence free vector field $w$, thus $v_1$ is a gradient, which we already knew, thanks to the dual formulation.
Then using that $\dt (\rho_2-\rho_1) + \nabla\cdot (J_2-J_1)=0$, with $\rho_2=\rho_1$ at times $0,T$, using $J_1/\rho_1=\nabla\phi_1$, and the Poisson equation $\Delta p = \rho-1$,  we get for the first order terms (second line)
\be
\int -\dt\phi_1 (\rho_2-\rho_1) - \demi\left|\nabla\phi_1\right|^2(\rho_2-\rho_1) - p_1(\rho_2-\rho_1) =0
\en
for all $\rho_2$. This gives the Hamilton-Jacobi equation
\be
\dt \phi_1 + \demi\left|\nabla\phi_1\right|^2 + p_1 = 0.
\en
(Note that thanks to the dual formulation, we had an approximate "$d\rho$ a.e." version of it.)
Taking the gradient of this equation gives the momentum equation (\ref{4euler-poisson})
\be
\dt v_1 + v_1\cdot\nabla v_1 = -\nabla p_1, \hspace{1cm} v_1=\nabla\phi_1.
\en
This will be the Proposition \ref{4prop2}.

Now at the critical point we keep only the third line of equation (\ref{4complete0}) (the second line vanishes for all $\rho_2$), and we get
\be
\int \demi \rho_2 |v_2-v_1|^2 + \demi\left|\nabla p_2-\nabla p_1\right|^2 = I(\rho_2,J_2) - I(\rho_1, J_1). 
\en
If we perturb $\rho_1$ in $\rho_2=(x+\delta(t,x))_{\#}\rho_1$ for a small smooth perturbation $\delta(t,x)$, 
and perturb $v_1$ in order to satisfy the conservation of mass (\ref{4continuite}),  we will obtain an estimate of the form
\be
&&\int \demi \rho_1|v_1(x+\delta(t,x)) -v_1(x)|^2 + \demi\left|\nabla p_1(x+\delta(t,x))-\nabla p_1(x)\right|^2\\
&\leq & C\int|\dt\delta|^2.
\en
This estimate is a kind of $H^1(d\rho_1)$ estimate for $v_1$, and also a $H^2$ estimate for $p_1$, 
therefore a $L^2$ estimate for $\rho_1$. 
We will then adapt techniques from \cite{DL} to obtain strong time regularity for $\rho_1$.
This will yield $\rho_1\in C(]0,T[; L^p(\Td))$ for $p<3/2$.
This will be the Proposition \ref{4prop3}.

However, all these calculations are formal, since we have to integrate $v_2$ against the measure $\rho_1$, and vice-versa, and we control $v_i$ only in $L^2(D, d\rho_i)$. Thanks to the dual formulation we will overcome this difficulty,
 by defining an extension of $v$ in $L^2_{loc}(]0,T[\times\Td)$.

\subsection{Second variation formula}

We first introduce a perturbation of the optimal path. We must perturb the optimal pair $(\rho, J)$ in such way that the conservation of mass (equation (\ref{4newweakconstr1})) is still satisfied. We proceed as follows:

Let \( \delta  \) and \( \eta  \) be two small parameters and take $\tau\in ]0,\frac{T}{2}[$.
Let \( \zeta(t)  \) be a smooth function compactly supported for $0<t<T$.
We choose $\eta$ small enough such that $t\rightarrow
t+\eta\zeta(t)$ is a diffeomorphism from $[0,T]$ to $[0,T]$. 
Let \( x\rightarrow w(x) \) be a smooth vector field and \((s,x)\to
e^{sw }(x) \) the flow associated to \( w(x) \) defined by
$$
\partial_s e^{sw }(x)=w(e^{sw }(x)) \mbox{ and } e^{0w }(x)=x,
$$ 
we can thus define $e^{\delta\zeta(t) w}(x)$. We  introduce, as in \cite{Br5}, the following measures:
\[
\rho^{\eta }(t,x)=\rho(t+\eta \zeta (t),x),\textrm{ }v^{\eta
}(t,x)=v(t+\eta \zeta (t),x)(1+\eta \dot\zeta (t)).\]
We check that the pair $(\rho^\eta, \rho^\eta v^\eta)$ satisfies the continuity equation  (\ref{4newweakconstr1}).
We define also \mbox{$p^\eta = p(t+\eta \zeta (t),x)$.}
Then we define the measures $(\rho^{\eta , \delta
},J^{\eta , \delta })$ so that for every \( f \) \( \in C(D) \) and 
$g\in (C(D))^d$ we have
\[
\int _{D}f(t,x)d\rho^{\eta ,\delta }(t,x) =\int _{D}f(t,e^{\delta \zeta
(t)w}(x))d\rho^{\eta }(t,x),
\]
and
\be
&&\int _{D}g(t,x)\cdot dJ^{\eta ,\delta }(t,x) \\
&=&\int _{D}g(t,e^{\delta \zeta
(t)w}(x))\cdot [(\partial _{t}+v^{\eta }(t,x)\cdot \nabla) e^{\delta \zeta
(t)w}(x)]d\rho^{\eta }(t,x).
\en
Note that this can be rewritten $\rho^{\eta,\delta}=e^{\delta \zeta(t)w}_{ \  \  \  \ \#}\rho^\eta$.
We check that the pair $(\rho^{\eta ,\delta },J^{\eta ,\delta })$ 
satisfies also the continuity equation (\ref{4newweakconstr1}).
Parameters $\eta,\,\delta$ being fixed, we will use the following notation:
\be
&&v^{\eta,\delta}(t,x)=(\partial _{t}+v^{\eta }(t,x)\cdot \nabla _{x})e^{\delta \zeta
(t)w}(x).
\en
(Note that we have $v^{\eta,\delta}(t,x)=v^\eta(t,x)+O(\delta)$.)
Therefore, $J^{\eta,\delta}$ can also be rewritten $e^{\delta \zeta(t)w}_{ \  \  \  \ \#}\rho^\eta v^{\eta,\delta}$.
Considering $(\phi\epu, p\epu)$ a smooth maximizing sequence for the dual problem, we have 
$$\rho^{\eta }(t,x)(\dt \phi_{\epsilon} + \frac{|\nabla\phi_{\epsilon}|^2}{2}+p_{\epsilon})(t,e^{\delta \zeta(t)w}(x))\leq 0,$$ 
and using (\ref{4supinf}) we can write:
\be
&&\frac{1}{2}\int_Dd\rho|v|^2+|\nabla p|^2\\ &\leq& \epsilon^2+ \int_D   d\rho(\dt\phi_{\epsilon}+p_{\epsilon})+ d\rho v\cdot\nabla\phi_{\epsilon}
+\nabla p\cdot\nabla p_{\epsilon} -|\nabla p_{\epsilon}|^2/2\\
&&-\int_Dd\rho^{\eta,\delta }(\dt \phi_{\epsilon} + \frac{|\nabla\phi_{\epsilon}|^2}{2}+p_{\epsilon}) .
\en
Then using the mass conservation equation (\ref{4newweakconstr1}) we have 
\be
\frac{1}{2}\int_D d\rho|v|^2+|\nabla  p|^2 \leq 
&&\epsilon^2+\int_D d\rho^{\eta}\left(v^{\eta,\delta}\cdot\nabla\phi_{\epsilon}(e^{\delta \zeta w})
-\frac{1}{2}|\nabla\phi_{\epsilon}|^2(e^{\delta \zeta w})\right)\\
+&&\int_Dd\rho p_{\epsilon} - d\rho^{\eta} p_{\epsilon}(e^{\delta \zeta
w} ) +\nabla  p\cdot\nabla p_{\epsilon} -|\nabla p_{\epsilon}|^2/2 ,
\en
so $\Delta p=\rho-1$ yields
\be
\frac{1}{2}\int_Dd\rho|v|^2+|\nabla p|^2 \leq && \epsilon^2 - \frac{1}{2}\int_Dd\rho^{\eta}\left|v^{\eta,\delta}-\nabla\phi_{\epsilon}(e^{\delta \zeta w})\right|^2 + \frac{1}{2}\int_Dd\rho^{\eta}|v^{\eta,\delta}|^2\\
+&&\int_D  p_{\epsilon}-p_{\epsilon}(e^{\delta \zeta w})\\
+&& \int_D D(e^{\delta \zeta w}):\nabla p^{\eta}\otimes\nabla p_{\epsilon}(e^{\delta \zeta w})-|\nabla p_{\epsilon}|^2/2,
\en
where $D(e^{\delta \zeta w})$ denotes the spatial derivative of $x\to e^{\delta \zeta w}(x)$. 
We obtain the complete formula:
\beq
&&\frac{1}{2}\int_Dd\rho^{\eta}\left|v^{\eta,\delta}-\nabla\phi_{\epsilon}(e^{\delta \zeta
w})\right|^2 + \frac{1}{2}\int_D\left|\nabla p^{\eta} -D(e^{\delta \zeta w})\nabla p_{\epsilon}(e^{\delta \zeta w})\right|^2 \nonumber\\
\leq &&  \epsilon^2 + \frac{1}{2}\int_Dd\rho^{\eta}|v^{\eta,\delta}|^2 -\frac{1}{2}\int_Dd\rho | v|^2\nonumber\\
+ &&  \frac{1}{2}\int_D |D(e^{\delta \zeta w})\nabla p_{\epsilon}(e^{\delta \zeta w})|^2- \frac{1}{2}\int_D |\nabla p_{\epsilon}|^2\nonumber\\
+&& \demi \int_D |\nabla p^{\eta}|^2 - \demi \int_D |\nabla p|^2\nonumber\\
+&&\int_D  p_{\epsilon}-p_{\epsilon}(e^{\delta \zeta w})\label{4complete}.
\enq

\section{Optimality equation}
In this section we prove the following:
\begin{prop}\label{4prop2}
The solution of Problem \ref{4infI2}  is a weak solution  of the Euler-Poisson system $(E-P)$ in the sense of Definition \ref{4defiweakep}.
The energy of the system defined for a.e. $t\in [0,T]$ by 
\be\label{4defeng}
E(t)=\frac{1}{2} \int_{\Td} d\rho(t,x)|v(t,x)|^2-|\nabla p(t,x)|^2 dx
\en
does not depend on time.
\end{prop}

{\it Remark.} The energy is a priori well defined in $L^1(0,T)$ since $I(\rho,v, p)$ is finite.

\subsection{Derivation of the momentum equation (\ref{4euler-poisson})}\label{4derivep}
In (\ref{4complete}), taking $\eta=0$, bounding the L.H.S. from below by $0$ and letting $\epsilon$ go to $0$ we get:

\be
0\leq&&\frac{1}{2}\int_Dd\rho |v^{0,\delta}|^2 -\frac{1}{2}\int_Dd\rho |  v|^2\\
+ &&  \frac{1}{2}\int_D |D(e^{\delta \zeta w} )\nabla p(e^{\delta \zeta  w} )|^2- \frac{1}{2}\int_D |\nabla p|^2\\
+&&\int_D  p-p(e^{\delta \zeta w} ).
\en

{\it Remark.} Notice that this inequality (or a similar one up to second order terms in $\delta$) could have been obtained directly by expressing that 
$$\tilde I(\rho^\delta, J^\delta,p^\delta )\geq \tilde I(\rho, J, p),$$
where $p^\delta=\Delta^{-1}(\rho^\delta -1)$. Notice however that the expression of $p^\delta$ is not so straightforward.

\bigskip

Expanding $e^{\delta \zeta(t) w}(x)=x+\delta \zeta(t) w(x)+O(\delta^2)$, we get
\be
&&\frac{1}{2}\int_D |D(e^{\delta \zeta w})\nabla p(e^{\delta \zeta w})|^2\\
=&&\frac{1}{2}\int_D \left|(I+\delta\zeta Dw)(e^{-\delta \zeta w})\nabla p\right|^2 \mathbf{J}(e^{-\delta \zeta w}) +O(\delta^2)\\
\en
with $I$ the identity matrix of order $d$ and $\mathbf{J}(e^{-\delta \zeta(t)w})$ the jacobian determinant of the mapping $x\rightarrow e^{-\delta \zeta(t)w}(x)$.
Using that $\mathbf{J}(e^{-\delta \zeta(t)w})=1-\delta\zeta(t)\nabla\cdot w +O(\delta^2),$ this is equal to
\be
\frac{1}{2}\int_D (|\nabla p|^2 + 2\delta\zeta Dw: \nabla p \otimes\nabla p -\delta\zeta|\nabla p|^2 \nabla\cdot w)+O(\delta^2).
\en
Then for $v$ we have
\be
&&\frac{1}{2}\int_Dd\rho|v^{\delta}|^2 -\frac{1}{2}\int_Dd\rho |v|^2\\
=&&\int_D d\rho v\cdot(\dt + v\cdot\nabla)\delta\zeta w+O(\delta^2)
\en
and 
\be
\int_D -p(e^{\delta\zeta w})=-\int_D p -\int_D \nabla p\cdot \delta\zeta w +O(\delta^2).
\en
This yields finally
\be
0\leq &&\delta\left[ \int_D  d\rho v\cdot(\dt+v\cdot \nabla) \zeta w+\zeta  Dw :\nabla p\otimes\nabla p \right.\\
&&\left. \  \  -\demi|\nabla p|^2\nabla\cdot \zeta w - \nabla p\cdot\zeta w \ dtdx  \right]+O(\delta^2).
\en
Thus for every $w$ smooth vector field on $\Td$, $\zeta\in C^{\infty}_c(0,T)$,  we have
\be
0&=&\int_D d\rho v\cdot(\dt+v\cdot \nabla)(\zeta w)+\zeta Dw: \nabla p\otimes\nabla p\\
&& -\demi|\nabla p|^2\nabla\cdot \zeta w -\nabla p\cdot \zeta w \ dtdx,
\en
hence we conclude that $(\rho,v,p)$ is a weak solution  of the Euler-Poisson system in the sense of Definition \ref{4defiweakep}. 
Note that in the proof of Proposition \ref{4prop1}, we have already shown that \linebreak  $\rho\in C([0,T]; {\cal P}-w*)$.

$\hfill\Box$

\subsection{Conservation of energy}
Here we shall deduce conservation of energy by using a wrinkle in time to perturb the minimizing path;
we take $\delta=0$ in (\ref{4complete}), minorize the LHS by 0 and let $\epsilon$ go to 0 to obtain
\be
0\leq &&\frac{1}{2}\int_D d\rho(t+\eta\zeta(t),x)(1+\eta\dot\zeta(t))^2|v(t+\eta\zeta(t),x))|^2 
\\
&+& \frac{1}{2}\int_D |\nabla p(t+\eta\zeta(t),x)|^2 \ dtdx\\
& -& \frac{1}{2}\int_D d\rho(t,x)|v(t,x)|^2 -\frac{1}{2}\int_D |\nabla p(t,x)|^2 \ dtdx .
\en
Changing variable in time $ t:=t+\eta\zeta(t)$, $dt:=dt(1+\eta\dot\zeta(t))$ we get
\be 
0\leq &&\frac{1}{2}\int_D d\rho(t,x)|v(t ,x)|^2 \eta\dot\zeta(t) + |\nabla p(t,x)|^2 (\frac{1}{1+\eta\dot\zeta(t)}-1) \ dtdx.
\en
Taking the first order term in $\eta$ we get 
\be 
\frac{1}{2} \int_D \left[d\rho(t,x)|v(t,x)|^2-|\nabla p(t,x)|^2 dtdx\right]\dot\zeta(t) \ dtdx  =0
\en
for any $\zeta\in C^{\infty}_c(0,T)$
which gives the conservation of energy, with 
$$E=\frac{1}{2} \int_{\Td} d\rho(t,x)|v(t,x)|^2-|\nabla p(t,x)|^2dx.$$
This ends the proof of Proposition \ref{4prop2}.

$\hfill \Box$

\section{Regularity properties of the minimizer}
In this section we obtain several regularity properties of solutions of Problem \ref{4infI1}. 
Those properties will follow from equation (\ref{4complete}), the rigorous analog of (\ref{4complete0}). Here we use in a crucial way that our solution is a minimizer of the action of the Lagrangian and not only a critical point. 
For points 1 and 2 we follow closely the method of Brenier in \cite{Br5} where similar results were obtained in the case of the Euler incompressible equation.
For the fourth point we use the first two points and a method close to the one used by DiPerna and Lions in \cite{DL}.
Similar results have also been obtained by different techniques in \cite{EG} for finite dimensional hamiltonian systems, using 
properties of a special Hamilton-Jacobi equation related to the Hamiltonian flow.
This section is thus devoted to the proof of the following:
\begin{prop}\label{4prop3}
The optimal solution $(\rho,J=\rho v)$ of Problem \ref{4infI2} has the following regularity properties:
\begin{enumerate}
\item The density $\rho$ belongs to $L^2_{loc}(]0,T[;L^2(\Td))$,
\item the velocity $v$ can be extended in all of $\Td$ to a function of $L^2_{loc}(]0,T[;L^2(\Td))$,
in such a way that for all $\tau$ in $]0,T/2]$, for all $y$ in $\R^d$,\\ 
$\displaystyle\int_{\Td}\int_{\tau}^{T-\tau}\rho(t,x) \left|v(t,x+y)-v(t,x)\right|^2 dt dx\leq C_\tau|y|^2$, 
\item the velocity potential $\phi$ can be chosen in $\Linf_{loc}(]0,T[\times\Td)$,
\item finally $\rho$ belongs to $C(]0,T[;L^p(\Td))$ for any $p\in [1,3/2[.$
\end{enumerate}
\end{prop}

\subsection{Spatial regularity: proof of points 1 and 2}
We are going to deduce spatial regularity using the time dependent uniform translation \linebreak $e^{\delta \zeta  w}(x)=x+\delta \zeta(t)y$ to perturb the minimizing path.
This corresponds to the case  $\eta=0$, $\zeta\equiv 1$ in $[\tau,T-\tau]$,  $w(x)=y$ fixed in inequality (\ref{4complete}). In this case $v^\delta(t,x)=v(t,x)+\delta\dot \zeta(t) y$, and $D(e^{\delta\zeta w})= I$, and  inequality (\ref{4complete}) becomes
\be
&&\frac{1}{2}\int_Dd\rho(t,x)|v^\delta(t,x)-\nabla\phi_{\epsilon}(t,x+\delta \zeta
(t)y)|^2  \\
&+& \frac{1}{2}\int_D|\nabla p(t,x) -\nabla p_{\epsilon}(t,x+\delta \zeta(t)y)|^2  \ dtdx\\
&\leq&  \epsilon^2 +  \demi\int_D d \rho (|v+\delta\dot{\zeta}y|^2 - |v|^2)\\
&=&  \epsilon^2 +  \demi\int_Dd \rho |\delta\dot{\zeta}y|^2.
\en
For the last line, we have used $\ds\Dt \int_{\Td} d\rho \, v \equiv 0$; indeed
take $\varphi=\zeta(t)y$ in  the momentum equation (\ref{4weakep}), and use $\ds\int_{\Td} \nabla p  \equiv 0$.
We have also
\be
\int_D d\rho|v^\delta -v|^2 = \int_Dd\rho |\delta\dot{\zeta}y|^2  \leq \frac{C}{\tau}\delta^2 |y|^2 
\en 
for a suitable choice of $\zeta$.
Hence 
\beq\label{4regepsilon}
&&\int_{\Td}\int_{\tau}^{T-\tau}d\rho(t,x)\left|\nabla\phi_{\epsilon}(t,x+y)-v(t,x)\right|^2 
+ \left|\nabla p_{\epsilon}(t,x+y)-\nabla p(t,x)\right|^2 \ dtdx \nonumber \\
&\leq& \epsilon^2 + \frac{C}{\tau}|y|^2\nonumber.
\enq
We let $\epsilon$ go to 0 and obtain 
\be
\int_{\Td}\int_{\tau}^{T-\tau} \left|\nabla p(t,x+y)-\nabla p(t,x)\right|^2 dt dx\leq \frac{C}{\tau}|y|^2,
\en
thus $D^2p$, and $\rho =1+\Delta p$ are in $L^2_{loc}(]0,T[;L^2(\Td))$. In particular, $\rho$ is absolutely continuous with respect to the Lebesgue measure of $D$.
We will also obtain that  $\nabla\phi_{\epsilon}$ is bounded in $L^2_{loc}(]0,T[;L^2(\Td))$. Indeed,
we get first from (\ref{4regepsilon}) that 
\beq\label{4epsy}
\int_{\Td}\int_{\tau}^{T-\tau}\rho(t,x)\left|\nabla\phi_{\epsilon}(t,x+y)\right|^2 \ dxdt \leq C(1+\frac{1}{\tau})
\enq
for $\epsilon\leq 1$. Integrating
this over $y\in \Td$ we get
\be
&&\int_{y\in\Td}\int_{x\in\Td}\int_{\tau}^{T-\tau} \rho(t,x+y)\left|\nabla\phi_{\epsilon}(t,x)\right|^2dxdydt\nonumber\\
=&&\int_{x\in\Td}\int_{\tau}^{T-\tau}\left|\nabla\phi_{\epsilon}(t,x)\right|^2dxdt\leq C(1+\frac{1}{\tau})\label{4gradphiL2},
\en 
thus we can, up to extraction of a subsequence,  define a weak limit for this sequence, as $\epsilon$ goes to 0,
$v=\nabla\phi\in L^2_{loc}(]0,T[\times \Td)$. 
However $v$ will  be uniquely defined only in the $d\rho$ a.e. sense.
Moreover for all $y\in \Td$, $\nabla\phi_{\epsilon}(\cdot + y)$ converges to $\nabla\phi(\cdot + y)$
in $L^2(d\rho)$ weak: indeed, from (\ref{4epsy}), the sequence $\nabla\phi_{\epsilon}(\cdot + y)$ is uniformly bounded in $L^2(d\rho)$,
and for all $\varphi \in C^{\infty}_c(]0,T[\times\Td)$, we have
\be
\int \rho(t,x) \nabla\phi_{\epsilon}(x + y)\cdot \varphi(t,x) \ dtdx \to \int \rho(t,x) \nabla\phi(x + y)\cdot \varphi(t,x) 
\ dtdx
\en
since $\rho\in L^2_{loc}(]0,T[\times\Td)$, and since $\nabla\phi\epu$ converges weakly to $v$ in $L^2_{loc}(]0,T[\times\Td)$.
Using that the $L^2(d\rho)$ norm is l.s.c. with respect to the weak $L^2(d\rho)$ convergence,
we will have 
\be
\int \rho(t,x) |\nabla\phi(x + y)|^2 \ dtdx \leq \liminf \int \rho(t,x) |\nabla\phi_{\epsilon}(x + y)|^2 \ dtdx. 
\en
Therefore \mbox{$(t,x) \to \nabla\phi(t,x +y)$} will be in $ L^2_{loc}(]0,T[;L^2(dx)\cap L^2(d\rho(t)))$ for any $y\in \Td$. 
We finally obtain
\beq\label{4spaceregul}
\int_{\Td}\int_{\tau}^{T-\tau}\rho(t,x)\left|v(t,x+y)-v(t,x)\right|^2 \ dtdx \leq C_{\tau}|y|^2.
\enq
This proves the first two points of Proposition \ref{4prop3}.

\subsection{$\Linf$ bound for the potential $\phi$: proof of point 3}
\label{4rho_linf}
We assume here that $d\leq 3$. We normalize $p$ so that its mean value is 0 for all time. 
Then since $\rho\in L^2_{loc}(]0,T[\times \Td)$ we have
$p  \in L^2_{loc}(]0,T[;H^2(\Td))$ which is continuously embedded in
$L^2_{loc}(]0,T[;C^{\demi}(\Td))$ and thus $\|p(t,.)\|_{\Linf}\in L^2_{loc}(]0,T[)$.
(We will see after that $\rho$ is in $\Linf_{loc}(]0,T[\times\Td)$ and we will be able to remove this assumption on the dimension.)

We take a regularization in $t,x$ of $\phi,p$:
on $]\tau, T-\tau[,$
\be
&&\phi\ep(t,x)=\eta\epu*\phi,\\
&& p\ep(t,x)=\eta\epu*p,\\
&&\eta\epu(t,x)=\frac{1}{\epsilon^{d+1}}\eta_1(\frac{t}{\epsilon},\frac{x}{\epsilon}),
\en
with $\eta_1$ compactly supported in $[-1,1]\times B(0,1)$ ,  and $0<\epsilon<\tau/2$. We first check that 
\beq\label{4ineqhjep}
\dt\phi\ep +\demi|\nabla\phi\ep|^2 +p\ep\leq 0.
\enq
Considering $\sigma \to \phi\ep(\sigma, \gamma(\sigma))$ 
with $\gamma\in C^1([0,T];\Td)$ we have
\be
\frac{d}{d\sigma}(\phi\ep(\sigma, \gamma(\sigma)))&=&\dt \phi\ep(\sigma, \gamma(\sigma))
+ \dot\gamma\cdot \nabla \phi\ep(\sigma, \gamma(\sigma))\\
&\leq& \dt \phi\ep(\sigma, \gamma(\sigma))+\demi|\nabla\phi\ep|^2(\sigma, \gamma(\sigma))+ \demi|\dot\gamma(\sigma)|^2\\
&\leq&-  p\ep(\sigma, \gamma(\sigma)) + \demi|\dot\gamma(\sigma)|^2
\en
using (\ref{4ineqhjep}), and we obtain 
\be
\phi\ep(t+s,x)&\leq& \inf_{\gamma \in \Gamma}\left\{\phi\ep(t,\gamma(t))+\int_t^{t+s}  -  p\ep(\sigma, \gamma(\sigma))+ \demi|\dot\gamma|^2(\sigma) \ d \sigma\right\},
\en
with $\Gamma$ the set of all continuous paths going from $[t,t+s]$ to $\Td$ such that $\gamma(t+s)=x$.
Then restricting the infimum to paths of the form $\gamma(\sigma)=\gamma(t)+\frac{\sigma-t}{s}(x-\gamma(t))$ 
and noticing that \linebreak $\|p\ep(t,.)\|_{\Linf}\in L^2_{loc}(]0,T[)$ (uniformly in $\epsilon$) implies that 
$\int_t^{t+s}| p\ep(\sigma,\gamma(\sigma))|d\sigma \leq C\sqrt{s}$,
we obtain the following upper bound:
\beq
\phi\ep(t+s,x)\leq \inf_{z\in\Td}\{\phi\ep(t,z)\}+  C(\frac{1}{s}+\sqrt{s})\label{4phiphi}.
\enq
A simple computation shows that
\beq
&&\int_{\Td}d\rho(t_2,x)\phi\ep(t_2,x) -\int_{\Td}d\rho(t_1,x)\phi\ep(t_1,x)  \label{4tendversK}\\
&&\to \int_{t_1}^{t_2}\int_{\Td} \demi |\nabla\phi|^2d\rho + |\nabla p|^2dtdx \nonumber 
\enq
as $\epsilon$ goes to 0, in particular it converges to a non-negative value, and will be greater than $-1$ for $\epsilon$ small enough.
It costs no generality to normalize $\phi\ep$ such  that $$\int_{\Td}d\rho(T/2,x)\phi\ep(T/2,x)=0.$$
Then by choosing $(t_1,t_2) = (\tau/2, T/2)$  and $(t_1,t_2) = (T/2, T-\tau/2)$ in (\ref{4tendversK}),
we get respectively that, for $\epsilon$ small enough, 
\beq
&&\int_{\Td} d\rho(\tau/2,x)\phi\ep(\tau/2,x) \leq 1,\label{4aa}\\
&&\int_{\Td} d\rho(T-\tau/2,x)\phi\ep(T-\tau/2,x) \geq -1.\label{4bb}
\enq
This implies 
\be
&&\inf_x \{\phi\ep(\tau/2,x)\} \leq 1,\\
&& \sup_x \{\phi\ep(T-\tau/2,x)\} \geq -1,
\en
and using (\ref{4phiphi}) we obtain
\beq
\|\phi\ep\|_{\Linf([\tau,T-\tau]\times\Td)}\leq C(\tau).\label{4Linfboundphi}
\enq
This leads to the conclusion, since $\phi\ep$ converges $dt dx$ a.e. to $\phi$.

{\it Remark 1.} This regularization of $\phi,p$ is needed to enforce the uniform bound on $p\ep$ in $H^2$. This might not have been satisfied by a smooth maximizing sequence of the dual problem.

{\it Remark 2: $L^2$ bound.}
Here we give also a bound when we do not assume $d\leq 3$. In this case we will have, from the $L^2_{loc}$ bound on $\rho$,  
$p\ep\in L^2_{loc}(]0,T[; L^k(\Td))$ for $k=\frac{2d}{d-4}>2$. We then consider, for $s,t\in ]0,T[$, two probability measures $\mu_s$ and $\mu_t$. We consider the geodesic path between $\mu_{s}$ and $\mu_t$, that we denote $\mu_{s'}, s'\in [s,t]$, obtained from Proposition \ref{mccann} (with proper time renormalization). The corresponding velocity
$v(s'), s'\in [s,t]$ satisfies $\dt \mu + \nabla \cdot(\mu v)=0$.
Combined with (\ref{4ineqhjep}), this will yield
\be
\int_{\Td} \phi\ep(t)d\mu(t) \leq \int_{\Td} \phi\ep(s) d\mu(s) + \int_s^t\int_{\Td} (\demi |v(s')|^2-p\ep(s'))d\mu(s')ds'.
\en
Using the properties of the geodesic path $\mu$ expressed in Proposition \ref{mccann},
we are able to obtain that
\be
\int_{\Td} \phi\ep(t)d\mu(t) - \int_{\Td} \phi\ep(s)d\mu(s) \leq C(\tau, |t-s|) \sup \{\|\mu(s)\|_{L^2},\|\mu(t)\|_{L^2}\},
\en
with $C(\tau, |t-s|)$ uniformly bounded for $\tau >\tau_0, \tau \leq s,t \leq T-\tau, t-s\geq \delta_0 >0$. (We have used in particular that $|v| \leq C(d)/|t-s|$.)
In (\ref{4aa}, \ref{4bb}) above, we choose $\tau$ such that $\rho(\tau/2), \rho(T-\tau/2)$ both belong to $L^2$ with norm bounded by some constant  $C(\tau)$. Since 
$\rho\in L^2_{loc}(]0,T[\times\Td)$,  this is true for almost every $\tau>0$.

From this, taking successively $\mu(s)=\rho(\tau/2), \mu(t)=\rho(T-\tau/2)$, we obtain that for all probability measure $\mu$, for all $t\in [\tau, T-\tau]$, $\left|\int \phi\ep(t)d\mu\right| \leq C(\tau) \|\mu\|_{L^2}$.
This then yields 
\beq\label{4L2boundonphi}
\|\phi\|_{\Linf([\tau, T-\tau];L^2(\Td))} \leq C(\tau).
\enq

\subsection{Strong time continuity of $\rho$: proof of point 4}
First let us notice that in inequality (\ref{4complete}), if we do not set $\eta=0$,  we obtain the mixed derivative estimate: 
\beq
&&\int_{t=\tau}^{T-\tau}\int_{\Td}\rho(t,x)\left|v(t+\eta,x+y)-v(t,x)\right|^2\\
&&+\left|\nabla p_{\epsilon}(t+\eta,x+y)-\nabla p(t,x)\right|^2 dt dx\nonumber\\
&\leq& C_{\tau}(|y|^2+|\eta|^2).\nonumber
\enq
By taking $y=0$ we get that $\partial_t \nabla p \in L^2_{loc}(]0,T[;L^2(\Td)).$ This implies
\begin{lemme}\label{4Ep}
The gravitational field $\nabla p$ belongs to  $C^{\demi}_{loc}(]0,T[;L^2(\Td))$. 
\end{lemme}

Now we prove the last part of Proposition \ref{4prop3}. We obtain that 
the density is strongly continuous with respect to time by showing some renormalization property,
that gives the continuity of some $L^p$ norm, with $p>1$.
Note that in the non-interacting case, since the functions of time $(k-1)\|\rho(t,\cdot)\|^k_{L^k}$ are convex, they are continuous in $]0,T[$ provided they are finite at $t=0$ or $t=T$.  
So we will prove:
\begin{lemme}\label{4renorm} 
Let $\alpha\in [1,3/2[$ and  $\ds G^\alpha(t)=\int_{\Td}\rho^\alpha(t,x) \ dx$, then  $G\in C(]0,T[)$.
\end{lemme}
We postpone the proof of the lemma after the proof of the last point of Proposition \ref{4prop3}.

\bigskip
\noindent
{\bf Proof of the last point of Proposition \ref{4prop3}} First we check the weak time continuity of $\rho$:
from the conservation of energy
and from Lemma \ref{4Ep}, $\int_{\Td}\rho |v|^2$ is uniformly bounded on $[\tau,T-\tau]$. Thus $\rho |v| =  \sqrt{\rho}\sqrt{\rho |v|^2} \in L^{\infty}([\tau,T-\tau];L^p)$ for some $p>1$ thanks to Lemma \ref{4renorm}.
It follows that from equation (\ref{4continuite}), $\dt \rho$ is bounded in $L^{\infty}([\tau,T-\tau],H^{-s})$ for some $s$.
Using classical arguments of functional analysis (see \cite{Li}) we can deduce that, for some $p>1$,
\be
\rho \in C(]0,T[; L^p - w).
\en 

Then Lemma \ref{4renorm} implies that $\rho \in C(]0,T[;L^p)$ for any $p\in [1,3/2[$:
indeed it is a classical fact that when a sequence $(u_n)_{n\in\N}$ converges weakly in $L^p, 1<p<\infty$,  to some $u$, if $\|u_n\|_{L^p}$ converges to $\|u\|_{L^p}$ the sequence converges strongly. The last point of Proposition \ref{4prop3} is thus proved. 

$\hfill \Box$

\bigskip
\noindent 
{\bf Proof of Lemma \ref{4renorm} } Let us prove the renormalization property when  the density and the velocity field are smooth: we use the identity
\be
\dt [\rho F(\rho)] + \nabla\cdot[\rho F(\rho) v] = -\rho^2F'(\rho)\nabla\cdot v.
\en
Integrating over $\Td$ we get that $\int_{\Td} \rho F(\rho)$ is continuous with respect to time as long as $\rho^2F'(\rho)\nabla\cdot v$ is in $L^1_{loc}([0,T]\times\Td)$. We will see that this is true for  $F(\rho) =\rho^\beta$, $\beta < \demi$ from the regularity property  (\ref{4spaceregul}).
\\
We introduce $\eta(x)=C\exp(-\frac{|x|^2}{\sqrt{1+|x|^2}})$ with C such that $\int_{\Rd}\eta(x)dx=1$.
Then as usual \linebreak $\eta_\epsilon(x)=\frac{1}{\epsilon^d}\eta(\frac{x}{\epsilon})$, and 
$\rho_\epsilon(x)=\rho * \eta_\epsilon(x)=\int_{\Rd}\rho(x-y)\eta_\epsilon(y)dy$,
$\rho$ being naturally extended to a $\Zd$ periodic function on all of $\Rd$.
With $(\rho, v)$ as before, we consider the pair $(\rho\epu,v\epu)$ defined by
\be
&&\rho\epu=\eta\epu* \rho,\\
&&v\epu=\eta\epu* (\rho v)/\rho\epu.
\en
We check that the pair $(\rho\epu, v\epu)$ still satisfies the mass conservation equation (\ref{4continuite}). 
Then we have the crucial property:
\begin{lemme}\label{4regulcarac}
For $\tau \in ]0,T/2[$, 
\be
\int_{\tau}^{T-\tau}\int_{\Td}\rho\epu|\nabla\cdot v\epu|^2\leq C(\tau).
\en
\end{lemme}
This lemma means that the spatial regularity property (\ref{4spaceregul}) is conserved through regularization.
Before proving this lemma, we conclude the proof of Lemma \ref{4renorm}:
since $\rho\epu>0$, for all $F \in C^1(]0,+\infty[)$ we have 
\be
&&\Dt\int \rho\epu F(\rho\epu)= - \int \rho\epu^2F'(\rho\epu)\nabla\cdot v\epu,\\
&&\Dt \int  \rho\epu^\alpha = (1-\alpha)\int \rho\epu^\alpha\nabla\cdot v\epu.
\en
Then using Lemma \ref{4regulcarac} and the fact that $\rho\epu \in L^2_{loc}(]0,T[\times\Td)$ we get that$$\int_{\tau}^{T-\tau}\int_{\Td}\rho\epu^{3/2}|\nabla\cdot v\epu|\leq C(\tau),$$
and also that for any $\alpha \in ]1/2,3/2[$, the sequence $\rho\epu^{\alpha}|\nabla\cdot v\epu|$
is equiintegrable on $[\tau,T-\tau]\times\Td$.
Thus the sequence $t\rightarrow \int_{\Td}\rho^{\alpha}\epu(t,x)dx$ is equicontinuous and the limit
$t\rightarrow \int_{\Td}\rho^{\alpha}(t,x)dx$  is continuous.

$\hfill \Box$
\\
{\bf Proof of Lemma \ref{4regulcarac}:}

\noindent
We have
\be
 \nabla\cdot  v\epu &=& -\frac{\nabla\rho\epu}{\rho\epu^2}\cdot (\rho v)* \eta\epu + 
\frac{1}{\rho\epu}(\rho v)* \nabla\eta\epu\\
&=&-\frac{\nabla\rho\epu}{\rho\epu^2}\cdot\int_{\Rd} \rho(x-y)(v(x-y)-v(x))\eta\epu(y)
-\frac{\nabla\rho\epu}{\rho\epu}\cdot v\\
&&+\frac{\nabla\rho\epu}{\rho\epu}\cdot v+\frac{1}{\rho\epu}\int_{\Rd}\rho(x-y)(v(x-y)-v(x))\cdot\nabla\eta\epu(y).
\en
We use the special shape of the regularization kernel:
there exists $C>0$ such that for all $x\in \Rd$, $|\nabla\eta(x)|\leq C \eta(x)$. This implies 
\be
\left|\frac{\nabla\rho\epu}{\rho\epu}\right|\leq \frac{C}{\epsilon}.
\en
We also have the usual bounds
\be
&&\int |y| |\nabla\eta\epu(y)|dy \leq C,\\
&& \int |y|^2\eta\epu(y)dy \leq C\epsilon^2.
\en
We define 
\be
A(x)=\frac{\nabla\rho\epu}{\rho\epu^2}\int_{\Rd} \rho(x-y)(v(x-y)-v(x))\eta\epu(y) \ dy.
\en
Then 
\be
A^2(x)\leq C\left|\int_{\Rd} \frac{\rho(x-y)\eta\epu(y)}{\rho\epu(x)}\frac{|v(x-y)-v(x)|}{\epsilon} \ dy\right|^2,
\en
and by Jensen's inequality this is less than
\be
C\int_{\Rd} \frac{\rho(x-y)\eta\epu(y)}{\rho\epu(x)}\frac{|v(x-y)-v(x)|^2}{\epsilon^2} \ dy,
\en
hence we obtain
\be
\int_{\Td} \rho\epu(x)A^2(x)dx \leq \int_{\Td}\int_{\Rd} \rho(x)\eta\epu(y)\frac{|v(x+y)-v(x)|^2}{\epsilon^2} \ dydx.
\en
For the next term
\be
B(x)=\frac{1}{\rho\epu}\int\rho(x-y)(v(x-y)-v(x))\nabla\eta\epu(y) \ dy,
\en we proceed by the same method.
We now have the bound
\be
&&\int_{\tau}^{T-\tau}\int_{\Td} \rho\epu |\nabla\cdot v\epu|^2(x) \ dxdt \\
&\leq& C \int_{\tau}^{T-\tau}\int_{\Td}\int_{\Rd} \rho(t,x)\frac{1}{\epsilon^2}|v(t,x+y)-v(t,x)|^2\eta\epu(y) \ dxdydt \\ 
&=& C \int_{\tau}^{T-\tau}\int_{\Td}\int_{\Rd} \rho(t,x)\frac{1}{|y|^2}|v(t,x+y)-v(t,x)|^2\eta\epu(y)\frac{|y|^2}{\epsilon^2} \ dxdydt \\
&\leq& C(\tau)
\en
using the spatial regularity property (\ref{4spaceregul}). The proof of Lemma \ref{4regulcarac} is  complete. 

$\hfill \Box$


\section{Consistency with smooth solutions of the Euler-Poisson system}
Here we show that the solution of the minimization problem co\"incides with a smooth potential solution of the Euler-Poisson system satisfying the boundary conditions, when the latter exists. We will try to make precise the smoothness required to reach this conclusion.  
\begin{theo}\label{4main3}
Let $(\nabla\phi,\rho,p)$ be the solution of Problem \ref{4infI2}, and let 
$(\nabla\psi,r,q)$ be such that
$$\psi\in W^{1,\infty}(D), \  r\in L^2([0,T]; H^{-1}(\Td))\cap \Linf([0,T]; L^1(\Td)), \ q\in \Linf(D).$$
Suppose that $(\nabla\psi,r,q)$ is a solution to 
\beq
&&\dt \psi + \frac{1}{2}|\nabla\psi|^2 + q \leq 0\label{4HJ},\\
&&r(\dt \psi + \frac{1}{2}|\nabla\psi|^2 + q)=0\label{4HJ=},\\
&&\dt r +\nabla\cdot (r\nabla\psi)=0\label{4contipsi} \text{ in }{\cal D}', \\
&&r|_{t=0}=\rho_0,\,r|_{t=1}=\rho_1\label{4boutpsi},\\
&&\Delta q = r-1\label{4poissonpsi}\text{ in }{\cal D}'.
\enq
Then $\rho= r$ and $\nabla\phi=\nabla\psi$ holds $d\rho$ a.e..
\end{theo}

{\it Proof.}
Since  $(\rho, \nabla\phi)$ satisfies the continuity equation (\ref{4continuite}) we have
\be
\int_D \rho \nabla\phi \cdot \nabla\psi  + \int_D \rho \dt \psi  = \int_{\Td} \rho_{T}\psi(T)-\int_{\Td}\rho_{0}\psi(0).
\en
Combining with (\ref{4HJ}) this yields
\be
\int_D \rho \nabla\phi \cdot \nabla\psi - \int_D \rho (\frac{1}{2}|\nabla\psi|^2 + q) \geq  \int_{\Td} \rho_{T}\psi(T)-\int_{\Td}\rho_{0}\psi(0). 
\en
Using the Poisson equation (\ref{4poisson}) we have  
\be
&&\frac{1}{2}\int_D \rho(-|\nabla\phi-\nabla\psi|^2 + |\nabla\phi|^2) + \int_D \nabla p \cdot \nabla q -q\\
 \geq &&
 \int_{\Td} \rho_{T}\psi(T)-\int_{\Td}\rho_{0}\psi(0) \\
= && \int_D r \nabla\psi \cdot \nabla\psi  + \int_D r \dt \psi \\
= &&-\int_D q r + \frac{1}{2}\int_D r |\nabla\psi|^2\\
= && \int_D |\nabla q|^2 -q + \frac{1}{2}\int_D r |\nabla\psi|^2,
\en
where we have used (\ref{4contipsi}, \ref{4boutpsi}) in the third line, (\ref{4HJ=}) in the fourth line,  and (\ref{4poissonpsi}) in the fifth line.
We finally get 
\be
&&\int_D \rho |\nabla\phi-\nabla\psi|^2 + \int_D |\nabla p-\nabla q|^2\\
&\leq & \int_D (\rho |\nabla\phi|^2+ |\nabla p|^2) - \int_D (r |\nabla\psi|^2 + |\nabla q|^2).
\en
Since $(\rho,\nabla \phi,p)$ is solution of the minimization problem, the RHS is non positive, and we obtain the expected result.

{\it Remark 1.} This is true in particular if $(\psi,q)$ is a $C^2(D)\times C^1(D)$ solution of 
$\dt \psi +\demi|\nabla\psi|^2 + q =0$ and thus shows the consistency with smooth solutions of the $(E-P)$ system.

{\it Remark 2.} From the results of Theorem \ref{4main2} 
depending on $\rho_0$ and $\rho_T$ the assumptions on $\psi, r,q$ can 
be weakened. For instance if $\rho_0$ and $\rho_T$ are in $\Linf$ then $\rho, \nabla\phi, \nabla p$ are in $\Linf([0,T]\times\Td)$ and one only needs $(1+r)|\nabla\psi|^2,  (1+ r)|\dt \psi|, |\nabla q|^2$ to be integrable to perform our computation. (Note that these assumptions imply that $\psi \in C([0,T]; L^1(\Td))$ and thus  $\int_{\Td} \rho_{T}\psi(T)-\int_{\Td}\rho_{0}\psi(0)$ is well defined.)
\\
This ends the proof of the Theorem \ref{4main3}.

$\hfill \Box$


\section{Path regularity for a time discretized interaction: Proof of Theorem \ref{4main2}}
In this section we prove several additional regularity properties for the variational solution.
The problems that we will treat are closely related to viscosity solutions of Hamilton-Jacobi equation.
In the remainder we will denote by HJ1, HJ2 the following operators:
\be
&&\textrm{ HJ1 }(\phi)=\dt\phi +\demi|\nabla\phi|^2, \\
&&\textrm{ HJ2 }(\phi,p)=\dt\phi +\demi|\nabla\phi|^2 + p, \\
\en
and we will make precise the sense in which they must be understood.
See \cite{E} for references about Hamilton-Jacobi equations.

\subsection{Formal bounds}
Solutions of our variational problem satisfy
\beq
&&\dt\phi +\demi|\nabla\phi|^2 +p = 0\label{41reg}  \  d\rho \ a.e.,\\
&&\dt\phi +\demi|\nabla\phi|^2 +p \leq  0\label{41'reg},\\
&&\dt\rho +\nabla\cdot(\rho\nabla\phi)=0\label{42reg},\\
&&\rho=1+\Delta p\label{43reg}.
\enq
If some $C^2(\Td)$ function $Z$ satisfies $Z(x_0)=0$  and $Z\leq 0$ in $\Td$ 
then $D^2(Z)(x_0)\leq 0$ in the sense of matrices and in particular this implies that
$\Delta Z(x_0)\leq 0$.
Using  (\ref{41reg}, \ref{41'reg}) and applying this result to $Z=\dt\phi +\demi|\nabla\phi|^2 +p$ one formally obtains
\be
\dt \Delta\phi +(\nabla\phi\cdot\nabla)\Delta\phi + \sum_{ij}|\partial_{ij}\phi|^2 + \Delta p\leq 0 \  d\rho \ a.e.,
\en
which combined with the inequality $|\Delta \phi|^2\leq d \sum_{i,j}|\partial_{ij}\phi|^2$ 
and with (\ref{43reg}) gives 
the inequality
\beq\label{4ineq}
\Dt \Delta\phi \leq 1-\rho - \frac{1}{d}|\Delta \phi|^2,
\enq
where the operator $\ds\Dt\cdot=\dt\cdot + \sum_{i=1}^d \partial_i \phi \  \partial_i\cdot  $ denotes the convective derivative along the flow generated by the velocity field 
$\nabla\phi(t,x)$.
We first deduce from this the following  upper bound for $\Delta\phi$ $d \rho$ a.e.: 
\be
\Delta \phi \leq C(d)(1+\frac{1}{t}),
\en
it is obtained by looking at the behavior of the differential inequality
\be
\dot {f} \leq 1-\frac{1}{d} f^2
\en
for large $f$. 
This "Oleinik-type" estimate is well known for viscosity solutions of HJ2=0 provided that $p$ (as is the case here) satisfies
$\Delta p \geq C$. However we don't know a priori that our solution is a viscosity solution, and moreover this bound is true in the sense of distributions. Here a complication is added by the fact that the solution $\phi, p$ satisfies HJ2$(\phi,p)=0$ only $d\rho$ a.e., and therefore {\it a priori} not in the viscosity sense.

Now notice  that our variational  solution exists on $t\in [0,T]$ and thus 
 $(\psi,q)(t)=(-\phi,p)(T-t)$ is also a solution to equations (\ref{41reg}) to (\ref{43reg}). 
Therefore we can obtain by the same way that  
\beq
&&\Delta \psi(t) \leq C(d)(1+\frac{1}{t})\label{4D1},\\
&&\Delta \phi(T-t) \geq -C(d)(1+\frac{1}{t})\label{4D2}.
\enq
This gives the following uniform bound
\be
\|\Delta \phi\|_{\Linf([\tau,T-\tau]\times\Td, \ d\rho)}\leq C(d)(1+\frac{1}{\tau}).
\en
Here the constant is universal, and it is only supposed that the solution exists from $t=0$ to $t=T$. 
This surprising (in the sense that it is not true for viscosity solutions) result comes from the fact that the transformation $(\phi,p)\rightarrow (-\phi,p)(T-t)$ does not necessarily transform a viscosity solution of (\ref{41reg}) into another viscosity solution but transforms a variational solution into another variational solution. Actually it will be proved that one can choose the variational solution to be a viscosity solution in one time direction, but it may only be a subsolution when reversing the time. However, it will co\"incide $d\rho$ a.e. with the viscosity going from $T$ to 0. This means that we have a reversibility property $d\rho$ a.e.. 

Moreover $\Delta \phi$ is the divergence of the velocity field, and we have
$\Delta \phi =-\Dt\log \rho$.
From (\ref{4ineq}) we have the following control on the second time derivative of $\rho$ along the flow:
\beq\label{4ineq2}
\Dtt \log \rho \geq \rho-1 + \frac{1}{d}|\Dt \log\rho|^2.  
\enq
Following the path of a ``particle'', the differential inequality satisfied by $\Theta=\log \rho$ is:  
\be
\ddot \Theta \geq \exp{\Theta} -1   + \frac{1}{d}|\dot\Theta|^2.
\en
We look for solutions of this equation that do not become infinite in $]0,T[$.
This condition implies that $\Theta(t)\leq C(\tau)$ for $\tau\leq t \leq T-\tau$, independently of the initial and final values of $\Theta$. 
We have thus  an interior unconditional bound for the $\Linf$ norm of $\rho$, namely that
 \be
\|\rho\|_{\Linf([\tau,T-\tau]\times\Td)}\leq C(\tau),  \  \ \forall\, 0<\tau<T/2.
\en
The above differential inequality will also yield  that some functionals of $\rho$ are convex along the displacement
induced by our variational problem:
indeed a formal computation gives the following:
\be
&&\frac{d^2}{dt^2} \int_{\Td}\rho\log \rho\geq 0,\\
&&\frac{d^2}{dt^2} \int_{\Td}|\rho|^k \geq 0\textrm{ for every } k\geq 1.
\en

{\it Remark.} This displacement convexity property is analog to the one found in \cite{Mc1} which was true for $k\geq 1-1/d$ (when multiplied by $k-1$).
In our case this is only valid down to $k\geq 1$ due to the gravitational term.
However,  note that this displacement does not induce a distance:
indeed take $\rho_0=\rho_T \not\equiv 1$ and check that the cost of the transportation of $\rho_0$ on $\rho_T$ following the Euler-Poisson flow is not 0.

In the next section we give a rigorous sense to the computations made above in order to obtain the Theorem \ref{4main2}.


\subsection{Rigorous proof of Theorem \ref{4main2}}
\paragraph{Outline of the proof}
The rigorous justification will be achieved in several steps.
We will first introduce a time discretization of the problem, where the potential energy term contributes only at the time $t_i=iT/N, i=1..N-1$. Between two time steps, the problem will be an optimal transport problem as in \cite{BB1}, \cite{Br1} and \cite{Mc1}. Then at each time step, the gravitational effect will be taken into account, and the velocities will 
have a "jump". From a Lagrangian point of view, the velocity of each particle will therefore be a piecewise constant function with respect to time.
This discrete problem will accept also a variational formulation. Then letting the time step go to 0, we will eventually recover the time continuous problem.
One of the advantages of this discrete formulation is that between two time steps, the velocity potential will be expressed
with convex functions, which offer several regularity properties, and that will allow us to deal rigorously with quantities such as the second derivatives of $\phi$. Moreover, we will show that we can choose the optimal $\phi$ to be the viscosity solution of the Hamilton-Jacobi equation HJ2$(\phi,p)=0$. Then the Hopf-Lax representation formula will be a precious tool.

This section is organized as follows:  we first introduce the discrete problem, and recall some facts about optimal transportation({\it i.e.} the non-interacting transport). Then we show that solutions of this discrete problem converge indeed to the time-continuous one when the time step goes to 0. The rest is devoted to the proof of the regularity properties. For this we show that we can choose a special solution of the discrete problem that will be a viscosity solution.
For this solution, we are able to perform rigorously our computations and conclude.

\subsubsection{Construction of a sequence of approximate solutions}
We introduce the discrete times $t_i= Ti/N, \ i= 1..N-1$,
and consider the functional
\be
I_N(\rho,v,p)=\demi\int_D d\rho(t,x)|v(t,x)|^2  +\frac{T}{2N}\sum_{i=1}^{N-1}\int_{\Td}|\nabla p(t_i,x)|^2 \ dx.
\en
We are now interested in solving the following variational problem:

\begin{problem}\label{4infIN}
Minimize 
\be
&&\tilde I_N(\rho,J,p)=\frac{T}{2N}\sum_{i=1}^{N-1}\int_{\Td}|\nabla p(t_i,x)|^2 \ dx\\
&&+\sup_{\begin{array}{c}\scriptstyle c,m \in C^0(D)\times (C^0(D))^d\\\scriptstyle c+|m|^2/2\leq 0\end{array}}\left\{\int_{D} c(t,x) d\rho(t,x)   + m(t,x)\cdot dJ(t,x)  \right\}
\en 
among all $(\rho,J,p)$ that satisfy
$\rho \in C([0,T];{\cal P}(\Td)-w*)$, $J\in ({\cal M}(D))^d$, $\nabla p(t_i) \in L^2(\Td)$ for all $1\leq i\leq N-1$, and 
\be
&&\dt\rho+\nabla\cdot J = 0,\\
&&\Delta p=\rho-1,\\
&&\rho(t=0)=\rho_0,\\
&&\rho(t=T)=\rho_T.
\en
We denote $K_N$ the value of this infimum.
\end{problem}
The interest of studying this problem is both to make rigorous the arguments of the previous section and to give a possible  numerical discretization of the Problem \ref{4infI1}.
It will also let appear some interesting links between optimal transportation, viscosity solutions of Hamilton-Jacobi equations, and transport equations.
\subsubsection{Basic facts on optimal transportation}\label{4optitrans}
We first recall the definition of the push-forward of a measure by a mapping:
\begin{defi}\label{4pushforward}
Let $\rho_0$ and $\rho_1$ be two probability measures on $\Td$ and let $X$ be a $d\rho_0$ measurable mapping from $\Td$ into itself. We say that $\rho_1$ is the push-forward of $\rho_0$ by $X$, that we denote by 
$\rho_1=X_{\#} \rho_0$, if the following holds:
\be
\forall f \in C^0(\Td), \int f(X(x))d \rho_0(x)= \int f(y) d \rho_1(y).
\en 
\end{defi} 
We denote here $D_i=[t_i, t_{i+1}]\times \Td$; the effect of the time discretization is that between two times $t_i$, once the measures $\rho(t_i)$ and $\rho(t_{i+1})$ have been chosen,  the problem becomes the following:
\begin{problem}\label{4monge-kanto}
Minimize 
\be
\tilde C(\rho,J)=\sup_{\begin{array}{c}\scriptstyle c,m \in C^0(D_i)\times (C^0(D_i))^d\\\scriptstyle c+|m|^2/2\leq 0\end{array}}\{\frac{1}{2}\int_{D_i} d\rho c + dJ\cdot m \}
\en among all $(\rho,J)$ that satisfy
$\rho \in C([t_i,t_{i+1}];{\cal P}(\Td)-w*)$, $J\in ({\cal M}(D_i))^d$, and 
\be
&&\dt\rho+\nabla\cdot J = 0,\label{4mk11}\\
&&\rho(t=t_i)=\rho_i,\label{4mk2}\\
&&\rho(t=t_{i+1})=\rho_{i+1}.\label{4mk3}
\en
The infimum is denoted $\underline C(\rho_i, \rho_{i+1}, |t_i-t_{i+1}|)$.
\end{problem}
{\it Remark.} Performing  a dilatation in the time variable we see that $\underline C(\rho_i, \rho_{i+1}, t)=\frac{1}{t}\underline C(\rho_i, \rho_{i+1}, 1)$.
The Wasserstein distance (of order 2) between $\rho_i$ and $\rho_{i+1}$, denoted $W_2(\rho_i,\rho_{i+1})$, is given by 
\be
\left[W_2(\rho_i,\rho_{i+1})\right]^2=\underline C(\rho_i, \rho_{i+1}, 1).
\en
This problem has been solved in \cite{BB1}, \cite{Br1} where it is shown that there exists a unique solution \linebreak  $(\rho,J= \rho v)$ ($v$ is only unique $d \rho$ a.e.)  that satisfies:
\be
&&v(t=t_i,x)=\frac{1}{t_{i+1}-t_i}(\nabla\varphi(x)-x) \ d\rho_i \ a.e.\label{4W2,1},\\
&& \dt (\rho v) + \nabla\cdot (\rho v \otimes v) =0 \nonumber,\\
&&\det D^2 \varphi(x) \rho_{i+1}(\nabla\varphi(x))=\rho_i(x)\nonumber,
\en
with $\varphi$ a convex function. 
The second equation means that the particle move with constant speed. The third equation
is the Monge-Amp\`ere equation that is satisfied in the following weak sense:
\be
\forall f \in C^0(\Td), \int f(\nabla\varphi(x))d \rho_i(x)= \int f(y) d \rho_{i+1}(y).
\en 
This means that $\nabla\varphi$ pushes  $\rho_i$ forward to $\rho_{i+1}$. 
The Wasserstein distance between two probability measures can be defined equivalently in the following ways:
\be
W^2_2(\rho_0,\rho_1)&=&\inf_{\rho, v }\int_{[0,1]\times\Td}  d\rho|v|^2/2\\
&=& \sup_{\Phi(x)+\Psi(y)\geq x\cdot y} \int_{\Td}d\rho_0(x)(|x|^2/2-\Phi(x)) + d\rho_1(y)(|y|^2/2-\Psi(y))\\
&=&\inf_{{\bf m}_{\#} \rho_0=\rho_1}\int_{\Td}\demi|x-{\bf m}(x)|^2 d\rho_0\\
&=& \sup_{\dt\phi +|\nabla\phi|^2/2 \leq 0} \int d\rho_1(x)\phi(1,x) -d\rho_0(x)\phi(0,x). 
\en
where the first infimum is taken over all the pairs $(\rho,v)$ satisfying 
\be
&&\dt\rho +\nabla\cdot [\rho v]=0,\\
&&\rho|_{t=0}=\rho_0,\rho|_{t=1}=\rho_1.
\en
Note that the last formulation is strongly reminiscent of the formulation of Proposition \ref{4propdual}.
Under some assumptions of absolute continuity of $\rho_0, \rho_1$ with respect to the Lebesgue measure (see \cite{Vi} for a complete reference), those three problems have a unique solution. For the first it has already been mentioned above. Then the optimal ${\bf m}$ is equal to $\nabla\varphi$, the optimal pair $(\Phi, \Psi)$ is equal (up to a constant) to $(\varphi, \varphi^*)$ with $\varphi^*$ the Legendre transform of $\varphi$ (see
the definition (\ref{4legendre}) below), and $\phi(0)=\Phi-|x|^2/2, \phi(1)=|x|^2/2-\Psi$.  For additional references about the Wasserstein distance
the reader can also refer to \cite{Mc1} and  \cite{O}.

\subsubsection{Existence of a minimizer for the approximate problem}
Following exactly the same method as in the first problem we can  show the existence of a unique minimizer to the Problem \ref{4infIN}. In this way we obtain the following:
\begin{prop}\label{4MN}
There exists a unique $\rho_N$ and a $d \rho_N$ a.e. unique $v_N=\nabla\phi_N$ solution of Problem \ref{4infIN}.
Moreover it satisfies:

\begin{enumerate}
\item  There exists $C$ such that for any $0<\tau<T/2$,
 $\ds \frac{T}{N}\sum_{\tau\leq t_i\leq T-\tau}\|\rho_N(t_i,.)\|_{L^2(\Td)} \leq \frac{C}{\tau} $.


\item   The solution $(\rho_N,v_N=\nabla\phi_N)$ is a weak solution of
\be
&&\dt(\rho_N v_N) +\nabla\cdot (\rho_N v_N\otimes v_N) = - \rho_N \frac{T}{N}\sum_{i=1}^{N-1}\delta_{t=t_i}\nabla p_N (t_i),\\
&&\dt\rho_N +\nabla\cdot (\rho_N v_N)=0,\\
&&\Delta p_N= \rho_N -1.
\en
where the product $\rho_N\nabla p_N$ is defined as in \ref{4defiweakep}.
The pair $(\phi_N, p_N)$ satisfies 
\be
\dt\phi_N + \demi|\nabla\phi_N|^2 + \frac{T}{N}\sum_{i=1}^{N-1}\delta_{t=t_i} p_N \leq 0,
\en
and also satisfies for any $t_i, 1\leq i \leq N-1$ 
\be
&&\phi_N(t_i^+,x)-\phi_N(t_i^-,x) \leq -\frac{T}{N} p_N(t_i,x) \ dx  \ a.e.,\\
&&\phi_N(t_i^+,x)-\phi_N(t_i^-,x)=-\frac{T}{N} p_N(t_i,x) \ d\rho_N(t_i) \  a.e..
\en

\end{enumerate}
\end{prop}
 
{\it Proof.} The proof is the same as the time continuous version therefore we will only sketch it briefly.

Concerning the existence of an admissible solutions for Problem \ref{4infIN} note that we don't need that either $\rho_0$ or $\rho_T$ is in any $L^p$ since two probability measures on $\Td$ are always at finite Wasserstein distance and thus
one can exhibit an admissible solution by transporting $\rho_0$ on $\rho(T/N)= 1$ between $t=0$ and $t=T/N$ , then letting $\rho(t_i)=1$ for $i\leq N-1$ and transporting $\rho(\frac{N-1}{N}T)$ on $\rho_T$. Note also that the solution of the continuous problem is admissible for the discrete problem.

Then having chosen an admissible solution $(\bar \rho, \bar J =\bar \rho \bar v, \bar p)$,  
the problem $\tilde I_N$ admits a dual problem:
\begin{prop}\label{4propdualN}
Let $K_N$ be the infimum of Problem \ref{4infIN}. Let
\be
D_N(\psi,q) &=& \int_{[0,T]\times \Td} d\bar\rho\dt\psi + d\bar J\cdot  \nabla\psi \\ 
&+&\frac{T}{N}\sum_{i=1}^{N-1}\int_{\Td} d\bar\rho(t_i)q(t_i) + (\nabla \bar p(t_i)\cdot\nabla q(t_i) - |\nabla q(t_i)|^2/2) dx.
\en
Then  $$K_n = \inf_{\rho,J,p} I_N =  \sup_{\psi,q} D_N,$$
where  the supremum is taken over all pairs $(\psi,q)\in C^1(D)$ such that 
$$\dt \psi +\demi |\nabla\psi|^2 +\frac{T}{N}\sum_{i=1}^{N-1}\delta_{t=t_i}q(t_i)\leq 0,$$ moreover the infimum is attained.
\end{prop}

{\it Remark.} The dual functional $D_N$ can be rewritten as
\be
D_N(\psi,q) &=& \int_{\Td} d\rho(T)\psi(T) - d\rho(0)\psi(0) +\frac{T}{N}\sum_{i=1}^{N-1}\int_{\Td} q(t_i)  - \demi |\nabla q(t_i)|^2 dx.
\en

The proof of this proposition uses the Fenchel-Rockafellar duality Theorem as the time continuous one. The optimal velocity will here also be potential.
Then taking $(\bar \rho, \bar v, \bar p)=(\rho_N, \nabla\phi_N, p_N) $ the optimal solution, 
 for any $\epsilon>0$ we find  $\psi\epu, q\epu$ such that $D_N(\psi\epu, q\epu)\geq K_N - \epsilon^2$
and we obtain:
\beq\label{4epsilonN}
&&\demi \int_D d\rho_N|\nabla\phi_N-\nabla\psi\epu|^2 + \demi \frac{T}{N}\sum_{i=1}^{N-1}\int_{\Td} |\nabla p_N(t_i)-\nabla q\epu(t_i)|^2 \ dx \\
&+& \int_D d\rho_N\bigl|\dt \psi\epu +\demi |\nabla\psi\epu|^2 +\frac{T}{N}\sum_{i=1}^{N-1}\delta_{t=t_i}q\epu(t_i)\bigr| \leq \epsilon^2.\nonumber
\enq
Then perturbing $(\rho_N, v_N)$ as we did in the time continuous problem, we find the optimality equation and the regularity properties.

Now we prove the last 2 points of Proposition \ref{4MN}, that link $\phi(t_i^-), \phi(t_i^+), p(t_i)$. 
Before this we point out that the $\Linf_{loc}(]0,T[; L^2(\Td))$ bound for $\phi$ established in Remark 2, paragraph \ref{4rho_linf}, is still valid for $\phi_N$ with minor adaptations.

Then, in every $]t_i, t_{i+1}[$, $\dt\phi_N \leq -|\nabla\phi_N|^2/2$ is  a negative measure, in the sense of ${\cal D}'(]t_i,t_{i+1}[\times\Td)$.
 In particular, for a.e. $x$, $\phi_N$ is decreasing with respect to time. It follows from the monotone convergence theorem that $t \to \phi_N(t,\cdot)$ has left and right (strong) limits in $L^2(\Td)$ at every  $t\in ]0,T[$. This gives sense to $\phi(t_i^-), \phi(t_i^+)$. 

Since $\phi_N$ is a subsolution to $
\dt\phi+ |\nabla\phi|^2/2 + \delta_{t_i}p(t_i)T/N=0$ on $]t_{i-1}, t_{i+1}[$, 
we have immediately that $\phi(t_i^+)- \phi(t_i^-)\leq -p(t_i)T/N$.

Then we take the maximizing sequence $\psi\epu$, and $\varphi$ a smooth function on $\Td$, for which we have
\be
\int_{t_i-\delta}^{t_i+\delta}\int_{\Td} \rho\partial_t\psi\epu \varphi
&=& \int_{t_i-\delta}^{t_i+\delta}\int_{\Td} -\rho\demi|\nabla\psi\epu|^2\varphi 
- \frac{T}{N}\int_{\Td} p(t_i) \varphi + O(\epsilon^2)\\
&=&\int_{\Td} [\rho_N(t_i+\delta) \psi\epu(t_i+\delta) - \rho_N(t_i-\delta) \psi\epu(t_i-\delta)]\varphi\\
&& - \int_{t_i-\delta}^{t_i+\delta}\int_{\Td} \rho_N v_N\cdot \nabla(\psi\epu\varphi). 
\en
For the first equality we have used equation (\ref{4epsilonN}) and we have used the mass conservation equation for the second one.

The functions $\psi\epu(t,x)$ are decreasing with respect to time in every $]t_i, t_{i+1}[$,  and converging in $L^2_{loc}(D)-w$ to $\phi_N$, therefore $\psi\epu(t) \to \phi_N(t)$ in $L^2 -w$ for almost every $t\in ]t_i, t_{i+1}[$ (actually at every $t$ where where $\phi_N$ is weakly continuous in time, thus almost everywhere, since $\phi_N$ is decreasing).

 Hence $\psi\epu(t\pm\delta) \to \phi_N(t\pm\delta)$ weakly for almost every $\delta$.
As $\delta\to 0$, $\rho(t_i \pm \delta)$ converges weakly to $\rho(t_i)$ in $L^2$: indeed, note that $\rho_N(t_i)\in L^2$ for $1\leq 1\leq N-1$, and between two $t_i$, the problem co\"incides with the optimal transportation; then from Proposition \ref{mccann}, we get $\rho_N\in C(]0,T[; L^2(\Td)-w)$).

Then, since $\phi_N$ is decreasing with respect to time and bounded in $L^2$, we get that $\phi_N(t_i \pm \delta)$ converges strongly to $\phi_N(t_i^\pm)$ in $L^2$ as $\delta$ goes to 0.
The other integrals go to $0$ when $\delta$ goes to 0, except the one involving $p(t_i)$,  and we obtain
\be
\phi_N(t_i^+)-\phi_N(t_i^-) = -\frac{T}{N}p_N(t_i) \ d\rho(t_i) a.e..
\en

\subsubsection{Regularity properties of solutions of the discretized Problem \ref{4infIN}}
Here we state the main result of this section, from which  Theorem \ref{4main2} will be deduced.
\begin{prop}\label{4regulFN}
Let $(t,x)\rightarrow (\phi_N ,J_n=\rho_N \nabla\phi_N, p_N)$ 
be solution of Problem \ref{4infIN}.

\begin{enumerate}
\item There exists $C$ depending only on $T$ and on the dimension such that for all $t$ in $]0,T[$, $\phi_N(t)$ is $d \rho_N(t)$ a.e. twice differentiable and satisfies 
\be
&& -C(1+\frac{1}{T-t}) \leq \Delta \phi_N(t,.) \leq C(1+\frac{1}{t})\; d\rho_N(t)\ a.e..
\en

\item The density $\rho_N$ is  bounded in $\Linf_{loc}(]0,T[\times \Td)$ uniformly with respect to $N$ and belongs to $\bigcap_{k>1}C(]0,T[; L^k(\Td))$.

\item  There exists $C$ such that for any $1\leq k \leq \infty $
\be
-C(1+ \frac{1}{t}) \leq \Dt \log\left(\|\rho_N(t,\cdot)\|_{L^k(\Td)}\right)\leq C(1+ \frac{1}{T-t}).
\en

\item The functions $\displaystyle\int_{\Td}[\rho_N(t,x)]^k \ dx , \int_{\Td} \rho_N\log\rho_N(t,x) \ dx$ are uniformly Lipschitz with respect to time in every interval $[\tau, T-\tau]$ for $\tau\in ]0,T/2]$, and converge as $N\to \infty$ to convex functions on $[0,T]$. 

\item The potential $\phi_N$ can be chosen such that for every $0<\tau\leq T/2$,
  $$\|\nabla\phi_N\|_{\Linf([\tau,T-\tau]\times\Td)}\leq C(\tau).$$

\item One can also choose $\phi_N$  to be the viscosity solution of 
$\dt \phi_N + \demi|\nabla\phi_N|^2 + \frac{T}{N}\sum_{i=1}^{N-1}\delta_{t=t_i}p_N =0$ in the sense of (\ref{4defvisco}).

\item All these results and bounds do not depend on $\rho_0$ neither on $\rho_T$, and are uniform with respect to $N$.

\end{enumerate}

\end{prop}
The proof of this proposition is postponed to the end of the paper. First we use it to show  
the convergence of solutions of  Problem \ref{4infIN} toward the solution of Problem \ref{4infI1}:

\subsubsection{Convergence of the solutions of the discretized Problem \ref{4infIN} to the solution of the continuous Problem \ref{4infI2}}

\begin{prop}\label{4convergence}
Let $(\rho_N,v_N, p_N)$ be as above, and  $(\rho,v,p)$ be solution of the minimization Problem \ref{4infI1},
with same initial and final densities in $L^{\frac{2d}{d+2}}$, then
\be
\lim_{N\rightarrow \infty}\int_{t=0}^T\int_{\Td} \left(\rho|v_N-v|^2 + |\nabla p_N -\nabla p|^2\right)(t,x)\ dxdt\ =0,
\en
moreover $\rho_N\nabla\phi_N$ converges strongly in $L^1_{loc}(]0,T[\times\Td)$ to $\rho\nabla\phi$.
\end{prop}
Then the last two propositions combined will yield the Theorem \ref{4main2} when passing to the limit.

\bigskip
\noindent
{\bf Proof of Proposition \ref{4convergence}}\\
Here we prove a slightly weaker version of Proposition \ref{4convergence} that allows us to pass to the limit in Proposition \ref{4regulFN} and obtain Theorem \ref{4main2}. Then we can get the 
full Proposition \ref{4convergence}.
We first choose $\tau \in ]0,T/2[$ and $k=k(\tau, N)$ such that $0<t_{k-1}\leq \tau < t_k < ..<t_{N-k} \leq T-\tau < t_{N-k+1}$. 
We set 
\be
&&  F^{\tau}(\rho,v)=\demi \int_{0}^{T}\int_{\Td}d\rho |v|^2  + \demi\int_{\tau}^{T-\tau}\int_{\Td}|\nabla p|^2 \ dxdt,\\
&&F_N^{\tau}(\rho,v)=\demi \int_{0}^{T}\int_{\Td}d\rho |v|^2 + \frac{T}{2N}\sum_{i=k}^{N-k}\int_{\Td}
|\nabla p(t_i,x)|^2 dx,
\en
understood that $p$ satisfies $\Delta p =\rho-1$.
We need to introduce those truncated functionals since we don't know a priori that the potential energy remains bounded near the boundary of the time interval and thus the convergence of the Riemann sum to the integral is not clear. We shall see after having proved the Theorem \ref{4main2} that this is the case when $\rho_0$ and $\rho_T$ are in $L^{\frac{2d}{d+2}}$.

It follows from Proposition \ref{4MN} that there exists a unique minimizer \mbox{$(\rho_N, v_N=\nabla\phi_N, p_N)$} for the functional $F_N^{\tau}$ under the constraints of Problem \ref{4infIN}. It can also be checked later in the proof that the regularity results of Proposition \ref{4regulFN} remain uniformly valid for $\tau\leq \tau_0$.
We consider $(\rho, v)$ solution of Problem \ref{4infI2}. From Lemma \ref{4Ep},  $\|\nabla p(t,\cdot)\|_{L^2(\Td)}^2$ is continuous in $]0,T[$ thus 
$$\forall \tau >0, \ F_N^\tau(\rho,v)\rightarrow F^\tau(\rho,v)\textrm{ as } N \to \infty.$$
Then $F^\tau(\rho,v)-F(\rho,v)\rightarrow 0$ when $\tau$ goes to 0,
and there exists a sequence $(\tau_N)_{N\in \N^*}$ with $\tau_N \downarrow 0$ such that
$$F_{N}^{\tau_N}(\rho,v) \rightarrow F(\rho,v).$$
In the remainder of this proof we set $F_N:=F_{N}^{\tau_N}$ (therefore $\tau,k$ will both depend on $N$) and $(\rho_N, v_N)$ the minimizer of $F_{N}^{\tau_N}$.
Thus
\be
F_N(\rho,v)\rightarrow F(\rho,v).
\en

The idea of the proof is to show that $F(\rho_N, v_N) \to F(\rho,v)$, and to use the coercivity of $F$ to deduce that
$(\rho_N, v_N ) \to (\rho,v)$. The first point will need some regularity property for $(\rho_N,v_N)$.

We have first $F_N(\rho_N,v_N)\leq F_N(\rho,v)$ and $\limsup F_N(\rho_N,v_N)\leq F(\rho,v)$.

We claim that we also have $\lim|F_N(\rho_N,v_N)-F(\rho_N,v_N)|=0$,  this will imply  $$\lim F_N(\rho_N,v_N) = \lim F(\rho_N, v_N)= F(\rho,v)=\lim F_N(\rho,v).$$ 
Indeed, since $\rho_0$ and $\rho_T$ are in $L^{\frac{2d}{d+2}}$, point 3 in Proposition \ref{4regulFN} implies that $\rho_N$ is in $\Linf([0,T];L^{\frac{2d}{d+2}})$.
It follows that  the sequence $\nabla p_N$ is uniformly bounded in $\Linf([0,T];L^{2})$ since $\Delta p_N =\rho_N -1$ and using Gagliardo-Nirenberg inequality (cf. proof of Lemma \ref{4ilexiste}).

Then, we have $\dt \int |\nabla p_N|^2/2 = -\int \rho_N v_N \cdot\nabla p_N$. Using the interior $\Linf$ bounds on $\rho_N, v_N$ (Proposition \ref{4regulFN}, points 2 and 5),  $\int|\nabla p_N|^2$ is uniformly Lipschitz in compact sets of $]0,T[$. 
Hence we can conclude that
\beq 
\lim_{N\rightarrow \infty} \int_0^T\int_{\Td} |\nabla p_N(t,x)|^2 dtdx - \frac{T}{N}\sum_{i=k_N}^{N-k_N}\int_{\Td}|\nabla p_N(t_i,x)|^2 \ dx =0.\label{4intpn}
\enq
(Remember that since $\tau_n \downarrow 0$, we also have $k_N/N \sim \tau_N \downarrow 0$.)
We can conclude that \linebreak $\lim F_N(\rho_N,v_N) =\lim F_N(\rho,v) = F(\rho,v)$ as $N\to \infty$. 

We show now that this implies that $(\rho_N, v_N)$ and $(\rho,v)$ are close to each other.
For this we use the dual formulation of the problem (Proposition \ref{4propdualN}). This formulation needs an admissible solution $(\bar \rho, \bar J=\bar\rho\bar v, \bar p)$ satisfying (\ref{4c1}, \ref{4c2}, \ref{4c3}, \ref{4c4}). We take here  $(\rho,v,p)$ the optimal solution of Problem \ref{4infI2} which is admissible for Problem \ref{4infIN}.
Using the fact that $F_N(\rho,v)$ is close to $F_N(\rho_N,v_N)$ for $N$ large, 
for any $\epsilon,
\delta>0$ there exists  $N,\psi_N\ep,q_N\ep$ (with $N=N(\delta)$) such that 
\be
&&\demi\int_D\rho|v|^2 \ dtdx+\frac{T}{2N}\sum_{i=k}^{N-k}\int_{\Td}|\nabla p(t_i)|^2 \ dx\\
&\leq& \demi\int_D\rho_N|v_N|^2 \ dtdx+\frac{T}{2N}\sum_{i=k}^{N-k}\int_{\Td}|\nabla p_N(t_i)|^2 \ dx + \delta\\
&\leq& \epsilon + \delta + \int_D d\rho \dt \psi\ep_N + d\rho v\cdot \nabla\psi\ep_N\\
&+& \frac{T}{N}\sum_{i=k}^{N-k}\int_{\Td} \rho(t_i)q\ep_N(t_i) + \nabla p(t_i)\cdot\nabla q\ep_N(t_i) -|\nabla q\ep_N(t_i)|^2/2 \ dx. 
\en
This eventually yields
\be
\int_D\demi \rho |v-\nabla\psi\ep_N|^2 \ dtdx + \frac{T}{2N}\sum_{i=k}^{N-k} \ dx
\int_{\Td}\demi|\nabla p(t_i)-\nabla q\ep_N(t_i)|^2 \ dx 
\leq \epsilon+ \delta. 
\en
For fixed $N$ , $(\nabla\psi_N\ep, q_N\ep)$ being a maximizing sequence for the dual problem, it  will converge to $(v_N, p_N)$ as $\epsilon\to 0$ (see (\ref{4epsilonN})), therefore we obtain
\be
\demi \int_D \rho |v-v_N|^2 \ dtdx  + \frac{T}{2N}\sum_{i=k}^{N-k} \int_{\Td} |\nabla p(t_i)-\nabla p_N(t_i)|^2  \ dx \leq \delta,
\en
and therefore this goes to 0 as $N\to \infty$.
Using the same procedure we can also get that
\be
\demi \int_D  \rho_N |v-v_N|^2   + |\nabla p_N-\nabla p|^2 \ dtdx \rightarrow 0 
\textrm{ as } N\rightarrow \infty.
\en
Now we show that the product $\rho_Nv_N$ converges to $\rho v$: using the equicontinuity property of the sequence $\|\rho_N(t,.)\|_{L^k}$ in $[\tau, T-\tau]$ for any $\tau\leq T/2$ and any $1\leq k <\infty$, (see Proposition \ref{4regulFN}), the sequence $\rho_N$ converges strongly in $L^k([\tau, T-\tau]\times \Td)$ for any $1\leq k <\infty$. 
Moreover remember that from Theorem \ref{4main}, $v\in L^2([\tau, T-\tau]\times \Td)$.
Then 
\be
&&\int_{[\tau, T-\tau]\times \Td} |\rho v - \rho_N v_N|\\
&\leq& \int_{[\tau, T-\tau]\times \Td}\rho_N|v_N-v| + |v|\, |\rho_N -\rho| \to 0.
\en
Hence, $\rho_Nv_N$ converges strongly to $\rho v$ in $L^1_{loc}(D)$.
Using the uniform $\Linf_{loc}(D)$ bound on $\rho_N,v_N$ (points 2 and 5 of Proposition \ref{4regulFN}), we get that $\rho_N v_N$ converges strongly to $\rho v$ in $L^p_{loc}(]0,T[\times\Td)$ for any 
$1\leq p <\infty$. 

\bigskip
\noindent
{\bf Proof of Theorem \ref{4main2}}: The theorem is  obtained passing to the limit in the  Proposition \ref{4regulFN}. The point 2, 3, 4, 5 remain true when we let $N$ go to $\infty$.
The other points will be shown at the end of the paper.

$\hfill \Box$

Then from Theorem \ref{4main2} if $\rho_0, \rho_T \in L^{\frac{2d}{d+2}}(\Td)$ 
we have $\rho\in\Linf([0,T];L^{\frac{2d}{d+2}}(\Td))$ and doing as in Lemma \ref{4ilexiste},  $\nabla p \in \Linf([0,T];L^{2}(\Td))$. This bound 
shows that the Riemann sum \linebreak  $\frac{T}{N}\sum_{i=1}^{N-1}\int_{\Td}|\nabla p(t_i,x)|^2 dx$ converges to 
$\int_D |\nabla p|^2 \ dt dx$ and this
allows us to take $\tau = 0$ in the previous proof and to conclude the proof of Proposition \ref{4convergence}. 

$\hfill \Box$

\subsection{Regularity properties of the time discretized solution: Proof of Proposition \ref{4regulFN}}
In this part $N$ is fixed and for sake of simplicity we drop the subscript $N$. We consider 
\\
\mbox{$(\rho=\Delta p+1,v=\nabla\phi)$} solution of the Problem (\ref{4infIN}).

\subsubsection{Preliminary: Construction of a special solution}
First we begin to show the consistency with the optimal transport Problem \ref{4monge-kanto}:
let $t_i\leq s,t \leq t_{i+1}$, we denote 
\beq\label{4defbfi}
\bfi_{s,t}(x)=(t-s)\phi(s,x)+|x|^2/2.
\enq
The function $\bfi_{s,t}$ goes from $\Rd$ to $\R$ if we extend $\phi$ to a periodic function on all of $\Rd$. Note also that for any $\vec p \in \Zd$, $\nabla \bfi_{s,t}(.+\vec p)= \nabla\bfi_{s,t}(.)  + \vec p.$
If $s=t_i$ and $t=t_{i+1}$ we denote $\bfi_{i,i+1}$ (resp. $\bfi_{i+1,i}$) instead of
$\bfi_{s,t}$ (resp. $\bfi_{t,s}$).
Note that $\phi$ is discontinuous at times $t_i$ and thus 
\be
&&\bfi_{i,i+1}(x)=|x|^2/2 + \frac{T}{N}\phi(t_i^+, x),\\
&&\bfi_{i,i-1}(x)=|x|^2/2 - \frac{T}{N}\phi(t_i^-, x).
\en
We also have
\beq\label{4defv}
v(s,x)=\nabla\phi(s,x)=\frac{1}{t-s}(\nabla\bfi_{s,t}(x)-x),
\enq
which is well defined on $\Rd /\Zd$.

In the first lemma, we will see that $\nabla\bfi_{s,t}$ pushes forward $\rho(s)$ on $\rho(t)$ minimizing the cost
 $\int_{\Td}\rho(s,x)|{\bf  m}(x)-x|^2dx$ among all ${\bf m}$ pushing forward $\rho(s)$ on $\rho(t)$ (that we denote hereafter \linebreak ${\bf m}_{\#}\rho(s)=\rho(t)$) and that $\bfi_{s,t}$ coincides with its convex hull $d\rho(s)$ a.e.. 

Then in the second lemma we will show that we can consider a solution for which every $\bfi_{s,t}$ is convex. This point that may seem to be a direct consequence of optimal transport (the fact that the optimal transport is given by the gradient of a convex function) needs from our point of view a 
careful proof: indeed we only know that the gradient of $\bfi_{s,t}$ will coincide $d\rho(s)$ a.e. with the gradient of a convex function, but since the optimality equation links $\bfi_{i-1,i}$, $\bfi_{i,i+1}$ and $p(t_{i})$ it must be checked that $\bfi_{i,i+1}$ can consistently be taken convex.
The convexity will then allow us to consider the second derivative of $\bfi_{i,i+1}$ since a convex function is almost everywhere twice differentiable, and then to make rigorous the inequality (\ref{4ineq}) and its consequences.

First for any $f:\Rd \rightarrow \R$ we denote $f^*$ its Legendre transform defined by
\beq
f^*(y)=\sup_{x\in\Rd}\left\{y\cdot x -f(x)\right\}.\label{4legendre}
\enq
The convex hull of $f$ is $(f^*)^*$ (or $f^{**}$ in short). It is equivalently defined as the supremum of all convex functions smaller than $f$. We will show the following lemma:
\begin{lemme}\label{4**}
Let $t_i\leq s,t \leq t_{i+1}$ and $\bfi_{t,s},\bfi_{s,t}$ be
defined as above. Then 
\be
&& \bfi_{s,t}\geq \bfi_{t,s}^* \textrm{ with equality } \ d \rho(s) \ a.e.,\\
&& \bfi_{s,t}=\bfi_{s,t}^{**} \ d \rho(s) \ a.e.,\\
&&\nabla\bfi_{s,t}=\nabla \bfi_{s,t}^{**} \ d \rho(s) \ a.e.,\\
&&\nabla \bfi_{s,t \, \#}^{**}\rho(s)=\rho(t).
\en
\end{lemme}

{\it Proof.} We  know that $\phi$ is the limit of a smooth sequence 
$\phi\epu$ satisfying the the constraint
\be
\dt\phi\epu +\demi|\nabla\phi\epu|^2 + \frac{T}{N}\sum_{i=1}^N\delta_{t=t_i}p\epu(t_i)\leq 0.
\en
Between $t_i$ and $t_{i+1}$ we have
\beq\label{4qwe}
\dt\phi\epu+\demi|\nabla\phi\epu|^2 \leq 0.
\enq 
Thus if $t>s$ considering $\displaystyle\gamma(\sigma)= x+(\sigma - s)\frac{y-x}{t-s}$ and using (\ref{4qwe})
we find
\be
\frac{d}{d\sigma}[\phi\epu(\sigma, \gamma(\sigma))]
&=&\dt\phi\epu(\sigma, \gamma(\sigma)) + \frac{y-x}{t-s}\cdot \nabla\phi\epu(\sigma, \gamma(\sigma))\\
&\leq&  \dt \phi\epu(\sigma, \gamma(\sigma)) +\demi  \left(\frac{|y-x|^2}{|t-s|^2} + |\nabla\phi\epu(\sigma, \gamma(\sigma))|^2\right)\\
&\leq&  \demi \frac{|y-x|^2}{|t-s|^2}.
\en
Integrating from $s$ to $t$ we find
\be\label{4rty}
&&\phi\epu(t,y)\leq \inf_x \{\phi\epu(s,x)+\frac{|y-x|^2}{2(t-s)}\},\\
&& \phi\epu(s,x)\geq \sup_y \{\phi\epu(t,y) -\frac{|y-x|^2}{2(t-s)}\}.
\en
Letting $\epsilon$ go to 0, it follows from (\ref{4defbfi}) that 
\beq\label{4>*}
&&\bfi_{t,s}(y) \geq \sup_x \{x\cdot y - \bfi_{s,t}(x)\}=\bfi^*_{s,t}(y),\\
&& \bfi_{s,t}(x)\geq\bfi^*_{t,s}(x).
\enq
This is the first point of the lemma.
The crucial point  is the following: for $d\rho(t)$ a.e. $y$ we have
\be
&&\phi(t,y)=\inf_x \{\phi(s,x)+\frac{|y-x|^2}{2(t-s)}\}\label{4cl1},
\en
or equivalently
\be
\bfi_{t,s}(y)=\sup_x\{ x\cdot y - \bfi_{s,t}(x)\} \ d\rho(t) \ a.e.. \label{4cl2}
\en 
Indeed take a smooth sequence $(\phi_\epsilon,p\epu)$ such that 
\be
\dt\phi\epu +\demi|\nabla\phi\epu|^2 + \frac{T}{N}\sum_{i=1}^N\delta_{t=t_i}p\epu(t_i)\leq 0
\en
that
maximizes the dual problem, {\it i.e.} such that
\be
&&\int_0^T\int_{\Td} \rho \dt\phi_\epsilon{}+\rho\nabla\phi\cdot \nabla\phi_\epsilon{} \ dt dx\\
&+&\frac{T}{N}\sum_{i=1}^{N-1} \int_{\Td}\rho p_\epsilon(t_i,x) 
+\nabla p \cdot\nabla p_\epsilon(t_i,x)
-\demi|\nabla p_\epsilon{}|^2(t_i,x) \  dx\\
&\geq& K_N-\epsilon^2.
\en
In view of (\ref{4epsilonN}), being a maximizing sequence of the dual problem implies the following: 
\be
\int_D \rho|\dt\phi\epu + \demi|\nabla\phi\epu|^2 + \frac{1}{N}\sum_{i=1}^N\delta_{t=t_i}p\epu(t_i)| \ dtdx \to 0 \textrm{ as } \epsilon \to 0, 
\en 
which in turn implies that 
\be
&&\int_{s}^t\int_{\Td} \rho\dt \phi_\epsilon{} + \rho\nabla\phi\cdot\nabla\phi_\epsilon{} \ dt' dx\\
&=&\int_{\Td}\rho(t,x)\phi_\epsilon(t,x)-\rho(s,x)\phi_\epsilon(s,x) \ dx \\
&\rightarrow&  \demi\int_s^t \int_{\Td}\rho|\nabla\phi|^2 \ dt'dx 
=\frac{1}{t-s}W^2_2(\rho(s),\rho(t)) \textrm{ as } \epsilon \to 0.
\en
where the first line comes from the mass conservation equation (\ref{4continuite}) satisfied by the pair$(\rho, v=\nabla\phi)$ and the last identity
comes from the fact that between $t$ and $s$ the problem co\"incides with the optimal transport Problem \ref{4monge-kanto}.

If ${\mathbf m}$ is a mapping realizing the optimal transport of $\rho(s)$ onto $\rho(t)$  then 
${\mathbf m}_{\#}\rho(s)=\rho(t)$ implies
\be
\int_{\Td}\rho(t,x)\phi_\epsilon(t,x)-\rho(s,x)\phi_\epsilon(s,x) \ dx = \int_{\Td} \rho(s,x)(\phi\epu(t,{\mathbf m}(x))-\phi\epu(s,x)) \ dx.
\en 
Using that 
\be
\phi\epu(t,x)-\phi\epu(s,y)\leq \frac{|x-y|^2}{2(t-s)},
\en
and that from the optimality of ${\bf m}$ we have   
\be
\frac{1}{2(t-s)}\int_{\Td} \rho(s,x)|x-{\mathbf m}(x)|^2 \ dx=\frac{1}{t-s}W^2_2(\rho(s),\rho(t)),
\en
we obtain by taking the limit $\epsilon\to 0$ that 
\be
\phi(t,{\mathbf m}(x))=\phi(s,x)+\frac{|x-{\mathbf m}(x)|^2}{2(t-s)} \  d\rho(s) \ a.e., 
\en
which is equivalent to 
\be
\phi(t,y)=\phi(s,x)+\frac{|y-{\mathbf m}^{-1}(y)|^2}{2(t-s)} \ d\rho(t) \ a.e.. 
\en
We remind the reader that ${\mathbf m}$ is invertible $d \rho(t)$ a.e. and ${\mathbf m}^{-1}$ can be defined as the ($d\rho(t)$ a.e. unique) mapping realizing the optimal transport of $\rho(t)$ onto $\rho(s)$.
This implies that 
if $t>s$, we have: 
\be
&&\phi(t,x)=\inf_y \{ \frac{|y-x|^2}{2(t-s)}+\phi(s,y) \} \  d \rho(t) \ a.e.,\\
&&\bfi_{t,s}=(\bfi_{s,t})^* \  d \rho(t) \ a.e.,
\en
the two lines being equivalent through equation (\ref{4defbfi}).
As a supremum of affine functions, any Legendre transform  is convex. Here, $\bfi_{t,s}$ 
coincides with a convex function  $d \rho(t)$ a.e.  and is above this function $dx$ a.e.  from (\ref{4>*}). Since $\bfi_{t,s}^{**}$ is the convex hull of $\bfi_{t,s}$  it follows that 
\be
\bfi_{t,s}^{**}=\bfi_{t,s} \ d\rho(t) \ a.e.,
\en
from which it can be deduced that
\be
&&\int_{\Td}(|x|^2/2-\bfi_{t,s}^{**})\rho(s,x) \  dx+\int_{\Td} (|y|^2/2-\bfi_{t,s}^{*}(y) )\rho(t,y) \ dy\\
&=& (t-s)\int_{\Td} \rho(s,x)\phi(s,x)-\rho(t,x)\phi(t,x) \ dx\\
&=&W^2_2(\rho(s),\rho(t)).
\en
This implies that $\nabla\bfi_{s,t\,\#}^{**}\rho(s)=\rho(t)$ (see \cite{Br1}).
Note also that if we set 
$$\tilde \phi(s,x)= \frac{1}{t-s}\left[\left(|\cdot|^2/2+(t-s)\phi(s,\cdot)\right)^{**}(x)-|x|^2/2\right],$$ we obtain 
that
\be
\int_{\Td}\rho(t,x)|\nabla \phi - \nabla\tilde\phi|^2(t,x)dx=0
\en
for a.e. $t$. The proof of Lemma \ref{4**} is complete. 

$\hfill \Box$

\bigskip
\noindent
We are going to use the previous lemma to construct a new sequence of solutions for which
the potentials $\bfi_{i,i+1}$ are convex. This will allow us 
to define $d \rho(t_i)$ a.e. the second derivative of  $\bfi_{i,i+1}$.
This special solution will turn out to be the viscosity solution of $\dt \psi +\demi|\nabla \psi|^2+\displaystyle\frac{T}{N}\sum_{i=1}^{N-1}\delta_{t=t_i}p=0$.
Remember that from Proposition \ref{4MN}, $\phi(t_i^+)$ satisfies
\beq
&&\phi(t_i^+,x)-\phi(t_i^-,x) \leq -\frac{T}{N} p(t_i,x) \ dx \  a.e.\label{4ti2},\\
&&\phi(t_i^+,x)-\phi(t_i^-,x)=-\frac{T}{N} p(t_i,x) \ d\rho(t_i) \  a.e.\label{4ti1}.
\enq
Consider the new solution $\psi$ defined by
\beq
&&\psi(t=0,x)=\frac{N}{T}\left[\left(\frac{|\cdot|^2}{2}+\frac{T}{N}\phi(t=0,\cdot)\right)^{**}(x)-|x|^2/2\right]\label{4new1},\\
&&\textrm{ on }]t_i,t_{i+1}[, \ \psi(t,x)=\inf_y\left\{\frac{|x-y|^2}{2 (t-t_i)}+\psi(t_i^+,y)\right\}
\label{4new2},\\
&&\psi(t_i^+,x)=\frac{N}{T}\left[\left(\frac{|\cdot|^2}{2}+\frac{T}{N}\psi(t_i^-,\cdot)-\frac{T^2 p(t_i,\cdot)}{N^2}\right)^{**}(x)-|x|^2/2\right]\label{4new3}. 
\enq
Exchanging $\rho_T$ and $\rho_0$, $(\rho(T-t), p(T-t), -\phi(T-t))$ will be the corresponding optimal solution, and we introduce also $\tilde \psi$ constructed by the same procedure from $(p(T-t), -\phi(T-t))$.

\begin{lemme}\label{4=**}
\noindent
\begin{enumerate}
\item For almost every $t\in [0,T]$, $\psi(t,.) $ coincides with $\phi(t,.)$  $d \rho(t)$ a.e. and $(\rho, v=\nabla\psi, p)$ is solution of Problem \ref{4infIN},

\item $\forall i\in [0..N-1]$, $(t,x)\rightarrow\psi(t,x)$ and $(t,x)\rightarrow -\psi(t_{i+1}+t_i-t,x)$  are both viscosity solutions (and subsolutions) of $\dt \psi +\demi|\nabla \psi|^2=0$ on $[t_i, t_{i+1}]$,

\item $\psi$ is the viscosity solution of $\dt \psi +\demi|\nabla \psi|^2+\displaystyle\frac{T}{N}\sum_{i=1}^{N-1}\delta_{t=t_i}p=0$ on $[0,T]$ in the sense of (\ref{4defvisco}),

\item finally $-\tilde\psi(T-t)= \psi(t)$ holds $d\rho(t)$ a.e..

\end{enumerate}

\end{lemme}

{\it Proof.}
We denote 
\beq\label{4defbpsi}
\bpsi_{s,t}(x)=|x|^2/2+(t-s)\psi(s,x) \textrm{ for }s,t\in[t_{i-1},t_{i+1}],
\enq
and $\bpsi_{i,i+1}, \bpsi_{i+1,i}$ as well.
Let us now prove by induction the following:
\be
\textit{ For all } 1\leq i\leq N-1, \  \psi(t_{i}^-,x)\geq \phi(t_{i}^-,x) \textit{ with equality } d \rho(t_{i})\, a.e..
\en

Equation (\ref{4new1}) implies that $\bpsi_{0,1}= \bfi_{0,1}^{**}$. 
Then (\ref{4new2}) implies that $\bpsi_{1,0}=\bpsi_{0,1}^{*}= \bfi_{0,1}^{***}=\bfi_{0,1}^{*}\leq \bfi_{1,0}$ with equality $d \rho(t_1)$ a.e. from Lemma \ref{4**}.
The equality  $\bfi_{0,1}^{***}=\bfi_{0,1}^{*}$ comes from the fact that $\bfi_{0,1}^{*}$ is convex  as a Legendre transform and that the the Legendre transform is an involution on convex functions.
Thus from (\ref{4defbpsi}) we get that $\psi(t_1^-) \geq \phi(t_1^-)$ with equality $d \rho(t_1)$ a.e.. This proves the property for $i=1$.
 
Suppose now that $\psi(t_{i}^-)\geq \phi(t_{i}^-)$ with equality $d \rho(t_{i})$ a.e..
Then we have, using  (\ref{4ti2})
\be
|x|^2/2+\frac{T}{N}\psi(t_{i}^-,x)-\frac{T^2}{N^2}p(t_i,x)&\geq& |x|^2/2+\frac{T}{N}\phi(t_{i}^-,x)-\frac{T^2 }{N^2}p(t_i,x)\\
&\geq&|x|^2/2+\frac{T}{N}\phi(t_{i}^+,x), 
\en
where both inequalities are equalities $d\rho(t_i)$ a.e., 
hence $\bpsi_{i,i+1}\geq \bfi_{i,i+1}^{**}$ from  (\ref{4new3}).
We now show that
$$\bpsi_{i,i+1}= \bfi_{i,i+1}^{**}  \ d \rho(t_i)\, a.e..$$
Indeed from Lemma (\ref{4**}), we know that $\bfi_{i,i+1}=\bfi_{i,i+1}^{**}$ $d \rho(t_i)$ a.e., therefore the convex function $\bfi_{i,i+1}^{**}$ is below $|x|^2/2+\frac{T}{N}\psi(t_{i}^-,x)-\frac{T^2}{N^2}p(t_i,x)$, with equality $d\rho(t_i)$ a.e.. We conclude since  $\bpsi_{i,i+1}$ is pinched between
$|x|^2/2+\frac{T}{N}\psi(t_{i}^-,x)-\frac{T^2}{N^2}p(t_i,x)$ and $\bfi_{i,i+1}^{**}$ that co\"incide $d\rho(t_i)$ a.e..  

Then 
\be
\bpsi_{i+1,i}=\bpsi_{i,i+1}^{*}\leq \bfi_{i,i+1}^{***}=\bfi_{i,i+1}^{*}\leq \bfi_{i+1,i},
\en
where we have used that the Legendre transform is decreasing with respect to functions, is  an involution on convex functions and from  Lemma \ref{4**} for the last inequality.
This implies $\psi(t_{i+1}^-)\geq \phi(t_{i+1}^-)$.

We have 
\be
&&\bpsi_{i,i+1}=\bfi^{**}_{i,i+1}  \ d\rho(t_i)\,  a.e.,\\
&&\nabla\bpsi_{i,i+1}=\nabla\bfi_{i,i+1}^{**}\ d\rho(t_i)\,  a.e.,\\
&&\nabla\bpsi_{i,i+1\,\#} \rho(t_i)= \nabla\bfi^{**}_{i,i+1\,\#} \rho(t_i)= \rho(t_{i+1}).
\en
This implies that $\bpsi_{i,i+1}^*=(\bfi^{**}_{i,i+1})^*$ $d\rho(t_{i+1})$ a.e.: 
for this use the identity 
\be
c(x) + c^*(\nabla c(x))= x\cdot\nabla c(x)
\en
that holds for any Lipschitz convex function $c$. This will also imply $\bpsi_{t_i,t}^*=(\bfi^{**}_{t_i,t})^*$ $d\rho(t)$ a.e. for all $t\in [t_i,t_{i+1}]$.  
Hence we have proved that $\phi(t_{i+1}^-)=\psi(t_{i+1}^-)$   $ d\rho(t_{i+1})$ a.e..
In view of Lemma \ref{4**}, we also have for all $t\in [0,T]$, $\psi(t)=\phi(t)$ $d\rho(t)$ a.e.. This proves the first point of the lemma.

On $]t_i,t_{i+1}[$ we have $\dt\psi + \demi|\nabla\psi|^2=0$ satisfied in the viscosity sense from the Hopf-Lax formula (\ref{4new2}). 
To see that $-\psi(t_{i+1}+ t_i -t,.)$ is also a viscosity solution to this equation, 
first note that
${\bpsi}_{i,i+1}=|x|^2/2 + \frac{T}{N}\psi(t_i,x)$ is convex.
Then we have 
\be
\psi(t,x)&=&\inf_y\{\psi(t_i,y)+\frac{|x-y|^2}{2(t-t_i)}\} \\
&=&\inf_y\{\psi(t_i,y)+\frac{|x-y|^2}{2(t_{i+1}-t_i)} +\gamma\frac{|x-y|^2}{2(t_{i+1}-t)}  \}
\en
for some $\gamma>0$. Now observe that $\alpha(y)= \psi(t_i,y)+\frac{|x-y|^2}{2(t_{i+1}-t_i)} $ is convex as well as
$\beta(y)=\gamma\frac{|x-y|^2}{2(t_{i+1}-t)}$ and thus we can apply the Fenchel-Rockafellar duality Theorem that says that if $\alpha, \beta$ are convex continuous we have
\be
\inf_y\{\alpha(y)+\beta(y)\}=\sup_x\{-\alpha^*(x)-\beta^*(-x)\}.
\en
Computing $\alpha^*$ gives $\frac{1}{t_{i+1}-t}\bpsi_{i+1,i}$ and we can check that
\be
\psi(t,x)=\sup_y\{\psi(t_{i+1}^{-},y)  - \frac{|x-y|^2}{2(t_{i+1}-t)}\},
\en
which says exactly that $t\rightarrow  -\psi(t_{i+1}+t_i-t)$ is the  viscosity solution of $\dt\psi + \demi|\nabla\psi|^2=0$ on $[t_{i}, t_{i+1}]$. 

Then we check that 
\be
\psi(t_{i+1}^+,x)&=& \frac{N}{T}[(|\cdot|^2/2+\frac{T}{N}\psi(t_{i+1}^-,\cdot)-\frac{T^2 p(t_i,\cdot)}{N^2})^{**}(x)-|x|^2/2] \\
&\leq & \psi(t_{i+1}^-,x)-\frac{T p(t_i,x)}{N}
\en
with equality $\rho_N(t_{i+1})$ a.e.,
which shows that $\psi$ is a subsolution of  
\beq\label{4HJdiscret}
\dt \psi +\demi|\nabla \psi|^2+\displaystyle\frac{T}{N}\sum_{i=1}^{N-1}\delta_{t=t_i}p=0.
\enq
To see that $\psi$ is the viscosity solution of (\ref{4HJdiscret}), 
just notice that our definition
of $\psi$ is the following:
\be
|x|^2/2-\frac{T}{N}\psi(t_{i+1}^-,x)&=&(|\cdot|^2/2+\frac{T}{N}\psi(t_{i}^+,\cdot))^*\\
&=&(|\cdot|^2/2+\frac{T}{N}\psi(t_{i}^-,\cdot)-\frac{T^2}{N^2}p(t_i,\cdot))^*,
\en
thus
\be
\psi(t_{i+1}^-,x)=\inf_{y}\{\psi(t_i^-,y)-\frac{T}{N}p(t_i,y)+\frac{|x-y|^2}{2\frac{T}{N}}\},
\en
and we obtain 
\be
&&\psi(t_{i+p}^-,x)=\\
&& \inf_{\gamma \in \Gamma} \left\{\psi(t_i^-,\gamma(t_i^-))+\int_{t_i}^{t_{i+p}^-}  \left[-  \frac{T}{N}\sum_{i=1}^{N-1}\delta_{\sigma=t_i}p(\sigma, \gamma(\sigma))+ \demi|\dot\gamma|^2(\sigma) \right]\,d \sigma \right\},
\en
where $\Gamma$ is the set of all continuous paths with $\gamma(t_{i+p})=x$,
and more generally that
\beq\label{4defvisco}
&&\psi(t,x)=\\
&& \inf_{\gamma(t)=x}\left\{\psi(s,\gamma(s))+\int_s^t  \left[-  \frac{T}{N}\sum_{i=1}^{N-1}\delta_{\sigma=t_i}p(\sigma, \gamma(\sigma))+ \demi|\dot\gamma|^2(\sigma) \right]\,d \sigma\right\}\nonumber.
\enq
This defines the viscosity solution of (\ref{4HJdiscret}), and proves the point 3 of the lemma.

The proof of the point 4 follows simply from the point 1: $\psi(t)$ and $\phi(t)$ co\"incide $d\rho(t)$ a.e., and $\tilde\psi$ is built from $-\phi(T-t)$ in the same way as $\psi$ is built from $\phi$, therefore $\tilde\psi(t)$ and $-\phi(T-t)$ co\"incide $d\rho(T-t)$ a.e., and we conclude. 
This achieves the proof of Lemma \ref{4=**}.

$\hfill \Box$

{\it Remark.} The time reversibility property is not valid for any viscosity solution. Actually this is true 
before occurrence of shocks. 
This is what one says when we decompose $\phi(0,y)+\frac{|x-y|^2}{2(t-t_i)}$ as the sum of two convex functions: this means that one can continue the rays  further without developing shocks.
Thus we see that our variational solution does not develop shocks in the interior of the time interval.  

\subsubsection{Proof of the bound on $\Delta \phi$}
The function $\psi$ satisfies
\be
&&\dt\psi +  \demi|\nabla\psi|^2 +\frac{T}{N}\sum\delta_{t=t_i}p(t_i)\leq 0,\\
&& \dt \rho + \nabla\cdot (\rho \nabla \psi)=0.
\en
From now we consider that $\phi:=\psi$ and thus $\bfi_{t,s}$ is convex for any $t_i\leq s,t \leq t_{i+1}$.
We are going to prove the following lemma:
\begin{lemme}\label{4d2}
The functions $\bfi_{i,i+1},\bfi_{i,i-1}$ are $C^{1,1}$  at every density 
point of $\rho(t_i)$ for $1\leq i \leq N-1$,  with $C^{1,1}$ norm bounded by $C(N,d)$.
Moreover there exists a set ${\cal E}_i$ of full measure for $\rho(t_i)$ such that everywhere in ${\cal E}_i$, $\phi(t_i^+)$ and $\phi(t_i^-)$ are twice 
differentiable  and  the following holds
\beq\label{4Delta_i}
\Delta\phi(t_i^+,x)-\Delta\phi(t_i^-,x)\leq \frac{T}{N}(1-\rho(t_i,x)).
\enq
\end{lemme} 

{\it Proof.}
Using (\ref{4ti1}, \ref{4ti2}) we get that 
\beq
 && \bfi_{i,i+1}(x) + \bfi_{i,i-1}(x)= |x|^2-\frac{T^2}{N^2}p(t_i,x) \ d \rho(t_i) \ a.e.,\label{4Deltati1}\\
&& \bfi_{i,i+1}(x) + \bfi_{i,i-1}(x) \leq |x|^2-\frac{T^2}{N^2}p(t_i,x) \ dx \ a.e.. \label{4Deltati2}
\enq 
Thus, $d \rho(t_i)$ almost everywhere, the convex function $\bfi_{i,i+1}+\bfi_{i,i-1}$ is  tangent from below to \\
\mbox{$|\cdot|^2-\frac{T^2}{N^2}p(t_i,\cdot)$}.
The Poisson equation satisfied by $p(t_i)$ being only true in the distribution sense we need to introduce a finite difference version of the Laplacian:

\begin{lemme}\label{4bolemme}
Let $c_d$ be the volume of the unit ball of $\Rd$ and $b_d$ be the $d-1$ dimensional Hausdorff measure of the unit sphere. Let $p$ be continuous and 
\be
&&\Delta_h p(x)=\frac{k_d}{h^2}\left[\frac{1}{b_d h^{d-1}} \int_{\partial B(x,h)}p(y)dy-p(x)\right],\\
&&\Delta_h^* p(x)=\frac{l_d}{h^2}\left[\frac{1}{c_d h^{d}} \int_{B(x,h)}p(y)dy-p(x)\right],
\en
with the constants $k_d, l_d$ chosen so that both operators converge to the Laplacian for smooth functions as $h\to 0$.
Then 
\begin{enumerate}
\item if $\Delta p \geq C$ in $\Td$ in the distribution sense then $\Delta_h p \geq C$   and $\Delta_h^* p \geq C$  everywhere in $\Td$,

\item if $\Delta p \in L^2(\Td)$ then, up to extraction of a subsequence in $h$, $\Delta_h^* p  \to \Delta p$   $dx$ almost everywhere,  
\end{enumerate}
finally, point 2 still holds even if $p$ is not continuous.
\end{lemme}

{\it Proof of claim 1.}  We take $f$ solution of 
\be
&&f-p|_{\partial B(x,h)}=0,\\
&&\Delta f=C \leq \Delta p.
\en
Observe that the boundary condition has a meaning since  $p$ is continuous.
Then from the maximum principle $f(x)\geq p(x)$ in $B(x,h)$, and thus $\Delta_h f\leq \Delta_h  p$ since $f=p$ on $\partial B(x,h)$. But \linebreak $f(x)=C |x|^2 /2d +g$ with $\Delta g=0$; using the fact that the average on a sphere of an harmonic function equals its value at the center of the sphere, we get that 
$\Delta_h  f=C$. To obtain the inequality for $\Delta_h^*$ just integrate over $h$.

{\it Proof of claim 2.}
If suffices to show that $\Delta_h^* p $ converges strongly in $L^1$ to $\Delta p$ if $\Delta p \in L^2(\Td)$. 
The Taylor formula gives for $dx$ almost every $x$: 
\be
\Delta_h^*p(x) &=& \frac{1}{c_d h^{d}}\frac{l_d}{h^2}\left[\int _{|y|\leq h}\int _{\sigma =0}^1(1-\sigma)y^t\cdot [D^2p(x+\sigma y)-D^2p(x)]\cdot y  \, d\sigma \, dy \right.\\
&+&  \left.\int _{|y|\leq h}\demi y^t\cdot D^2p(x)\cdot y \, dy \right]\\
&=& \Sigma_1(x) +\Sigma_2(x).
\en
Then $\Sigma_2(x)$ is equal $dx$ almost everywhere to $\Delta p (x)$ and $\Sigma_1$ converges strongly to 0 in $L^2$ since
$p\in W^{2,2}$ from elliptic regularity.
  
$\hfill \Box$

\bigskip
\noindent
Hence, since $\Delta p(t_i) = \rho(t_i)-1 \in L^2(\Td)$ (from Proposition \ref{4MN}), up to extraction of a subsequence in $h$, for almost every $x\in \Td$, $\Delta_h^* p(t_i,x)$ converges to $\rho(t_i,x)-1$.
Applying $\Delta_h^*$ to (\ref{4Deltati1}, \ref{4Deltati2}) at a point where equality
(\ref{4Deltati1}) holds 
we get that 
\beq
\Delta_h^* \bfi_{i,i+1} +\Delta_h^* \bfi_{i,i-1} 
\leq 2d -\frac{T^2}{N^2} \Delta_h^* p(t_i).\label{4ineqDelta}
\enq
Therefore we have  $d\rho(t_i)$ a.e. 
$$\limsup_{h\to 0} \Delta_h^* \bfi_{i,i+1} +\Delta_h^* \bfi_{i,i-1} \leq 2d + \frac{T^2}{N^2}(1-\rho(t_i)).$$ 
A convex function, say $f$, such that $\Delta_h^*f(x)$ is bounded as $h$ goes to 0 is $C^{1,1}$ at $x$, since  the trace controls the norm of a positive matrix. 
Moreover the left hand side is nonnegative, implying that  $d\rho(t_i)$ a.e. 
$$\rho(t_i,x)\leq 1+ \frac{2dN^2}{T^2}.$$
This yields a $\Linf$ bound on $\rho$ (that depends on $N$), and implies from elliptic regularity that $p(t_i)$ is $C^{1,\alpha}$ for every $0<\alpha<1$.
Since convex functions are continuous and  $p(t_i)$ is continuous, equality (\ref{4Deltati1}) holds now on a closed set of full measure for $\rho(t_i)$ and thus at every density point of $\rho(t_i)$.
Now we can apply point 1 of Lemma \ref{4bolemme} to obtain that $\Delta_h^* \bfi_{i,i+1} +\Delta_h^* \bfi_{i,i-1} $ is bounded at every density point of $\rho(t_i)$, implying the expected $C^{1,1}$ bound.
Convex functions are twice differentiable almost everywhere, and thus almost everywhere
$\Delta_h^*\bfi_{i,i+1} \to \Delta\bfi_{i,i+1}$ and $\Delta_h^*\bfi_{i,i-1} \to \Delta\bfi_{i,i-1}$.
We define  ${\cal E}_i$ to be the set of all points $x$ such that 
\begin{enumerate}
\item $x$ is a point of Lebesgue differentiability for $\rho(t_i)$ where $\rho(t_i,x)>0$,

\item $\Delta_h^*p(t_i,x) \to \rho(t_i,x)-1$ as $h\to 0$, 

\item $\bfi_{i,i+1}$ and $\bfi_{i,i-1}$ are twice differentiable, 

\item equality (\ref{4Deltati1}) holds.

\end{enumerate}
The set ${\cal E}_i$ is of full measure for $\rho(t_i)$ and by removing sets of 0 measure to ${\cal E}_i$ one can impose the condition
$\nabla\bfi_{i,i+1}({\cal E}_i)={\cal E}_{i+1}$.
Then at every point of ${\cal E}_i$ we have
\beq
&&0\leq \Delta\bfi_{i,i+1}(x)+\Delta\bfi_{i,i-1}(x)\leq 2d +\frac{T^2}{N^2}(1-\rho(t_i,x))\label{4Deltabfi},\\
&&-\frac{2dN}{T}\leq \Delta\phi(t_i^+,x)-\Delta\phi(t_i^-,x)\leq \frac{T}{N}(1-\rho(t_i,x)).
\enq
This implies the following  bound for $\rho,\phi$:
\be
&&\|\rho(t_i)\|_{\Linf}\leq C(d)N^2,\\
&& \textrm {for every } x\in{\cal E}_i,\, \|\phi(t_i,x)\|_{C_x^{1,1}}\leq C(d)N/T.
\en
The proof of Lemma \ref{4d2} is complete.

$\hfill \Box$

\paragraph{Construction of the characteristics}
Remember that $\bfi_{i,i+1}$ is given by:
\be
\bfi_{i,i+1}(x)= |x|^2/2  +\frac{T}{N} \phi(t_i^+,x),
\en
and that for $s \in [t_i, t_{i+1}]$ (resp. $s \in [t_{i-1}, t_{i}]$) ,   $\bfi_{t_i,s}$  is given by
\be
\bfi_{t_i,s}(x)&=& |x|^2/2  +(s-t_i) \phi(t_i^+,x),\\
(\textrm{ resp. }\bfi_{t_i,s}(x)&=& |x|^2/2  +(s-t_i) \phi(t_i^-,x)). 
\en
We know that in ${\cal E}_i$,  $\bfi_{i,i+1}, \bfi_{i,i-1}$ are convex and twice differentiable.
Thus for any $s \in ]t_{i-1}, t_{i+1}[$, $\bfi_{t_i,s}$ is twice differentiable and $D^2\bfi_{t_i,s}$ is invertible. This implies (see the Appendix on convex functions of \cite{Mc1}) that $D^2\bfi_{s,t_i}$ exists at point $\nabla\bfi_{t_i,s}(x)$ for $x\in{\cal E}_i$.
Thus we can define ${\cal E}_s= \nabla\bfi_{t_i,s}({\cal E}_i)$, this definition makes sense pointwise in
${\cal E}_i$, and pointwise in ${\cal E}_s$, $\phi_s$ is twice differentiable. Note also that
from Lemma \ref{4=**} the definition ${\cal E}_s= \nabla\bfi_{t_i,s}({\cal E}_i)$ and ${\cal E}_s= \nabla\bfi_{t_{i+1},s}({\cal E}_{i+1})$ are consistent since $\nabla\bfi_{i,i+1}({\cal E}_i)={\cal E}_{i+1}$.
Then we can define a trajectory $(x_s, s\in [0,T])$ as follows: starting from $X_{s_0}\in{\cal E}_{s_0}$
for any  $s_0\in ]0,T[$, we define
$x_s=\nabla\bfi_{s_0, s}(x_0)$ for any $s$ in the same interval $[t_i, t_{i+1}]$ as $s_0$, and we proceed similarly in other intervals. Thus we define a flow $\Xi(s,t,x)$ that gives at time $s$ the position of the particle located in $x$ at time $t$. This flow $\Xi(s,t,x)$ is defined everywhere on ${\cal E}_t$, and $\Xi(s,t,{\cal E}_t)={\cal E}_s$.
We may denote $(x_s,\, s\in[0,T])$ a trajectory and it will be understood that $x_t\in {\cal E}_t,  \ \forall t\in ]0,T[$.

\paragraph{Conclusion of the proof}
Now we bound $\Delta\phi$ along a trajectory: for $(t,s)\in [t_i,t_{i+1}]$, we have
\be
&&(t-s)\left( \Delta\phi(t)(x_t)-\Delta\phi(s)(x_s) \right)\\
&=&(t-s)\left( \Delta\phi(t,\nabla\bfi_{s,t}(x_s))-\Delta\phi(s,x_s) \right) \\
&=&2d-\Delta\bfi_{s,t}(x)-\Delta\bfi_{t,s}(\nabla\bfi_{s,t}(x_s)),
\en 
but this is negative since we have the relation 
$$D^2\bfi_{s,t}(x)=[D^2\bfi_{t,s}]^{-1}(\nabla\bfi_{s,t}(x)),$$
and therefore
\be
2d-\Delta\bfi_{s,t}(x)-\Delta\bfi_{t,s}(\nabla\bfi_{s,t}(x_s))=2d -\sum (\lambda_i+1/\lambda_i)\leq 0,
\en
where the $\lambda_i$ are the eigenvalues of $D^2\bfi_{s,t}(x_s)$ well defined on ${\cal E}_s$. 
Thus we conclude first  that for every $x_s\in {\cal E}_s$ and $x_t=\nabla\bfi_{s,t}(x_s)$  \beq\label{4Delta} 
(\Delta \phi(t,x_t)-\Delta \phi(s,x_s))\cdot(t-s)\leq 0 \ \textrm{ for } t_i\leq s,t \leq t_{i+1}.
\enq
Then we obtain a quantitative estimate of the decay of $\Delta \phi$, between $t_i$ and $t_{i+1}$: we take $s=t_i, t=t_{i+1}$ in the previous inequality, 
from the convexity of $ x\rightarrow x+1/x$ we have
\be
\sum_{j=1}^d(\lambda_j+1/\lambda_j)&\geq& d(\Delta/d + d/\Delta)
\textrm{ where } \Delta=\sum _{j=1}^d\lambda_j=\Delta\bfi_{i,i+1}\\
&=& \frac{1}{\Delta}(\Delta -d)^2 +2d\\
&=& \frac{(\frac{T}{N}\Delta \phi(t_i^+,x_i))^2}{\frac{T}{N}\Delta \phi(t_i^+,x_i)+d} +2d,
\en
since $\frac{T}{N}\Delta \phi(t_i^+)+d = \Delta\bfi_{i,i+1}$.
Hence 
\be
\Delta\phi(t_{i+1}^-,x_{i+1})-\Delta\phi(t_i^+,x_i) \leq -\frac{T}{N}\frac{(\Delta \phi(t_i^+,x_i))^2}{\frac{T}{N}\Delta \phi(t_i^+,x_i)+d}.
\en
Using (\ref{4Delta_i}) we obtain 
\be
\Delta\phi(t_{i+1}^+,x_{i+1})\leq \Delta\phi(t_{i}^+,x_{i})+\frac{T}{N}\left(1-\frac{(\Delta\phi(t_{i}^+,x_{i}))^2}{d+\frac{T}{N}\Delta\phi(t_{i}^+,x_{i})}\right).
\en
We know from (\ref{4Deltabfi}) that
$\Delta\bfi_{i,i+1}\leq 2d + T^2/N^2$ thus 
$\Delta\phi(t_{i}^+,x_{i})\leq \frac{N}{T}(d + T^2/N^2)$. It follows that
\mbox{$d+\frac{T}{N}\Delta\phi(t_{i}^+,x_{i}) \leq 2d + \frac{T^2}{N^2} \leq 3d$} for $N$ large enough.
We finally obtain the following bound for $N$ large enough:
\be
\Delta\phi(t_{i+1}^+,x_{i+1})\leq \Delta\phi(t_{i}^+,x_{i})+\frac{T}{N}\left(1-\frac{(\Delta\phi(t_{i}^+,x_{i}))^2}{3d}\right).
\en
This is a discrete version of the differential inequality
$\dot \Theta \leq 1- \frac{1}{3d}\Theta ^2$ and we will conclude that
\be
\Delta\phi(t_i) \leq C(d)(1+\frac{1}{t_i})  \text{ in } {\cal E}_i \text{ for any } 1\leq i \leq N-1.
\en
This will be a consequence of the following lemma:

\begin{lemme}
Let $(X_n)_{n\in [0..N]}$ be a sequence defined by $X_0\in \R $ and such that 
\be
X_{n+1}\leq X_n + \frac{T}{N}(1-X_n^2/C^2).
\en
Then for $n\geq 1$,
$$X_n \leq C+ \frac{C^2}{t_n}, \  \ t_n=\frac{Tn}{N}.$$
\end{lemme}

{\it Proof.} The proof is by induction.
For $n=1$, we look for 
\be
\max_{X_0} \{X_0 + \frac{T}{N}(1-X_0^2/C^2)\}.
\en
This yields $$X_i \leq T/N + NC^2/(4T) \leq C+ C^2/t_n=C+C^2N/(nT)$$ for $N \geq T/C$, $1\leq n \leq 4$.
Then suppose that $n\geq 4$, $X_n\leq C+C^2/t_n$. The function \linebreak \mbox{$x\to x+\frac{T}{N}(1-x^2/C^2)$} is increasing for $x\leq NC^2/(2T)$. For $N\geq 4T/C$ and $n\geq 4$, we have $$NC^2/(2T) \geq C + C^2/t_n= C + NC^2/(nT).$$
Therefore we have (for $N\geq 4T/C, n\geq 4$) 
\be
X_{n+1} &\leq& C+C^2/t_n + \frac{T}{N}(1-(C+C^2/t_n )^2/C^2)\\
        &=& C + C^2/t_n \, (1-2T/(CN) - T/(Nt_n)).
\en
In order to conclude we need $1-2T/(CN) - T/(Nt_n) \leq n/(n+1)$.
But we have \linebreak $1-T/(Nt_n) = 1 - 1/n \leq n/(n+1)$, therefore we conclude the proof of the lemma.

$\hfill \Box$

Using then the transformation $\phi(t,x) \rightarrow -\phi(-t,x)$, 
that transforms the solution of Problem \ref{4infIN} in another solution of \ref{4infIN}
exchanging $\rho_0$ and $\rho_T$, we get that
$\Delta\phi_N(t_i) \geq -C(d)(1+\frac{1}{T-t_i})$  in ${\cal E}_i$ for any $1\leq i \leq N-1$.
Since we know from (\ref{4Delta}) that $t \to \Delta\phi(t,x_t)$ is decreasing between $t_{i}$ and $t_{i+1}$ 
we can  conclude that there exists for each $t$ a set of full 
measure for $d \rho(t)$  on which $\phi(t,.)$ is twice differentiable, and where the following equality
holds:
\beq\label{4boundDeltaphi}
-C(d)(1+\frac{1}{T-t}) \leq \Delta\phi_N(t,x) \leq  C(d)(1+\frac{1}{t}).
\enq
The first part of Proposition \ref{4regulFN} is proved.


\subsubsection{Proof of the $\Linf_{loc}(]0,T[\times \Td)$ bound on $\rho$}
We begin by writing the Monge-Amp\`ere equation that links $\rho(t_i)$ to $\rho(t_{i+1})$
\be
&&\rho(t_{i+2}, x_{i+2})\det (I+\frac{T}{N} D^2\phi(t_{i+1}^+, x_{i+1}))=\rho(t_{i+1}, x_{i+1}),\\
&&\rho(t_i,x_{i})\det (I-\frac{T}{N} D^2\phi(t_{i+1}^-, x_{i+1}))=\rho(t_{i+1}, x_{i+1}).
\en
This equation makes sense since, from Lemma \ref{4d2}, at $x_s$, $\phi$ is twice differentiable, and we use 
\linebreak \cite[Proposition A2]{Mc1}.

Now using the domination of the geometric mean by the arithmetic mean we have
\be
\det (I+\frac{T}{N} D^2\phi(t_{i+1}^+))\leq (1 + \frac{T}{dN} \Delta\phi(t_{i+1}^+))^d,
\en
hence
\beq
&&\frac{\rho(t_{i+1}, x_{i+1})}{\rho(t_{i+2}, x_{i+2})}\leq (1+\frac{T}{dN} \Delta\phi(t_{i+1}^+,x_{i+1}))^d\label{4dtrho1},\\
&&\frac{\rho(t_{i+1}, x_{i+1})}{\rho(t_i,x_i)}\leq (1-\frac{T}{dN} \Delta\phi(t_{i+1}^-,x_{i+1}))^d\label{4dtrho2}.
\enq
We deduce first the following:
\be
\frac{1}{(1+\frac{T}{dN} \Delta\phi(t_{i}^+,x_{i}))^d}\leq      \frac{\rho(t_{i+1}, x_{i+1})}{\rho(t_i,x_i)}\leq (1-\frac{T}{dN} \Delta\phi(t_{i+1}^-,x_{i+1}))^d.
\en
Note that we also have
\beq\label{4liprho}
\frac{1}{(1+\frac{t-s}{d} \Delta\phi(s,x_s))^d}\leq      \frac{\rho(t,x_t)}{\rho(s,x_s)}\leq (1-\frac{t-s}{d} \Delta\phi(t,x_t))^d
\enq
for $t_i < s,t < t_{i+1}$, $x_s \in {\cal E}_s$, and this implies using (\ref{4boundDeltaphi}) that  along a trajectory $s\to x_s$, $\log(\rho(s,x_s))$ is Lipschitz: for all $t_1,t_2\in[\tau,T-\tau]$,

\beq\label{4trajlip}
|\log(\rho(t_1,x_{t_1})-\log(\rho(t_2,x_{t_2})|\leq C(1+\frac{1}{\tau(T-\tau)})|t_2-t_1|.
\enq
Taking the logarithm of (\ref{4dtrho1}, \ref{4dtrho2}) we obtain
\beq
&&\log(\rho(t_{i+2}, x_{i+2}))+\log(\rho(t_i,x_i))-2\log(\rho(t_{i+1}, x_{i+1}))\nonumber\\
&\geq &-d\log (1+\frac{T}{dN} \Delta\phi(t_{i+1}^+,x_{i+1})) - d\log (1-\frac{T}{dN} \Delta\phi(t_{i+1}^-,x_{i+1}))\nonumber\\
&\geq&-\frac{T}{N} (\Delta\phi(t_{i+1}^+,x_{i+1})-\Delta\phi(t_{i+1}^-,x_{i+1}))
\nonumber\\
&\geq& \frac{T^2}{N^2}(\rho(t_{i+1}, x_{i+1})-1)\label{4soiqrho},
\enq
where at the third line we have used the concavity of the $\log$ and at the last line we have used (\ref{4Delta_i}):
\be
\Delta\phi(t_{i+1}^+,x)-\Delta\phi(t_{i+1}^-,x)\leq \frac{T}{N}(1-\rho_{i+1}(x)) \ d\rho(t_{i+1}) \, a.e..
\en
\\
We fix $\tau \in ]0, T/2[$. For any trajectory $x_s, s\in]0,T[$ with $x_s\in {\cal E}_s$ for all $s$, 
 $\log(\rho(s,x_s))$ is uniformly Lipschitz with respect to $s$  in $[\tau, T-\tau]$ from (\ref{4trajlip}) and $\log(\rho(s,x_s))$ remains finite in $[\tau, T-\tau]$.
Moreover (\ref{4soiqrho}) holds at every time $t_i$.
Using this we claim  an unconditional bound for $\rho(s,x_s)$ for $\tau\leq s \leq T-\tau$.
\\
Proof of claim:
The sequence $(\Theta(t_i)=\log \rho(t_i,x_i))_{1\leq i \leq N-1}$ satisfies a discretization of  the differential inequality 
\beq\label{4ineqdif}
\ddot \Theta \geq \exp \Theta -1. 
\enq
From (\ref{4trajlip}), we have the {\it a priori} bound $$\left|\Theta(t_{i+1})-\Theta(t_i)\right|\leq C(\tau) T/N$$ for $t_i$ in $[\tau, T-\tau]$.
We argue by contradiction: 
take $0<\tau <T/4$ and suppose $\Theta(t_{i_0})\geq M$, with $t_{i_0} \in [\tau, T/2]$. Then choose $M$ so large that $M-C(\tau)(T-2\tau)\geq M/2$.  Thus, on 
$[\tau, T-\tau]$ we have $\Theta \geq M/2$ from the {\it a priori} bound above. This implies that 
$$[\Theta(t_{i+1})- \Theta(t_i)]  -  [\Theta(t_{i})- \Theta(t_{i-1})]   \geq  T^2/N^2 (\exp(M/2)-1)$$ on $[\tau, T-\tau]$. We then choose $M$ large enough so that $\Theta_{[(T-\tau) N]/T}- \Theta_{[(T-\tau) N]/T -1}  > C(\tau)T/N$ in 
contradiction with the {\it a priori} bound above. 

Applying the same argument after having changed $t$ in $T-t$ gives the bound on $[T/2, T-\tau]$.
Note that this proof does not depend on the initial and final 
values of $\Theta$.
 
\subsubsection{Time continuity of $\rho$} Remember that in Theorem \ref{4main} we have proved that $\rho \in C(]0,T[; L^p(\Td))$ for any $p\in [1,\frac{3}{2}[$. We have now  an unconditional bound on $\rho$ in
 $\Linf_{loc}(]0,T[\times \Td)$. Thus the strong time  continuity in every $L^p(\Td)$,  $1\leq p <\infty$  follows and the point 2 of Proposition \ref{4regulFN} is proved.

\subsubsection{Lipschitz bound for $\log\left(\|\rho(t,\cdot)\|_{L^k(\Td)}\right)$}
Since $\phi$ is twice differentiable at $x_s\in{\cal E}_s$, we can use the identity
\be
\rho(t,x_t)\det (I+(t-s)D^2\phi(s,x_s))=\rho(s,x_s),
\en
and it implies
\be
\Dt \int_{\Td}[\rho(t,x)]^k \ dx = -(k-1)\int_{\Td}[\rho(t,x)]^k \Delta\phi(t,x) \ dx.
\en
We have thus 
\be
\Dt \|\rho(t,\cdot)\|_{L^k} \leq \frac{k-1}{k} \|\Delta\phi(t,\cdot)\|_{\Linf(d\rho(t))} \|\rho(t,\cdot)\|_{L^k}, 
\en
hence using point 1 of Proposition \ref{4regulFN} we get that 
\be
-C(d)(1+\frac{1}{t})\leq \Dt \log \left( \|\rho(t,\cdot)\|_{L^k(\Td)} \right) \leq C(d)(1+\frac{1}{T-t}).
\en
This proves the point 3 of Proposition \ref{4regulFN}.

\subsubsection{Displacement convexity of functionals of $\rho$} 
Here we show that given $\rho_N(t,x)$ solution of Problem $\ref{4infIN}$, the functions $\int_{\Td}\rho_N(t,x)\log(\rho_N(t,x)) \ dx,$  $\int_{\Td}[\rho_N(t,x)]^k \ dx$, $k\in [1, +\infty[$ converge to convex functions of $t\in [0,T]$.
We drop the subscript $N$.
Let  $\frac{d}{dt}$ denote the convective derivative $\dt .+ \nabla \phi \cdot \nabla .$ and $\frac{d^2}{dt^2}=(\frac{d}{dt})^2$. The density $\rho$ satisfies (\ref{4soiqrho}) which is the finite difference version  of 
\beq\label{4d2dt2logrho}
\frac{d^2}{dt^2}\log \rho(t,x) \geq \rho(t,x) -1.
\enq
Note that we do not include all terms, since one could show that our solutions satisfies
\be
\frac{d^2}{dt^2}\log \rho(t,x) \geq \frac{1}{d}\left|\Dt \rho(t,x)\right|^2 +  \rho(t,x) -1.
\en
The term we have omitted here allows to extend the displacement convexity of $(k-1)\int [\rho(t,x)]^k \ dx$ down to $k\geq 1-1/d$ in the non-interacting case (cf \cite{Mc1}). However in the gravitating case, the additional term will limit this property to positive values of $k-1$.
What could probably be proved in this case is that between $k=1$ and $k=1-1/d$ the functional $(k-1)\int [\rho(t,x)]^k \ dx$ is semi-convex.

\paragraph{Formal proof}
We first give a formal proof, eluding the fact that  we only have a discrete version of the differential inequalities.
The lemma given below will justify the calculations.
Using (\ref{4d2dt2logrho}) with the identities
\be
&&\frac{d^2}{dt^2}\rho(t,x)  =\frac{1}{\rho(t,x)}|\Dt \rho(t,x)|^2 +\rho(t,x)\frac{d^2}{dt^2}\log \rho(t,x),\\
&&\frac{d^2}{dt^2}[\rho(t,x)]^k=k(k-1)[\rho(t,x)]^{k-2}|\Dt \rho(t,x)|^2  + k [\rho(t,x)]^{k-1} \frac{d^2}{dt^2}\rho(t,x), \en
we obtain for $k\geq 0$ 
\be
\frac{d^2}{dt^2}[\rho(t,x)]^k \geq k^2[\rho(t,x)]^{k-2}|\Dt \rho(t,x)|^2 + k[\rho(t,x)]^k(\rho(t,x)-1).
\en
Noticing that for all smooth $F$ we have
\be
\frac{d^2}{dt^2}\left[\int_{\Td}\rho(t,x) F(\rho(t,x)) \ dx\right] =\int_{\Td}\rho(t,x)\frac{d^2}{dt^2}(F(\rho(t,x))) \ dx,
\en
and applying this to $F(\rho(t,x))=[\rho(t,x)]^k,\, k\geq 0$,
we get
\be
\frac{d^2}{dt^2}\int_{\Td} [\rho(t,x)]^{k+1} \ dx\geq\int_{\Td} k ([\rho(t,x)]^{k+2}-[\rho(t,x)]^{k+1}) \ dx \geq 0.
\en
Indeed, using Jensen's inequality twice we have 
\be
\int_{\Td} [\rho(t,x)]^{k+1} \ dx &\geq& \left(\int_{\Td} [\rho(t,x)]^k \ dx\right)^{\frac{k+1}{k}}\\
&=&  \left(\int_{\Td} [\rho(t,x)]^k \ dx\right)^{\frac{1}{k}} \int_{\Td} [\rho(t,x)]^k \ dx \\
&\geq& 1 \int_{\Td} [\rho(t,x)]^k \ dx
\en for $k\geq 1$ to conclude.

This convexity property combined with the unconditional bound for $\rho$ in $\Linf([\tau, T-\tau]\times\Td)$ yields a uniform Lipschitz bound for 
$\|\rho(t,\cdot)\|_{L^k}$ in $[\tau, T-\tau]$ for any $1\leq k \leq \infty$. (For the case $k=+\infty$ this is because of (\ref{4trajlip})).

\paragraph{Rigorous proof}
\begin{lemme}
Let $k$ be greater than 1. Then $\bar\rho$, the limit of the sequence $\rho_n$ solutions of Problem \ref{4infIN},
satisfies 
\be
\frac{d^2}{dt^2} \int [\bar\rho(t,x)]^k \ dx \geq (k-1)\int [\bar\rho(t,x)]^{k+1} -[\bar\rho(t,x)]^k \ dx
\en
in the sense of ${\cal D'}(0,T)$.
\end{lemme}

{\it Proof.} We consider $\rho$ solution of Problem \ref{4infIN}. We find after some elementary calculations that
\be
&& \rho_{i+1}^{k}(x_{i+1})+ \rho_{i-1}^{k}(x_{i-1})-2 \rho_{i}^{k}(x_{i})\\
&\geq&k\rho_{i}^{k-1}(x_{i})(\rho_{i+1}(x_{i+1})+ \rho_{i-1}(x_{i-1})-2 \rho_{i}(x_{i}))\\
&+& \demi k(k-1) \rho_i^{k-2}(x_i)\left((\rho_{i+1}(x_{i+1})- \rho_{i}(x_{i}))^2+ (\rho_{i-1}(x_{i-1})- \rho_{i}(x_{i}))^2\right)\\
&+& O(N^{-3}).
\en
The term $O(N^{-3})$ depends on $\|\Delta\phi(t)\|_{\Linf(d\rho(t))}, \|\rho(t)\|_{\Linf}$ and is therefore uniformly bounded in compact sets of $]0,T[$.
Then we use (\ref{4soiqrho}), 
to get also by simple calculations (using the concavity of the logarithm)
\be
\rho_{i+1}(x_{i+1})+ \rho_{i-1}(x_{i-1})-2 \rho_{i}(x_{i}) &\geq& \frac{T^2}{N^2}\rho_i(x_i)(\rho_i(x_i)-1)\\
&+&\demi\left(\frac{(\rho_{i+1}(x_{i+1})- \rho_{i}(x_{i}))^2}{\rho_i(x_i)}+ \frac{(\rho_{i-1}(x_{i-1})- \rho_{i}(x_{i}))^2}{\rho_i(x_i)}\right).
\en
Combining these inequalities, we get
\be
&& \rho_{i+1}^{k}(x_{i+1})+ \rho_{i-1}^{k}(x_{i-1})-2 \rho_{i}^{k}(x_{i})
\geq \frac{T^2}{N^2}k\rho_{i}^{k}(x_{i}))(\rho_i(x_i)-1) + k^2 Q + O(N^{-3}),
\en
with $Q$ a positive quantity. The term $O(N^{-3})$ is uniform on every interval $[\tau, T-\tau]$, for $\tau > 0$.
We  integrate over $\Td$, this yields
\be
\int \left(\rho_{i+1}^{k+1} + \rho_{i-1}^{k+1}-2 \rho_{i}^{k+1}\right) \ dx
\geq  \frac{T^2}{N^2}\int k\left(\rho_{i}^{k+2} - \rho_i^{k+1}\right) \ dx + O(N^{-3}).
\en
We now choose a test function $\varphi\in C^{\infty}_c(]0,T[)$.
Denoting $D^2_h$ the second incremental quotient \mbox{($D^2_h f(x)=h^{-2} (f(x+h)+f(x-h)-2f(x))$}, with $h=T/N$), we have
\be
\frac{T}{N}\sum_{i=1}^{N-1}\left(D^2_h \int [\rho(t_i,x)]^k\right) \varphi(t_i)& =& \frac{T}{N}\sum_{i=1}^{N-1}\int [\rho(t_i,x)]^k D^2_h\varphi(t_i) \\
&\to& \int_D [\bar \rho(t,x)]^k \varphi''(t) \ dtdx
\en
which proves the lemma.

$\hfill \Box$

The convexity of $\int \rho\log\rho$ follows the same lines, and we skip the proof.

\subsubsection{Proof of the $W^{1,\infty}([\tau,T-\tau]\times\Td)$ bound for $\phi$}

Hereafter we use again the subscript $N$ for the solution of Problem \ref{4infIN} while $\nabla\phi,p$ is the solution of Problem \ref{4infI2}.
\\
From Lemma \ref{4=**}, $\phi_N$ can be given by the Hopf-Lax formula (\ref{4defvisco}) that we recall here:
\beq\label{4defvisco2}
\phi_N(t,x)=
\inf_{\gamma(t)=x}\left\{\phi_N(s,\gamma(s))+\int_s^t  \left[-  \frac{T}{N}\sum_{i=1}^{N-1}\delta_{\sigma=t_i}p_N(\sigma, \gamma(\sigma))+ \demi|\dot\gamma|^2(\sigma) \right]\,d \sigma\right\},
\enq
where we still use $t_i=\frac{Ti}{N}$.
This formula is valid for any $0\leq s\leq t\leq T$.
We are going to prove the following lemma that will yield the result when letting $N$ go to $+\infty$.
\begin{lemme}\label{4lemmelipvisco}
Let $\phi_N$ be defined by (\ref{4defvisco2}).
If $|\nabla p_N(t_i,\cdot)|\leq l(\tau)$ for any $t_i\in [\tau, T-\tau]$
with $l(\tau)<+\infty$ for $\tau>0$, and $p_N(t_i,\cdot)$ has mean value 0, 
 then 
\be
&&|\nabla\phi_N(t,x)|\leq C(\tau),\\
&&|\phi_N(t,x)-\phi_N(s,x)| \leq C(\tau) \left(\frac{T}{N}(1+{\bf E}(\frac{N(t-s)}{T})) + |t-s|\right),
\en
for any $0<\tau\leq t,s \leq T-\tau$, $dx$ a.e. $x\in \Td$,  where $\bf{E}(\cdot)$ denotes the integer part.
\end{lemme}

{\it Proof.} Let $\gamma$ be a minimizer in the infimum (\ref{4defvisco2}). Let $\tilde\gamma(\sigma) = \gamma (\sigma) + \frac{\sigma-s}{t-s}(z-x)$. Then 
\be
\phi_N(t,z)&\leq& \phi_N(s,\tilde\gamma(s))+\int_{s}^t  \left[-  \frac{T}{N}\sum_{i=1}^{N-1}\delta_{\sigma=t_i}p_N(\sigma, \tilde\gamma(\sigma))+ \demi|\dot{\tilde{\gamma}}|^2(\sigma) \right]\,d \sigma\\
&\leq& \phi_N(s,\gamma(s)) + \int_{s}^t  \left[-  \frac{T}{N}\sum_{i=1}^{N-1}\delta_{\sigma=t_i}p_N(\sigma, \gamma(\sigma))+ \demi|\dot{\gamma}|^2(\sigma) \right]\,d \sigma\\
&+& \int_{s}^t |\dot {\gamma}(\sigma)|\frac{|z-x|}{|t-s|} + \frac{|z-x|^2}{2|t-s|^2} \ d\sigma + (1+\frac{T}{N}{\bf E}(\frac{N(t-s)}{T}))l(\tau) |x-z|. 
\en
We choose here $s,t$ such that $t-s\geq \tau/2$, $s\geq \tau/2$ and $t\leq T-\tau$. We will have $\dot \gamma$  uniformly bounded in $L^2([s,t])$ by a constant $C(\tau)$. 
The second line is equal  to $\phi_N(t,x)$, the third line is bounded by
$C(\tau)\frac{|z-x|}{t-s}$
and we get that 
\be
\|\nabla\phi_N\|_{\Linf([\tau,T-\tau]\times\Td)}\leq C(\tau).
\en

Now we choose $s,t \in [\tau, T-\tau]$. By taking $\gamma(\sigma)\equiv x$ in (\ref{4defvisco2}), we obtain that \linebreak $\phi_N(t,x)\leq \phi_N(s,x) + C(\tau)\frac{T}{N}(1+{\bf E}(\frac{N(t-s)}{T}))$.
Using now that $\nabla\phi_N$ is bounded by $C(\tau)$ we get
\be
&&\phi_N(t,x) - \phi_N(s,x)\\
 &\geq& \inf_{y\in \Td}\left\{\phi_N(s,y)-\phi_N(s,x) - \frac{T}{N}(1+{\bf E}(\frac{N(t-s)}{T}))\|p_N\|_{\Linf([s,t]\times\Td)} + \frac{|y-x|^2}{2(t-s)} \right\}\\
&\geq& \inf_{y\in \Td}\left\{- \frac{T}{N}(1+{\bf E}(\frac{N(t-s)}{T}))\|p_N\|_{\Linf([s,t]\times\Td)} - C(\tau)|y-x| + \frac{|y-x|^2}{2(t-s)} \right\}\\
&\geq& -C'(\tau)\left(\frac{T}{N}(1+{\bf E}(\frac{N(t-s)}{T}))+|t-s|\right).
\en
This proves the lemma. $\hfill \Box$

We can use this lemma to conclude the $W^{1,\infty}([\tau,T-\tau]\times\Td)$ bound for $\phi$. Indeed, we know already that
$|\rho_N(t,x)|\leq C(\tau)$ for $t\in [\tau, T-\tau]$. Since $\Delta p_N = \rho_N-1$, we thus have $p_N(t)\in W^{2,k}(\Td)$ for $1\leq k <\infty$, 
uniformly on $[\tau, T-\tau]$.  The Sobolev embedding Theorem then yields  $$\|p_N\|_{\Linf([\tau, T-\tau]; C^{1}(\Td))} \leq C(\tau).$$ 
We let then $N$ go to $+\infty$ so that 
$\frac{T}{N}E(\frac{N(t-s)}{T})\to t-s$. We already know that $\rho_N$ converges to $\rho$ in $C(]0,T[; L^p(\Td))$ for $p<\infty$.
From the bounds obtained we will also have $p_N$ converging uniformly to $p$ on compact sets of $]0,T[\times \Td$.
Finally, extracting a subsequence if necessary, $\phi_N$ will converge to $\phi$ uniformly on compact sets of $]0,T[\times \Td$.

\subsubsection{Convergence to viscosity solutions}
For a given smooth $\gamma: [s,t]\to \Td$ compute
\be
&&\phi_N^\gamma(t,x) = \phi_N(s,\gamma(s))+\int_s^t  \left[-  \frac{T}{N}\sum_{i=1}^{N-1}\delta_{\sigma=t_i}p_N(\sigma, \gamma(\sigma))+ \demi|\dot\gamma|^2(\sigma) \right]\,d \sigma,\\
&&\phi^\gamma(t,x) = \phi(s,\gamma(s))+\int_s^t  \left[- p(\sigma, \gamma(\sigma))+ \demi|\dot\gamma|^2(\sigma) \right]\,d \sigma.
\en
Combining the mass conservation equation (\ref{4continuite})  and the Poisson equation, we get
\be
\dt p_N = \Delta^{-1}(-\nabla\cdot(\rho_N \nabla\phi_N)).
\en
From the $\Linf_{loc}(]0,T[\times\Td)$ bounds on $\rho_N, \nabla\phi_N$ and elliptic regularity, we get that for any $d<k<\infty$ 
\be
\dt p_N(t,.) \in W^{1,k}(\Td) \subset C^\alpha(\Td), \  \alpha = 1-\frac{d}{k},
\en uniformly for $t\in [\tau, T-\tau]$. We also already know that $p_N(t)$ is uniformly Lipschitz in space for $t\in [\tau, T-\tau]$. 
Hence we have $\phi_N^\gamma(t,x) \to \phi^\gamma(t,x)$ since $(\phi_N, p_N)$ converge uniformly to $(\phi,p)$ in every compact set of $]0,T[\times\Td$ and from the bound just obtained on $p$.  Since we can choose $\dot \gamma$ to remain bounded in $L^2([s,t])$ and thus $\gamma$ bounded
in $C^{\demi}([s,t])$ we  conclude that 
\be
\inf_{\|\dot\gamma\|_{L^2[s,t]}\leq C, \gamma(t)=x}\{\phi_N^\gamma(t,x)\} 
- \inf_{\|\dot\gamma\|_{L^2[s,t]}\leq C, \gamma(t)=x}\{\phi^\gamma(t,x)\} \to 0
\en
as $N \to \infty$. Hence  we get that
\beq\label{4defvisco2bis}
\phi(t,x)=\inf_{\gamma(t)=x} \left\{\phi(s,\gamma(s))+\int_s^t  \left[- p(\sigma, \gamma(\sigma))+ \demi|\dot\gamma|^2(\sigma) \right]\,d \sigma\right\}
\enq
and $\phi$ is the viscosity solution of $\dt\phi+\demi|\nabla\phi|^2 + p=0$ on every $[s,t]\subset]0,T[$.

{\it Remark.} In particular, $\phi$ is a subsolution, therefore admissible for the time continuous dual problem. Hence, if $\phi_N$ maximizes the time discretized dual problem, it converges to the maximizer of the time continuous dual problem, uniformly on compact sets of $]0,T[\times\Td$.

\subsubsection{Reversibility and $C^{1, L\log L}$ regularity for $\phi$}
We can define the backward solution $\psi$ starting from $\phi(T)$ by $\psi(T)=\phi(T)$ and 
\beq\label{4defvisco3}
\psi(s,x)= \sup_{\gamma(s)=x}\left\{\psi(t,\gamma(t)) - \int_s^t  \left[- p(\sigma, \gamma(\sigma))+ \demi|\dot\gamma|^2(\sigma) \right]\,d \sigma\right\}
\enq
for any $0\leq s \leq t\leq T$. 
This amounts to take for $\psi$ the limit of the sequence $-\tilde\psi_N(T-t)$ in Lemma \ref{4=**}. Therefore $t\to-\psi(T-t,\cdot)$ is the viscosity solution of 
$\dt\varphi+ \demi|\nabla\varphi|^2+q=0$ with $q(t,\cdot)=p(T-t,\cdot)$.
Then from Lemma \ref{4=**}, we know that $\psi=\phi$ $d\rho$ a.e..
Since $\psi$ and $\phi$ are Lipschitz continuous, the set $\phi=\psi$ is a closed set of full measure for $\rho(t)$ for each $t\in ]0,T[$. This means that on the set that contains the dynamics ({\it i.e.} the support of $\rho$), the solution of our Hamilton-Jacobi equation is reversible, a property which is generically false for viscosity solutions. This property 
is close to several regularity properties obtained in \cite{EG} on the Aubry set in weak KAM theory.
This implies regularity on the support of $\rho$:

Using that in the time interval $[s,t]$, $\rho=1+\Delta p$ is uniformly bounded, we have from elliptic regularity (see \cite{GT} for example) $\nabla p$ uniformly log-Lipschitz, (that we will denote by $p\in C^{1, L\log L}$) {\it i.e.} $\forall (x,y)\in \Td \text{ with } |x-y| \leq \demi, \forall t\in [\tau, T-\tau]$, 
$$ \left|\nabla p(t,x)-\nabla p(t,y)\right| \leq C(\tau) |x-y| \log(\frac{1}{|x-y|}).$$
Then consider an optimal path $\gamma$ in (\ref{4defvisco2bis}). Perturb it in $\tilde\gamma(\sigma) = \gamma(\sigma) + \frac{\sigma-s}{t-s}(z-x)$.
We have 
\be
\phi(t,z)\leq T^a_{t,x}(z)=\phi(t,x) &+& \int_s^t \dot\gamma(\sigma)\cdot \frac{z-x}{t-s} + \demi \frac{|z-x|^2}{|t-s|^2} \ d\sigma  \\ &-&\int_s^t \left[ p(\sigma, \gamma+ \frac{\sigma-s}{t-s}(z-x)) - p(\sigma, \gamma)\right] \ d\sigma.
\en   
The second integral is $C^{1, L\log L}$ with respect to $z$, while the first is a quadratic polynomial in $z$. Moreover we have $T^a_{t,x}(z) \geq \phi(t,z)$ and  
$T^a_{t,x}(x) = \phi(t,x)$, which means that $T^a_{t,x}$ is tangent from above to 
$\phi(t,\cdot)$ at $x$.

In a similar way, we use  (\ref{4defvisco3}) to obtain at every point $(x,t)$ a $C^{1, L\log L}$ function $ T^b_{x,t}$ tangent from below to $\psi$.
The log-Lipschitz constants are uniformly bounded in $[\tau, T-\tau]$.

If follows that at every point of the set $\{\psi=\phi\}$, $\phi$ (or $\psi$) is pinched between two functions whose gradients are log-Lipschitz.
 This implies the following, which is the point 7 of Theorem \ref{4main2}:
\begin{prop}
For every $t\in [\tau, T-\tau], \ \tau > 0$, there exists a closed set ${\cal S}_t$ of full measure for $\rho(t)$ such that
$\phi$ is differentiable at every point of ${\cal S}_t$, and
\be
\forall \, (x,y) \in {\cal S}_t, |x-y|\leq 1/2, \  \  \   \left| \nabla\phi(t,x)-\nabla\phi(t,y) \right| \leq C(\tau) |x-y| 
\log(\frac{1}{|x-y|}).
\en
\end{prop}

{\it Proof.} We define the set ${\cal S}_t=\{x\in \Td, \phi(t,x)=\psi(t,x)\}$. We choose two points $(x,y)\in {\cal S}_t$. 
Using that $\phi$ is pinched at $x$ between two $C^{1, L\log L}$ functions, 
 we have at $x$, for all $z\in \Td$,
\be
 \left|\phi(t,z) - \phi(t,x) -\nabla\phi(t,x)\cdot (z-x)\right| \leq C|x-z|^2 \log(\frac{1}{|x-z|}),
\en
and we have the same changing $x$ in $y$. 
Combining the two inequalities, we get that 
\be
&& \left| (\nabla\phi(t,y)-\nabla\phi(t,x))\cdot(z-x)\right|\\
&\leq& C
\left( |x-z|^2 \log(\frac{1}{|x-z|}) + |y-z|^2 \log(\frac{1}{|y-z|}) + |y-x|^2 \log(\frac{1}{|y-x|})    \right).
\en
Taking the maximum value of the first line among all $z$ such that $\{|z-x|=|y-x|\}$ (note that this implies
$|y-z|\leq 2 |y-x|$), we get that
\be
\left|\nabla\phi(t,y)-\nabla\phi(t,x)\right|\leq C'|x-y| \log(\frac{1}{|x-y|}).
\en
This proves the proposition.

$\hfill \Box$

\subsubsection{Estimates up to the boundary}
If $\rho_T\in L^k(\Td)$ with $k>d$, then $p$ remains Lipschitz in space up to $t=T$, hence
one can choose $\phi \in W^{1,\infty}([\tau,T]\times\Td)$. 
Instead of (\ref{4defvisco2bis}) one can have $\phi$ defined by (\ref{4defvisco3}),
hence if $\rho_0\in L^k(\Td)$, one can
choose $\phi \in W^{1,\infty}([0,T-\tau]\times\Td)$.
If $\rho_0$ and $\rho_T$ are both in $L^k(\Td)$, with $k>d$ then, using (\ref{4defvisco3})
one can choose first $\phi$ such that 
$\phi(0,\cdot)\in W^{1,\infty}(\Td)$. Then we obtain that 
$\phi\in W^{1,\infty}([0,T]\times\Td)$ using the following result that can be found in \cite{E}:
\begin{prop}
Let $\phi$ be the viscosity solution on $[0,T]\times\Td$ of 
\be
&&\dt\phi+ \demi|\nabla\phi|^2+p=0,\\
&&\phi(t=0,\cdot)=\phi_0,
\en
with $p\in \Linf([0,T]; W^{1,\infty}(\Td))$ and $\phi_0\in W^{1,\infty}(\Td)$. Then
$\phi \in W^{1,\infty}([0,T]\times\Td)$.
\end{prop}
This last result ends the proof of Theorem \ref{4main2}.

$\hfill \Box$

\bibliography{ep-variationnel-biblio}

\end{document}